\documentclass{article}

\usepackage{bm}
\usepackage{amsmath}
\usepackage{amsthm}
\usepackage{amscd}
\usepackage{amssymb}
\usepackage[top=2cm,bottom=2cm,left=1.5cm, right=2cm]{geometry}
\usepackage{tikz}
\usepackage{pgfplots}
\usetikzlibrary{calc,3d}
\usepackage{psfrag}
\usepackage{pgfpages}
\usepackage{epstopdf}
\usepackage{graphicx}
\usepackage{booktabs} 
\usepackage{caption}
\usepackage{subcaption} 
\usepackage{placeins}
\usepackage[makeroom]{cancel}
\usepackage{enumerate} 

\usepackage{tikz}

\usetikzlibrary{shapes.geometric}
\usepgflibrary{plotmarks}
\usetikzlibrary{patterns}
\usetikzlibrary{snakes}
\usetikzlibrary{calc,3d}
\usetikzlibrary{backgrounds}
\usetikzlibrary{automata}
\usetikzlibrary{fadings}
\usetikzlibrary{arrows}
\tikzstyle{nod}=[circle,solid,draw=opiswykr!80,fill=opiswykr!40,scale=.45]
\tikzstyle{nod2}=[circle,solid,draw=opiswykr!80,fill=opiswykr!80,scale=.5]

\newcommand{\delete}[1]{\textcolor{blue}{\st{#1}}}
\renewcommand{\delete}[1]{}

\newcommand{\Real}{\mathbb{R}}

\DeclareMathOperator*{\diverg}{div}

\DeclareMathOperator{\range}{\mathcal{R}}
\DeclareMathOperator{\rank}{rank}
\DeclareMathOperator{\Span}{span}

\DeclareMathOperator{\Tr}{tr}

\newtheorem{theorem}{Theorem}
\newtheorem{conjecture}{Conjecture}

\newcommand{\myeqref}[1]{\eqref{#1}}
\newcommand{\pfdg}{{PF-DG{}}}
\newcommand{\PFDG}{\pfdg}

\newlength{\figuredim}
\newcommand{\bilinearoper}{a}
\newcommand{\bilinearoperoper}{A}
\newcommand{\linearoper}{l}
\newcommand{\vspan}[1]{\Span\{#1\}}

\newcommand{\oinnerprod}[2]{\left( #1 \right)_{#2} }
\newcommand{\binnerprod}[2]{\left\langle #1 \right\rangle_{#2} }
\newcommand{\boundarysign}[1]{{#1}_b}

\newcommand{\thematrix}{\vm {A}}

\newcommand{\rvector}{\vm b}

\newcommand{\matrows}{n}
\newcommand{\matcols}{m}
\newcommand{\solsset}{S}
\newcommand{\transfmatrix}{\vm W}
\newcommand{\transfvector}{\vm s}
\newcommand{\ttransfvectorb}{\dftemp_b}
\newcommand{\ttransfmatrix}{\overline{\transfmatrix}}
\newcommand{\residvector}{\vm r}
\newcommand{\nsubdivision}{s}
\newcommand{\freevar}{\vm z}

\newcommand{\rarg}{\vm x}

\newcommand{\temperature}{u}
\newcommand{\temperatureh}{\temperature_h}
\newcommand{\temperaturep}{\temperature^h}
\newcommand{\temperaturec}{\temperature_0}
\newcommand{\temperatureb}{\temperature_b}
\newcommand{\auxvariable}{\phi}

\newcommand{\cmtr}{d}
\newcommand{\cmtrboundary}{b}
\newcommand{\cmtroper}{D}
\newcommand{\cmtx}{\vm D}

\newcommand{\tcmtx}{\overline{\cmtx}}
\newcommand{\cmtvb}{\vm b}
\newcommand{\tcmtvb}{\overline{\cmtvb}}
\newcommand{\dftemp}{\vm q}
\newcommand{\dftempfree}{\dftemp_0}
\newcommand{\dftempsegment}{\vm a}

\newcommand{\funspace}{V}

\newcommand{\funspacedual}{\funspace^*}

\newcommand{\fpk}{\bm P}
\newcommand{\spk}{\bm S}
\newcommand{\Jacob}{\bm F}
\newcommand{\GLstrain}{\bm E}
\newcommand{\RightCG}{\bm C}
\newcommand{\unitmatrix}{\bm I}
\newcommand{\mattangent}{\mathbb{C}}

\newcommand{\realspace}{\mathbb{R}}
\newcommand{\mnull}{\ker}
\newcommand{\projectionp}{\mathcal{P}_\porder}

\newcommand{\paramone}{m_1}
\newcommand{\paramtwo}{m_2}
\newcommand{\paramthree}{m_3}

\newcommand{\funspacecz}{{\funspace^{c}_{0}}}
\newcommand{\funspacecb}{{\funspace^{c}_{b}}}
\newcommand{\approxfunspace}{{\funspace_{\porder}}}

\newcommand{\approxfunspaceelem}{{\funspace_{\porder}^e}}
\newcommand{\vectorspacecz}{{Z^{c}_{0\porder}}}
\newcommand{\vectorspacecb}{{Z^{c}_{b\porder}}}
\newcommand{\tvectorspacecz}{{\overline{Z}}^{c}_{0\porder}}
\newcommand{\tvectorspacecb}{{\overline{Z}}^{c}_{b\porder}}

\newcommand{\Gshapesc}{{\Psi}}
\newcommand{\GGshapesc}{\varPhi}
\newcommand{\Gshapev}{{\vm \Gshapesc}}
\newcommand{\GGshapev}{{\vm \GGshapesc}}
\newcommand{\Gshapesegment}{\vm \Phi}
\newcommand{\segmentmatrix}{\vm M}
\newcommand{\segmentmatrixred}{\overline{\segmentmatrix}}

\newcommand{\porder}{p}
\newcommand{\bporder}{{\overline{\porder}}}

\newcommand{\iinbbr}[1]{\left\{ #1 \right\}}
\newcommand{\inbbr}[1]{ [ \negthinspace[ #1 ] \negthinspace ] }

\newcommand{\displusk}{u}
\newcommand{\displu}{\mac{u}}

\newcommand{\norm}[1]{\lVert#1\rVert}

\newcommand{\abs}[1]{\lvert #1 \rvert}

\newcommand{\fvu}{\mac{v}}

\newcommand{\fvt}{{v}}

\newcommand{\fvtt}{{w}}

\newcommand{\scalarprod}{\cdot}

\newcommand{\Domain}{\varOmega}
\newcommand{\Boundary}{\varGamma}
\newcommand{\Boundarydirich}{\Boundary_D}
\newcommand{\Boundaryneum}{\Boundary_N}

\newcommand{\Domainfem}{\Domain_h}
\newcommand{\Domainelem}{\Domain^e}

\newcommand{\mint}{\int\limits}
\newcommand{\omint}{\mint_{\Domain}}
\newcommand{\intend}{\, \mathrm d}
\newcommand{\dOm}{\, d \bm{x}}
\newcommand{\gaint}[1]{\mint_{\Boundary_{#1}}}
\newcommand{\gaintt}[1]{\mint_{{#1}}}
\newcommand{\dGa}{\, dS}

\newcommand{\eye}{\bm I}
\newcommand{\smas}{\mac{b}}
\newcommand{\spow}{\mac{t}}
\newcommand{\youngsk}{E}
\newcommand{\young}{\mac{\youngsk}}

\newcommand{\strain}{\vmac{\varepsilon}}

\newcommand{\strainsk}{\varepsilon}

\newcommand{\hg}{f}
\newcommand{\hf}{h}

\newcommand{\nablav}{{\boldsymbol \nabla}}
\newcommand{\gradoper}{\nablav}
\newcommand{\gradoperh}{\gradoper_h}

\newcommand{\Vnormal}{\mac{n}}
\newcommand{\VVnormal}{\mac{N}}

\newcommand{\ssig}{\sigma}

\newcommand{\bsig}{\bs \ssig}
\newcommand{\sig}{\bsig}

\newcommand{\stiffmatrix}{\vm K}

\newcommand{\rightvector}{\vm F}

\newcommand{\dfu}{\mac{\check{u}}}

\newcommand{\dftempsc}{q}

\newcommand{\skeletonsign}{s}

\newcommand{\skeletonnormal}{\Vnormal}

\newcommand{\thefunction}{g}

\newcommand{\skeletonfem}{\Boundary_{\skeletonsign}}

\newcommand{\elemsize}{h}

\newcommand{\lamemu}{\mu}
\newcommand{\lamelambda}{\lambda}

\newcommand{\youngmodul}{E}

\newcommand{\cargone}{x}
\newcommand{\cargtwo}{y}

\newcommand{\eqcomma}{,}  
\newcommand{\eqdot}{{.}}

\newcommand{\bs}{\boldsymbol}
\newcommand{\mac}[1]{\mathbf {#1}}
\newcommand{\vm}[1]{\bm{#1}}

\newcommand{\vmac}[1]{\boldsymbol {#1}}
\newcommand{\tran}{^{\mathrm{T}}}

\newcommand{\logLogSlopeTriangle}[5]
{
	
	\pgfplotsextra
	{
		\pgfkeysgetvalue{/pgfplots/xmin}{\xmin}
		\pgfkeysgetvalue{/pgfplots/xmax}{\xmax}
		\pgfkeysgetvalue{/pgfplots/ymin}{\ymin}
		\pgfkeysgetvalue{/pgfplots/ymax}{\ymax}
		
		\pgfmathsetmacro{\xArel}{#1}
		\pgfmathsetmacro{\yArel}{#3}
		\pgfmathsetmacro{\xBrel}{#1-#2}
		\pgfmathsetmacro{\yBrel}{\yArel}
		\pgfmathsetmacro{\xCrel}{\xArel}
		
		\pgfmathsetmacro{\lnxB}{\xmin*(1-(#1-#2))+\xmax*(#1-#2)} 
		\pgfmathsetmacro{\lnxA}{\xmin*(1-#1)+\xmax*#1} 
		\pgfmathsetmacro{\lnyA}{\ymin*(1-#3)+\ymax*#3} 
		\pgfmathsetmacro{\lnyC}{\lnyA-#4*(\lnxA-\lnxB)}
		\pgfmathsetmacro{\yCrel}{\lnyC-\ymin)/(\ymax-\ymin)} 
		
		\coordinate (A) at (rel axis cs:\xArel,\yArel);
		\coordinate (B) at (rel axis cs:\xBrel,\yBrel);
		\coordinate (C) at (rel axis cs:\xCrel,\yCrel);
		
		\draw[#5]   (A)-- node[pos=0.5,anchor=south] {1}
		(B)--
		(C)-- node[pos=0.5,anchor=west] {#4}
		cycle;
	}
}

\newcommand{\logLogSlopeTTriangle}[5]
{
	
	\pgfplotsextra
	{
		\pgfkeysgetvalue{/pgfplots/xmin}{\xmin}
		\pgfkeysgetvalue{/pgfplots/xmax}{\xmax}
		\pgfkeysgetvalue{/pgfplots/ymin}{\ymin}
		\pgfkeysgetvalue{/pgfplots/ymax}{\ymax}
		
		\pgfmathsetmacro{\xArel}{#1}
		\pgfmathsetmacro{\yArel}{#3}
		\pgfmathsetmacro{\xBrel}{#1+#2}
		\pgfmathsetmacro{\yBrel}{\yArel}
		\pgfmathsetmacro{\xCrel}{\xArel}
		
		\pgfmathsetmacro{\lnxB}{\xmin*(1-(#1-#2))+\xmax*(#1-#2)} 
		\pgfmathsetmacro{\lnxA}{\xmin*(1-#1)+\xmax*#1} 
		\pgfmathsetmacro{\lnyA}{\ymin*(1-#3)+\ymax*#3} 
		\pgfmathsetmacro{\lnyC}{\lnyA-#4*(\lnxB-\lnxA)}
		\pgfmathsetmacro{\yCrel}{\lnyC-\ymin)/(\ymax-\ymin)}
		
		\coordinate (A) at (rel axis cs:\xArel,\yArel);
		\coordinate (B) at (rel axis cs:\xBrel,\yBrel);
		\coordinate (C) at (rel axis cs:\xCrel,\yCrel);
		
		\draw[#5]   (A)-- node[pos=0.5,anchor=north] {1}
		(B)--
		(C)-- node[pos=0.5,anchor=east] {#4}
		cycle;
		
	}
}

\title{Penalty-free discontinuous Galerkin method}
\author{Jan Ja{\'s}kowiec\footnote{jan.jaskowiec@pk.edu.pl, Faculty of Civil Engineering, Cracow University of Technology, Cracow Poland}
	\and 
	N. Sukumar\footnote{nsukumar@ucdavis.edu, Department of Civil and Environmental Engineering, University of California, Davis, USA}} 
\date{}

\begin{document}
		
\maketitle

\begin{abstract}
	In this paper, we present a new high-order discontinuous Galerkin (DG) method, in which neither a penalty parameter nor a stabilization parameter is needed. We refer to this method as penalty-free DG (\PFDG). In this method, the trial and test functions belong to the broken Sobolev space, in which the functions are in general discontinuous on the mesh skeleton and do not meet the Dirichlet boundary conditions. However, a  subset can be distinguished in this space, where the functions are continuous and satisfy the Dirichlet boundary conditions, and this subset is called admissible. The trial solution is chosen to lie in an \emph{augmented}  admissible subset, in which a small violation of the  continuity condition is permitted. This subset is constructed by applying special augmented constraints to the linear combination of finite element basis functions. In this approach, all the advantages of the DG method are retained without the necessity of using stability parameters or numerical fluxes. Several benchmark problems in two dimensions (Poisson equation, linear elasticity, hyperelasticity, and biharmonic equation) on polygonal (triangles, quadrilateral and weakly convex polygons) meshes as well as a three-dimensional Poisson problem on hexahedral meshes are considered.  Numerical results are presented that affirm the sound accuracy and optimal convergence of the method in the $L^2$ norm and the energy seminorm. 
\end{abstract}
%

\section{Introduction}
The discontinuous Galerkin (DG) method has a long history that dates back to the early seventies of the previous century~\cite{Reed_Hill_1973_techrep}.
Since that time, the method has undergone constant development and many monographs have appeared on this topic~\cite{opac-b1096869,di2011mathematical,Hesthaven:2007,Riviere2008}.
It is a well-known property of DG methods that they are applicable to discretizations that produce very general meshes. The finite elements that are used in the DG method can be arbitrary (convex or nonconvex) polygonal and polyhedral shapes. The approximation space in these methods is populated by polynomials of any degree on an element without the a priori need to link the polynomials on adjacent elements~\cite{jaskowiec_dg_finel_2016,jaskowiec_cm_2017,jaskowiec_cmame_2017,cangiani2017hp}. In the DG methods, unlike in the finite element method (FEM), the approximation functions across contiguous finite elements are discontinuous. Consequently, integration along the mesh skeleton (inter-element boundaries) is needed. Depending on the version of the DG method used, the continuity of the final solution is enforced by applying different techniques to the mesh skeleton. In the interior penalty DG (IPDG) method, the most popular one in the literature~\cite{Mu2014432,FLD:FLD3653,e2008discontinuourivis}, continuity is obtained via Nitsche's method~\cite{Nitsche1971}, which utilizes numerical fluxes applied to the skeleton and an additional penalty-like term, referred to as the stabilization parameter. 
An adaptive stabilization strategy in the DG method for nonlinear elasticity problems has been proposed~\cite{EYCK20083605,EYCK20082989}, in which for better approximation, the stabilization term is chosen to adapt to the problem.
In the local DG (LDG) methods~\cite{cockburn_1998,CASTILLO20061307,NUM:NUM22150}, the Bassi-Rebay DG~\cite{BASSI1997267,Eyck2006}, Baumann-Oden DG~\cite{BAUMANN1999311}, or mixed formulation DG~\cite{Dautov2013}, two numerical fluxes are used for the primary and secondary fields. This is due to the ultra-weak problem formulation in which a second-order problem is written as two first-order equations. These numerical fluxes utilize additional parameters that must be evaluated. An alternative approach is adopted in the hybridized DG (HDG) method~\cite{cockburn_2004,Cockburn2009,Samii2016,Cockburn2016,KABARIA2015303}, discontinuous Petrov-Galerkin (DGP) method~\cite{Demkowicz2014,KEITH2017107}, or DG method with Lagrange multipliers (DGLM)~\cite{Borker2017}, where the numerical fluxes are treated as additional unknowns. Among the aforementioned DG methods, the IPDG methods are the most widely used. The symmetric version of the method, SIPDG, dates back to a paper by Douglas~\cite{Douglas1982} and has been thoroughly studied~\cite{Riviere2008,BIRD201978}. The main drawback of this approach is the necessity to evaluate the penalty (or stability) parameters. In previous studies, significant emphasis has been placed on the stability of various versions of the IPDG method~\cite{EPSHTEYN2007843,BREZZI20063293,Cyril2013,Feng2018}. The necessity of setting the stability parameters is the weakest link in penalty-based DG methods. For further details, the interested reader can refer to a few pertinent references on DG` \cite{douglas_2002,Kirby2005,BREZZI20063293,Antonietti2016}.

In this paper, we introduce a penalty-free discontinuous Galerkin (PF-DG) method that is applicable to standard finite elements as well as general polygonal and polyhedral meshes. Use of such general meshes has appeared in polygonal and polyhedral finite element methods,~\cite{Hormann:2017:GBC} DG methods on polytopes,~\cite{cangiani2017hp} and the virtual element method.\cite{Beirao:2013:BPV} Unlike classical DG methods and the virtual element method that require stabilization, the PF-DG method proposed herein does not need a penalty nor stabilization parameter, while retaining the flexibility of being applicable on polytopal meshes.
In the PF-DG method, the finite-dimensional broken Sobolev space is used to search for approximate solutions in which some equality constraints are defined, which allows one to choose only those functions that are continuous and satisfy the inhomogeneous Dirichlet boundary conditions. These constraints are applied using the least squares approach on the mesh skeleton as well as the outer boundary, and as a consequence, a subset of admissible functions is obtained, 
which consists only of continuous functions that meet the boundary conditions on the Dirichlet boundary. The admissible functions are constructed by solving a least-square constraints problem. However, numerical computation of these constraints on polygonal meshes proves to be overly restrictive, which causes numerical instabilities. 
Thus, an \emph{augmented} admissible subset is used so that very small discontinuities on the mesh skeleton are allowed. The constraint problem leads to a linear system of equations that is solved prior to the main discrete problem. The matrix in the constraints problem is rectangular and singular, which admits a nonunique solution, and the null-space method~\cite{scott2022null,scott2022solving} is applied to obtain the complete solution. It should be noted here that in Oden at al.~\cite{ODEN1998491}, a penalty-free DG method is proposed, but this method has not been widely adopted due to its low level of accuracy when compared to penalty-based approaches. In this paper, we compare the penalty-free approach of Oden et al.\ to the one proposed herein on one benchmark problem. 
Another DG method without penalty parameters was presented by John et al.~\cite{john2016stable}, where the lifting operator on the mesh skeleton is introduced in the discrete model. The method exploits the properties of the piecewise Raviart–Thomas–N{\'e}d{\'e}lec element in which the lifting operator increases the polynomial space by one.
Other penalty-free versions of the DG method are available~\cite{Liu2011,Riviere2014}, but these are limited to Stokes flow. Continuity of the approximate solution is enforced using the method of Lagrange multipliers~\cite{FARHAT20031389, BROGNIEZ201349, KIM2015488}, however, it leads to a larger non-positive-definite problem. 
Even though there are some commonalities, the approach pursued in the \PFDG{} method is distinct from the enhanced strain method (ESM)~\cite{reddy1995stability,wriggers1996note}. In both these methods, special approximation subspaces are selected---in \PFDG{} the trial space must be (nearly) continuous and in the ESM the stress and enhanced strain subspaces are mutually orthogonal.
The \PFDG{} method presented in this paper is applicable on high-order finite elements and can be applied to different classes of boundary-value problems.

The outline of the remainder of this paper follows. In Section~\ref{sec:pfdg_basic}, a detailed derivation for the \PFDG{} method is presented for the Poisson problem. The new DG method requires the prior solution of a constrained algebraic problem. The details on how to construct the constraints and how to solve them are presented in Section~\ref{sec:constriants} (also see Appendices A--C). In  Sections~\ref{sec:pfdg_elasticity} and~\ref{sec:pfdg_4order}, the \PFDG{} method is presented for linear elasticity and the fourth-order biharmonic equation, respectively. Six two-dimensional benchmark problems that involve Poisson equation, linear and nonlinear elasticity, and the biharmonic equation are solved in Section~\ref{sec:benchmarks}. In addition, numerical results for a three-dimensional Poisson problem are presented in Section~\ref{ssec:4order_example}. For all problems, we show that the \PFDG{} method converges optimally in Sobolev norms. We close with some final remarks and conclusions in Section~\ref{sec:conclusions}.

\section{Discontinuous Galerkin method without stability parameter}\label{sec:pfdg_basic}

In this section, the PF-DG method---a DG method without a stability parameter---is presented. The mathematical formulation for the continuous problem is described in Section~\ref{ssec:mathematical}. The finite-dimensional formulation of the method is presented in Section~\ref{ssec:dg_three}.

\subsection{Mathematical formulation}\label{ssec:mathematical}
The penalty-free DG method is presented for the model Poisson problem that is defined in the domain $\Domain \subset \Real^d$ ($d = 1,\,2,\,3$) with boundary $\Boundary := \partial \Omega$. We seek the scalar field $\temperature : \Real^d \to \Real$ that solves the following Poisson boundary-value problem:
\begin{subequations}\label{eq:heatcond1}
	\begin{align}
		\label{eq:heatcond1_a}
		-\nabla^2\temperature &= \hg \ \ \text{in }\Domain , \\
		\label{eq:heatcond1_b}
		\temperature &= \boundarysign{\temperature} \ \ \text{on }\Boundarydirich , \\
		\label{eq:heatcond1_c}
		\gradoper\temperature \scalarprod \Vnormal &= \boundarysign{\hf} \ \ \text{on }\Boundaryneum \eqcomma 
	\end{align}
\end{subequations}
where $\hg\in L^2(\Domain)$ is the source density function, $\boundarysign{\temperature}$ and $\boundarysign{\hf}$ are prescribed values on the Dirichlet and Neumann boundaries $\Boundarydirich$ and $\Boundaryneum$, respectively,and $\Vnormal$ is the unit outward normal vector on $\Boundaryneum$. Throughout the paper standard notation for the $L^2$ inner product over the domain and its boundary is used:
\begin{align}
	\oinnerprod{\temperature,\fvt}{\Domain} = \omint \temperature \fvt \, d \bm{x} ,
	\quad
	\binnerprod{\temperature,\fvt}{\Boundary} = \gaintt{\Boundary} \temperature \fvt \, dS \eqdot
\end{align}

It is assumed that $\Domain$ is covered by the set $\Domainfem$ that consists of non-overlapping cells $\Domainelem\in\Domainfem$, which are called  finite elements. The interior boundaries of the $\Domainfem$ create the mesh skeleton $\skeletonfem$. On $\Domain$, the broken Sobolev space is defined as:
\begin{align}\label{eq:broken_sobolev}
	\funspace = 
	\left\{ w \in L^2(\Domain):\, w\Bigl|_{\Domain^e} \in H^1(\Domain^e) \ \ \forall \Domain^e \in \Domain_h \right\} \eqdot
\end{align}
The space $\funspace$ consists of element-wise functions that are in $H^1(\Omega^e)$ but can be discontinuous on the mesh skeleton. Two operators on the mesh skeleton are defined that are associated with the discontinuity: the jump operator $\inbbr{\bullet}$ and the mid-value operator $\iinbbr{\bullet}$, which are defined for an auxiliary function $\thefunction$ as follows:
\begin{align}
	\label{eq:def-disc}
	\begin{aligned}
		&
		\inbbr{\thefunction} = \lim_{\epsilon \rightarrow 0} \, \Bigl(	\thefunction(\bm x + \epsilon \, \skeletonnormal) - \thefunction(\bm x - \epsilon \, \skeletonnormal) \Bigr) \eqcomma
		\\&
		\iinbbr{\thefunction} = \frac{1}{2} \lim_{\epsilon \rightarrow 0}\, \Bigl( \thefunction(\bm x + \epsilon \, \skeletonnormal)
		+
		\thefunction(\bm x - \epsilon \, \skeletonnormal) \Bigr) \eqcomma
	\end{aligned}
\end{align}
where $\skeletonnormal$ is the unit  normal vector to $\skeletonfem$. 

Two  disjoint subsets in $\funspace$ are distinguished, which both consist of continuous functions and functions that are homogeneous  on the Dirichlet boundary in the first set (test functions) and satisfy the boundary conditions in~\eqref{eq:heatcond1} in the second set (trial functions): \begin{subequations}\label{eq:function_sets}
	\begin{align}
		&\funspacecz = \left\{\fvtt \in \funspace: \inbbr{\fvtt}=0 \text{ on } \skeletonfem, \
		\fvtt = 0 \text{ on } \Boundarydirich \right\} \eqcomma
		\\
		&\funspacecb = \left\{\fvtt \in \funspace: \inbbr{\fvtt}=0 \text{ on } \skeletonfem , \
		\fvtt = \boundarysign\temperature \text{ on } \Boundarydirich \right\} \eqcomma
		\\
		& \funspacecz \cap \funspacecb = \emptyset \eqdot
	\end{align}
\end{subequations}
In the case when $\boundarysign\temperature\equiv 0$ on $\Boundarydirich$ then $\funspacecb \rightarrow \funspacecz$. Among the functions in $\funspacecb$, we want to find the solution of the boundary-value problem, which we refer to as an \emph{admissible subset}. Note that an arbitrary function in $\funspacecb$ is a linear combination of any function in $\funspacecz$ and a continuous function $u_b$ whose trace on the Dirichlet boundary is $\bar{u}$. It should be observed that $\funspacecz = H^1_0(\Domain)$, 
which indicates that $\funspacecz$ is the Sobolev space inside the broken Sobolev space $\funspace$.

The formulation of the DG method starts with a Galerkin weak form of the strong form in~\myeqref{eq:heatcond1} with $\fvt$ as the test function: Find $\temperature \in \funspacecb$ such that
\begin{align}	\label{eq:heatweak1}
	-\oinnerprod{\gradoper^2\temperature,\fvt}{\Domain} - \oinnerprod{\hg, \fvt}{\Domain} = 0 \ \ \forall  \fvt \in \funspace \eqdot
\end{align}
On using the divergence theorem in~\eqref{eq:heatweak1} and taking into account the discontinuity of functions in $\funspace$ on the mesh skeleton, we get
\begin{align*}
	\oinnerprod{\gradoper\temperature, \gradoper\fvt}{\Domain} 
	+ \gaint{\skeletonsign} \inbbr{\gradoper\temperature\fvt}\scalarprod \skeletonnormal\intend s
	- \binnerprod{\gradoper\temperature\scalarprod\Vnormal, \fvt}{\Boundarydirich} 
	- \oinnerprod{\hg, \fvt}{\Domain} = 0
	\eqdot
\end{align*}
On using the property of the jump operator that $\inbbr{gh} = \iinbbr{g}\inbbr{fh}+\inbbr{g}\iinbbr{h}$, we have
\begin{align}\label{eq:heatweak2}
	\oinnerprod{\gradoper\temperature, \gradoper\fvt}{\Domain} 
	+ \binnerprod{\iinbbr{\gradoper\temperature}\scalarprod\skeletonnormal, \inbbr{\fvt}}{\skeletonfem}
	+ \cancelto{0}{\binnerprod{\inbbr{\gradoper\temperature}\scalarprod\skeletonnormal, \iinbbr{\fvt}\skeletonnormal}{\skeletonfem}}
	- \binnerprod{\gradoper\temperature\scalarprod\Vnormal, \fvt}{\Boundarydirich} 
	- \binnerprod{\boundarysign{\hf}, \fvt}{\Boundaryneum} 
	- \oinnerprod{\hg, \fvt}{\Domain} = 0
	\eqdot
\end{align}
The third component in~\eqref{eq:heatweak2} goes to zero due to the continuity of the exact solution $\temperature$ and its derivatives on the mesh skeleton. The fifth component uses the Neumann boundary condition in~\eqref{eq:heatcond1_c}. In the finite-dimensional problem, the same is applied to the approximate solution $\temperaturep\approx \temperature$ even though $\inbbr{\gradoper\temperaturep} \scalarprod \skeletonnormal \neq 0$ and $\gradoper\temperaturep\scalarprod \Vnormal \neq \boundarysign{\hf}$ on the mesh skeleton and Neumann boundary, respectively.  In this case, a weak enforcement of the gradient continuity and Neumann boundary condition of the approximate solution ensues. It results in the following weak formulation of the Poisson problem: Find $\temperature \in \funspacecb$ such that
\begin{subequations}\label{eq:heatweak2a}
	\begin{align}
		\bilinearoper(\temperature,\fvt) &= \linearoper(\fvt) \ \ \forall \fvt \in \funspace \eqcomma \label{eq:heatweak2a-a}
		\intertext{where}
		\bilinearoper(\temperature, \fvt) & = \oinnerprod{\gradoperh\temperature, \gradoperh\fvt}{\Domain} 
		+ \binnerprod{\iinbbr{\gradoperh\temperature}, \inbbr{\fvt}\skeletonnormal}{\skeletonfem}
		- \binnerprod{\gradoperh\temperature, \fvt\Vnormal}{\Boundarydirich} \eqcomma
		\label{eq:heatweak2a-b}
		\\
		\linearoper(\fvt) &= \oinnerprod{\hg, \fvt}{\Domain} + \binnerprod{\boundarysign{\hf}, \fvt}{\Boundaryneum} \eqdot
		\label{eq:heatweak2a-c}
	\end{align}
\end{subequations}
where $\gradoperh$ denotes the element-wise gradient operator.

The bilinear functional $\bilinearoper(\fvt,\temperature)$ is not symmetric due to the integrals along the mesh skeleton and the Dirichlet boundary. The function $\temperature$ that solves~\eqref{eq:heatweak2a} belongs to the space $\funspacecb$, which suggests that following the procedure presented in Oden et al.:\cite{ODEN1998491} the two zero-value components, namely $\binnerprod{\inbbr{\temperature} \skeletonnormal , \iinbbr{\gradoperh\fvt}}{\skeletonfem}$ and $\binnerprod{\left(\boundarysign{\temperature} - \temperature \right)\skeletonnormal, \gradoperh\fvt}{\Boundarydirich}$, can be added to~\eqref{eq:heatweak2a-b} to make the functional $\bilinearoper(\fvt,\temperature)$ symmetric or skew-symmetric. The resulting formulation is well-defined and can be used to find the approximate solution. Later in Section~\ref{sec:benchmarks}, we compare the method of Oden at al.\cite{ODEN1998491} with the proposed \PFDG{} method on a benchmark problem.

The functions in the space $\funspace$ are generally discontinuous on the mesh skeleton and their values on the Dirichlet boundary differ from $ \boundarysign{\temperature}$. In order to solve the problem defined in~\eqref{eq:heatweak2} the subset $\funspacecb$ has to be first distinguished from the space $\funspace$. To find the set of functions that satisfy  both conditions, we require that a function $\auxvariable\in\funspace$ meets the following least-squares conditions:
\begin{subequations}\label{eq:side_conditions}
	\begin{align}
		\label{eq:side_conditions_a}
		&\arg\min_{\auxvariable\in\funspace} \binnerprod{\inbbr{\auxvariable},\inbbr{\auxvariable}}{\skeletonfem}
		\, \Rightarrow \,
		\binnerprod{\inbbr{\auxvariable},\inbbr{\fvt}}{\skeletonfem} = 0 \ \ \forall \fvt \in \funspace
		\eqcomma
		\\&\label{eq:side_conditions_b}
		\arg\min_{\auxvariable\in\funspace} \binnerprod{\auxvariable - \boundarysign{\temperature},\auxvariable - \boundarysign{\temperature}}{\Boundarydirich}
		\, \Rightarrow \,
		\binnerprod{\auxvariable,\fvt}{\Boundarydirich} - \binnerprod{\boundarysign{\temperature},\fvt}{\Boundarydirich} = 0
		\ \  \forall \fvt \in \funspace
		\eqdot
	\end{align}
\end{subequations}

Since $\skeletonfem \cap \Boundarydirich$ is of measure zero,  the problem written in~\eqref{eq:side_conditions_a} and~\eqref{eq:side_conditions_b} can be placed in a single equation: Find $\auxvariable \in \funspace$ such that
\begin{subequations}\label{eq:dg_bd_cond_equation}
	\begin{align}\label{eq:dg_bd_cond_equation_a}
		\cmtr(\auxvariable,\fvt) = \cmtrboundary(\fvt) \ \ \forall \fvt \in \funspace \eqcomma
	\end{align}
	where
	\begin{align}\label{eq:dg_bd_cond_equation_b}
		&	\cmtr(\auxvariable,\fvt) = \binnerprod{\inbbr{\auxvariable},\inbbr{\fvt}}{\skeletonfem} +  \binnerprod{\auxvariable,\fvt}{\Boundarydirich}
		,\quad
		\cmtrboundary(\fvt) = \binnerprod{\boundarysign{\temperature},\fvt}{\Boundarydirich} \eqdot
	\end{align}
\end{subequations}

The solution to the problem defined in~\eqref{eq:dg_bd_cond_equation} is nonunique. It can be shown that all possible solutions of the problem define the subset $\funspacecb$ in $\funspace$  and the solution of the homogeneous version of~\eqref{eq:dg_bd_cond_equation} defines $\funspacecz\subset\funspace$. 

The following operators are introduced:
\begin{subequations}
	\begin{align}
		&&\bilinearoperoper : \funspace \rightarrow \funspacedual\,,\quad
		\left( \bilinearoperoper\temperature,\fvt \right) &= \bilinearoper(\temperature,\fvt) 
		\qquad \forall \fvt\in\funspace ,
		\\
		&&\cmtroper: \funspace \rightarrow \funspacedual\,,\quad
		\left( \cmtroper\temperature,\fvt \right) &= \cmtr(\temperature,\fvt) 
		\qquad \forall \fvt\in\funspace .
	\end{align}
\end{subequations}

The operator $\cmtroper$ is neither surjective nor injective. The operator $\bilinearoperoper$ need not be surjective or injective; however, the common null space of $\cmtroper$ and $\bilinearoperoper$ is trivial, $\mnull(\bilinearoperoper)\cap \mnull(\cmtroper)=\{0\}$. The operator $\cmtroper$ is not surjective, so its range is the subspace of the dual space $\range(\cmtroper) \subset \funspacedual$.

The following two theorems show the relations between the $\funspacecz$ and $\funspacecb$ and the~\eqref{eq:dg_bd_cond_equation}.

\begin{theorem}\label{th:theorem1}
	The null space of the $\cmtroper$ is the same as the space $\funspacecz$, i.e. $\mnull(\cmtroper)=\funspacecz$.
\end{theorem}

\begin{proof}
	In order to show that $\mnull(\cmtroper) = \funspacecz$ it is sufficient to establish that 
	$\funspacecz\subset\mnull(\cmtroper)$ and $\mnull(\cmtroper)\subset\funspacecz$.
	
	\begin{enumerate}[(i)]
		\item First, we show that $\funspacecz\subset\mnull(\cmtroper)$.
		Consider a function $\fvtt\in\funspacecz$. Let us assume that there is $g\in\cmtroper^\ast$ such that
		\begin{align*}
			\cmtroper\fvtt=g 
			\quad \iff \quad
			\cmtr(\fvtt,\fvt) = g(\fvt) \quad \forall \, \fvt \in \funspace .
		\end{align*} 
		Using the definition of the operator $\cmtr$ in~\eqref{eq:dg_bd_cond_equation_b}, we have
		\begin{align*}
			\binnerprod{\inbbr{\fvtt},\inbbr{\fvt}}{\skeletonfem} +  \binnerprod{\fvtt,\fvt}{\Boundarydirich} = g(\fvt) .
		\end{align*}
		The function $\fvtt$ is continuous in the domain and goes to zero on the Dirichlet boundary, so 
		\begin{align*}
			&
			\binnerprod{\inbbr{\fvtt},\inbbr{\fvt}}{\skeletonfem}=0
			\quad \text{and} \quad
			\binnerprod{\fvtt,\fvt}{\Boundarydirich} = 0
			\\&
			\Longrightarrow \quad
			g(\fvt)= 0  
			\quad \Longrightarrow \quad
			\cmtroper\fvtt=0
			\quad \Longrightarrow \quad
			\fvtt \in \mnull(\cmtroper) .
		\end{align*}
		
		\item Secondly, we show that $\mnull(\cmtroper)\subset\funspacecz$. 
		Consider a function $\fvtt\in\mnull(\cmtroper)$, so we have
		\begin{align*}
			&
			\cmtroper\fvtt=0
			\quad \iff \quad
			\cmtr(\fvtt,\fvt) = 0 \quad \forall \, \fvt \in \funspace 
			\\&
			\binnerprod{\inbbr{\fvtt},\inbbr{\fvt}}{\skeletonfem} +  \binnerprod{\fvtt,\fvt}{\Boundarydirich} = 0  \quad \forall \, \fvt \in \funspace .
		\end{align*}
		Because $\skeletonfem\cap \Boundarydirich = \emptyset$, it implies that
		\begin{align*}
			\left.
			\begin{aligned}
				&
				\binnerprod{\inbbr{\fvtt},\inbbr{\fvt}}{\skeletonfem} = 0 \quad \forall \, \fvt \in \funspace 
				\quad \iff \quad
				\inbbr{\fvtt} = 0 \quad \text{on } \skeletonfem
				\\&
				\binnerprod{\fvtt,\fvt}{\Boundarydirich} = 0  \quad \forall \, \fvt \in \funspace
				\quad \iff \quad
				\fvtt = 0 \quad \text{on } \Boundarydirich  
			\end{aligned}
			\,
			\right\} \quad \Longrightarrow \quad \fvtt\in\funspacecz  .
		\end{align*}
	\end{enumerate}
	Hence, $\mnull(\cmtroper)$ is the Sobolev space $H^1_0(\Domain)$.
\end{proof}

\begin{theorem}\label{th:theorem2}
	There exists a function $ \fvtt \in \funspacecb$ 
	if and only if
	$\cmtroper\fvtt = \cmtrboundary$.
\end{theorem}
The proof of Theorem~\ref{th:theorem2} 
follows similar lines as that carried out for Theorem~\ref{th:theorem1} and hence is omitted.

Using Theorems~\ref{th:theorem1} and~\ref{th:theorem2}, practical definitions for $\funspacecb$ and $\funspacecz$ can now be stated:
\begin{subequations}
	\begin{align}
		&\funspacecb =\left\{	\fvtt \in \funspace: \cmtr(\fvtt,\fvt) - \cmtrboundary(\fvt) = 0 \ \ \forall  \fvt \in \funspace \right\} \eqcomma
		\\&
		\funspacecz =\left\{	\fvtt \in \funspace: \cmtr(\fvtt,\fvt) = 0 \ \ \forall  \fvt \in \funspace  \right\} \eqdot
	\end{align}
\end{subequations}

At this moment, it is clear that on solving the problem in~\eqref{eq:dg_bd_cond_equation}, we can obtain functions that lie in $\funspacecb$.  Let $\temperature_b$ be a particular function that is chosen from the solutions contained in $\funspacecb$. The function $\temperature$, which we want to find can be now expressed as the  superposition of  two  functions:
\begin{align}\label{eq:temp_decomp}
	\temperature = \temperaturec + \temperatureb \quad \text{ in } \Domain \eqcomma
\end{align}
where $\temperaturec \in \funspacecz$ is the function to be found and $\temperatureb\in\funspacecb$ is a known function. Now, the problem defined in \eqref{eq:heatweak2} is reformulated: Find $\temperaturec \in \funspacecz$ such that
\begin{align}\label{eq:heatweak3}
	\bilinearoper(\temperaturec, \fvt) = \linearoper(\fvt) - \bilinearoper(\temperatureb, \fvt) \ \ \forall \fvt \in \funspacecz \eqdot
\end{align}
It should be noted that the following two components in $\bilinearoper(\temperature,\fvt)$ vanish due to fact that the test function belongs to $\funspacecz$:
\begin{align*}
	&\binnerprod{\iinbbr{\gradoperh\temperature}, \inbbr{\fvt}\skeletonnormal}{\skeletonfem} = 0
	,\quad
	\binnerprod{\gradoperh\temperature, \fvt\Vnormal}{\Boundarydirich} = 0 \eqdot
\end{align*}

Subsequently, we show the uniqueness of the solution for the problem presented in~\eqref{eq:heatweak3}.

\begin{theorem}\label{th:well_posedness_inf}
	Assume that the following holds:
	\begin{enumerate}[(i)]
		\item \label{it:coercivity} Coercivity in $\mnull(\cmtroper)$: There exists $c_1$ such that
		\begin{align*}
			\bilinearoper(\temperature,\temperature) \ge c_1 \norm{\temperature}_{\funspace}^2
			\quad \forall\, \temperature \in \mnull(\cmtroper) \eqdot
		\end{align*}
		
		\item Continuity: \label{it:continuity} There exist $c_2,\,c_3 >0$ such that: 
		\begin{align*}
			\bilinearoper(\temperature,\fvt)  &\le c_2 \norm{\temperature}_{\funspace}\, \norm{\fvt}_{\funspace}  	\quad \forall\, \temperature,\,\fvt \in \funspace \eqcomma
			\\
			\linearoper(\fvt) &\le c_3 \norm{\fvt}_{\funspace} \quad \forall\, \fvt \in \funspace \eqdot
		\end{align*}

		\item \label{it:infsup} Inf-sup condition for constraints: There exist $c_4$ such that
		\begin{align*}
			&
			\inf_{\temperature\in\funspace} \sup_{\fvt\in\funspace} 
			\frac{\cmtr(\temperature,\fvt)}{\norm{\temperature}_{\funspace} \norm{\fvt}_{\funspace}}
			\ge c_4
			\quad \rightarrow \quad
			\norm{\cmtroper\temperature}_{\funspacedual} \ge c_4\norm{\temperature}_{\funspace} \quad \forall\,\temperature\in\funspace
			\eqdot
		\end{align*}
	\end{enumerate}
	
	Then there exists a unique solution $\temperature\in \mnull(\cmtroper)$ of the problem in~\eqref{eq:heatweak3} and
	\begin{align}\label{eq:sol_ineq}
		\norm{\temperature}_{\funspace} \leq \frac{1}{c_1}
		\left(\norm{\linearoper}_{\funspacedual} + c_2 \norm{\temperatureb}_{\funspace} \right) \eqdot
	\end{align}
\end{theorem}

\begin{proof}
	Since $\range(\cmtroper) \ne \emptyset$, there exists $\temperatureb \in \funspace$ such that
	\begin{align*}
		\cmtroper\temperatureb=\cmtrboundary \quad \rightarrow \quad
		\cmtr(\temperatureb,\fvt) = \cmtrboundary(\fvt) \quad \forall \fvt \in \funspace
	\end{align*}
	The operator $\cmtroper$ is not injective, so $\temperatureb$ is not unique. According to the inf-sup condition~\eqref{it:infsup} we have
	\begin{align}
		c_4\norm{\temperatureb}_{\funspace} \leq 
		\sup_{\fvt\in\funspace}  \frac{\cmtr(\temperatureb,\fvt)}{ \norm{\fvt}_{\funspace}}
		= \sup_{\fvt\in\funspace}  \frac{\cmtrboundary(\fvt)}{ \norm{\fvt}_{\funspace}} =
		\norm{\cmtrboundary}_{\funspacedual} 
		\quad \rightarrow \quad
		\norm{\temperatureb}_{\funspace} \leq  \frac{1}{c_4}\norm{\cmtrboundary}_{\funspacedual} 
	\end{align}
	
	Let us consider $\linearoper - \bilinearoperoper\temperatureb \in \funspacedual$. The coercivity condition~\eqref{it:coercivity} means that $\bilinearoperoper$ is an isomorphism from $\mnull(\cmtroper)\subset\funspace$ to $\mnull(\cmtroper)^\ast\subset\funspacedual$, and therefore we can find a unique $\temperaturec\in\mnull(\cmtroper)$ such that
	\begin{align}\label{eq:proof_zeqproblem}
		\bilinearoperoper\temperaturec = \linearoper - \bilinearoperoper\temperatureb 
		\quad \text{in} \quad \mnull(\cmtroper)^\ast \eqcomma
	\end{align}
	and
	\begin{align*}
		&
		c_1 \norm{\temperaturec}_{\funspace} \leq \frac{\bilinearoper(\temperaturec,\temperaturec)}{\norm{\temperaturec}_{\funspace}}
		\leq \sup_{\fvt\in\mnull(\cmtroper)}
		\frac{\bilinearoper(\temperaturec,\fvt)}{\norm{\fvt}_{\funspace}}
		= \norm{\bilinearoperoper\temperature}_{\mnull(\cmtroper)^\ast} 
		= \norm{\linearoper - \bilinearoperoper\temperatureb }_{\mnull(\cmtroper)^\ast} 
		\leq \norm{\linearoper - \bilinearoperoper\temperatureb }_{\funspacedual}
		\leq \norm{\linearoper}_{\funspacedual} + \norm{\bilinearoperoper\temperatureb}_{\funspacedual}
		\leq \norm{\linearoper}_{\funspacedual} + c_2 \norm{\temperatureb}_{\funspace}
		\\&
		\rightarrow \quad
		\norm{\temperaturec}_{\funspace} \leq \frac{1}{c_1}
		\left(\norm{\linearoper}_{\funspacedual} + c_2 \norm{\temperatureb}_{\funspace} \right)
		\eqcomma
	\end{align*}
	where we have used that $f\in\funspacedual$ such that 
	$\norm{f}_{\mnull(\cmtroper)^\ast} \leq \norm{f}_{\funspacedual}$ by definition of the dual norm and the fact that $\mnull(\cmtroper)\subset\funspace$ and so $\mnull(\cmtroper)^\ast\subset\funspacedual$.
	
	Now set $\temperature=\temperaturec+\temperatureb$, and due 
	to~\eqref{eq:proof_zeqproblem}, the following weak statement of equilibrium is satisfied:
	$\bilinearoper(\temperature,\fvt) = \linearoper(\fvt) \ \ \forall \fvt \in \funspacecz$.
\end{proof}

The admissible space is the null space of the operator $\cmtroper$, which consists of continuous functions that vanish on the Dirichlet boundary. We consider an augmented admissible space where some small levels of discontinuities are allowed and $\overline{\mnull}(\cmtroper)$ is defined as
\begin{align}\label{eq:augmented_ker}
	&
	\overline{\mnull}(\cmtroper) = \left\{ \auxvariable\in\funspace: \abs{\cmtr(\auxvariable,\fvt)} < \epsilon\norm{\auxvariable}_{\funspace} \ \ \forall \fvt \in \funspace\right\} ,
\end{align}
where $\epsilon>0$ is a small number. It is evident that
$\mnull(\cmtroper)\subset \overline{\mnull}(\cmtroper)\subset\funspace$.

The bilinear operator $\bilinearoper(\cdot,\cdot)$ consists of two components: the first on the domain and the second on the mesh skeleton and Dirichlet boundary, respectively:
\begin{align}\label{eq:bilinoper_augment}
	{\bilinearoper}(\temperature,\fvt) = \bilinearoper_{\Domain}(\temperature,\fvt) - \bilinearoper_{\Boundary}(\temperature,\fvt).
\end{align}%
The operator for the Poisson problem with the two components is presented in~\eqref{eq:heatweak2a-b}. The component $\bilinearoper_{\Boundary}$ vanishes for $\cmtr(\temperature,\fvt)=0$, which implies that
\begin{align}
	\sup_{\fvt \in \funspace}\frac{\abs{\bilinearoper_{\Boundary}(\temperature,\fvt)}}{\norm{\fvt}_{\funspace}} \leq \gamma(\epsilon)\norm{\temperature}_{\funspace}
	\quad \forall \, \temperature \in \overline{\mnull}(\cmtroper)
	\,,\quad \gamma(\epsilon) \xrightarrow[\epsilon\rightarrow 0]{} 0 .
\end{align}
The number $\gamma>0$ depends on $\epsilon$ from~\eqref{eq:augmented_ker} and $\gamma$ goes to zero when $\epsilon$ is zero. The operator $\bilinearoper$ is coercive in the admissible space, but can also remain coercive for the augmented admissible space for some small values of $\epsilon$, which is shown in the following conjecture.

\begin{conjecture}\label{th:augmentet_coerc}
	When the operator $\bilinearoper(\temperature,\fvt)$ is coercive in $\mnull(\cmtroper)$ with the constant $c_1$ then the operator remains coercive in $\overline{\mnull}(\cmtroper)$ for some small $\epsilon$.
\end{conjecture}

\begin{proof}

	It can be observed that 
	\begin{align*}
		\overline{\mnull}(\cmtroper) \xrightarrow[\epsilon\rightarrow 0]{} \mnull(\cmtroper)
		\,,\quad
		\overline{\mnull}(\cmtroper) \xrightarrow[\epsilon\rightarrow \infty]{} \funspace
	\end{align*}
	Consequently, the coercivity of the $\bilinearoper_{\Domain}$ operator changes with $\epsilon$ in the following manner:
	\begin{align*}
		\bilinearoper_{\Domain}(\temperature,\temperature) \ge c(\epsilon)\norm{\temperature}_{\funspace}^2
		\,,\quad
		c(\epsilon) \xrightarrow[\epsilon\rightarrow 0]{} c_1
		\,,\quad
		c(\epsilon) \xrightarrow[\epsilon\rightarrow \infty]{} 0.
	\end{align*} 
	
	The function $\c(\epsilon)$ is monotonically decreasing while the function $\gamma(\epsilon)$ is monotonically increasing as $\epsilon$ increases. According to~\eqref{eq:bilinoper_augment}, for any $\temperature\in\overline{\mnull}(\cmtroper)$ we have
	\begin{align*}
		{\bilinearoper}(\temperature,\temperature) = \bilinearoper_{\Domain}(\temperature,\temperature) - \bilinearoper_{\Boundary}(\temperature,\temperature)
		\ge
		\bilinearoper_{\Domain}(\temperature,\temperature) - \abs{\bilinearoper_{\Boundary}(\temperature,\temperature)}
		\ge
		c(\epsilon)\norm{\temperature}_{\funspace}^2 - \gamma(\epsilon)\norm{\temperature}_{\funspace}^2
		= \bigl(c(\epsilon) - \gamma(\epsilon) \bigr) \norm{\temperature}_{\funspace}^2.
	\end{align*}
	The operator $\bilinearoper(\temperature,\temperature)$ is coercive in $\overline{\mnull}(\cmtroper)$  so long as $\epsilon$ is small enough such that $c(\epsilon) > \gamma(\epsilon)$. \qedhere
\end{proof}

In the next section, the algorithm to compute the discrete bilinear form in the augmented admissible space is presented. The coercivity of the discrete bilinear 
form is numerically confirmed through the numerical example in Section~\eqref{ssec:poissons_dg} where we show that the stiffness matrix is well-conditioned.

%

\subsection{Approximate solution in the \PFDG{} method}\label{ssec:dg_three}
Now we apply the approach presented in the previous section to find an approximate solution of the boundary-value problem  in~\eqref{eq:heatweak3}. First, we need to rewrite constraint equations~\eqref{eq:dg_bd_cond_equation} in algebraic form. Then we present the construction of the finite-dimensional constraints and obtain the discrete system of equations that correspond to the weak form in~\eqref{eq:heatweak3}.

The domain $\Domain$ is covered by a finite element mesh $\Domainfem$, and the elements in the mesh $\Domainelem \in \Domainfem$ can be arbitrary polygons in 2D or polyhedrons in 3D with affine edges and planar faces, respectively. The approximate solution is to be found in the approximation space $\approxfunspace \subset \funspace$, which is constructed on a set of complete, polynomial basis functions up to order $\porder$ in each element $\approxfunspaceelem$, where $\approxfunspace=\bigcup\limits_e \approxfunspaceelem$. The linear subspace $\approxfunspaceelem$ is constructed by the basis functions $\{\Gshapesc^e_k\}_{k=1}^{m_p}$, where $\Gshapesc^e_k$ are the basis functions in the $e$-th element and the support of the basis functions are limited to the $e$-th element. The numbers of basis functions depend on the approximation order $p$ and the spatial dimension: $m_p=\porder+1$ in 1D, $m_p = \frac{1}{2}(\porder+1)(\porder+2)$ in 2D, and $m_p = \frac{1}{6}(\porder+1)(\porder+2)(\porder+3)$ in 3D. The polynomial basis functions can be: monomials, Legendre, Chebyshev, Lagrange polynomials, or any other basis. The approximation of a function $\temperature$ in this space is given by
\begin{align} \label{eq:glob_approx}
\temperature\approx\temperaturep = \Gshapev_\porder\dftemp \eqcomma
\end{align}
where 
$\Gshapev_\porder$ is a row vector that consists of all the basis functions in the mesh and $\dftemp \in \realspace^{m_p}$ is the vector of degrees of freedom for this approximation. In general, this approximation is not continuous on the mesh skeleton and does not satisfy the Dirichlet boundary conditions. The discrete equivalents of the sets $\funspacecz$ and $\funspacecb$ are defined as the following vector sets $\vectorspacecz$ and $\vectorspacecb$, respectively:
\begin{subequations}\label{eq:discrete_sets}
\begin{align}
	\label{eq:discrete_sets_a}
	&\vectorspacecz = \left\{ \vm v \in \realspace^{m_\porder}\,: \inbbr{\Gshapev_\porder} \vm v = 0 \ \text{on} \ \skeletonfem
	\ \text{and}\  \Gshapev_\porder\vm v = 0 \  \text{on} \ \Boundarydirich \right\} \eqcomma
	\\
	\label{eq:discrete_sets_b}
	&\vectorspacecb = \left\{ \vm v \in \realspace^{m_\porder}\,: \inbbr{\Gshapev_\porder} \vm v = 0 \   \text{on} \ \skeletonfem
	\ \text{and} \  \Gshapev_\porder\vm v = \projectionp \boundarysign{\temperature} \ \text{on} \ \Boundarydirich \right\}
	\eqdot
\end{align}
\end{subequations}
where the vector $\vm v$ in~\eqref{eq:discrete_sets} represents all vectors that meet the particular conditions. In addition, $\projectionp \boundarysign{\temperature}$ is the projection of the function $\boundarysign{\temperature}$ to $\approxfunspace\Bigr|_{\Boundarydirich}$, such that $\gaintt{\Boundarydirich} \Gshapev_\porder\tran \left(\projectionp \boundarysign{\temperature} -  \boundarysign{\temperature} \right) \dGa = \vm 0$.

The vector set $\vectorspacecb$, in combination with the basis functions, constitutes the finite admissible subset. In order to find the vector sets $\vectorspacecb$ and $\vectorspacecz$, the discrete form of~\eqref{eq:dg_bd_cond_equation} has to be solved, which reads
\begin{subequations}\label{eq:discrete_conds}
\begin{align}
	\cmtx \vm v &= \cmtvb \eqcomma 
	\qquad \cmtx \in \realspace^{m_p} \times \realspace^{m_p}\,,\, \cmtvb \in \realspace^{m_p}\,,\, \vm v \in \realspace^{m_p} \eqcomma
	\intertext{where}
	\cmtx &= \gaintt{\skeletonfem} \inbbr{\Gshapev_\porder}\tran\inbbr{\Gshapev_\porder} \dGa  +  
	\gaintt{\Boundarydirich} \Gshapev_\porder\tran\Gshapev_\porder \dGa,
	\quad
	\cmtvb =  \gaintt{\Boundarydirich} \Gshapev_\porder\tran \boundarysign{\temperature}\dGa \eqdot
\end{align}
\end{subequations}

The discrete version of Theorems~\ref{th:theorem1} and~\ref{th:theorem2} is presented in the two following Theorems.

\begin{theorem}
There is a vector $\vm v \in \vectorspacecb$ 
if and only if
$\cmtx \vm v = \cmtvb$.
\end{theorem}

\begin{proof}
\begin{enumerate}
	\item $\vm v \in \vectorspacecb$ $\Longrightarrow$ 	$\cmtx \vm v = \cmtvb$.
	The vector is applied to \eqref{eq:discrete_conds} and we have
	\begin{align*}
		\gaintt{\skeletonfem} \inbbr{\Gshapev_\porder}\tran\inbbr{\Gshapev_\porder} \dGa  \vm v + 
		\gaintt{\Boundarydirich} \Gshapev_\porder\tran\Gshapev_\porder \dGa \vm v
		=  \gaintt{\Boundarydirich} \Gshapev_\porder\tran \boundarysign{\temperature}\dGa \eqcomma
	\end{align*}
	or
	\begin{align*}
		\gaintt{\skeletonfem} \inbbr{\Gshapev_\porder}\tran\inbbr{\Gshapev_\porder} \vm v \dGa  + 
		\gaintt{\Boundarydirich} \Gshapev_\porder\tran \left( \Gshapev_\porder \vm v  - \boundarysign{\temperature}\right)\dGa =0 \eqcomma
	\end{align*}
	and since 	$ \inbbr{\Gshapev_\porder} \vm v = 0$
	and $ \Gshapev_\porder \vm v  = \projectionp \boundarysign{\temperature}$ and from~\eqref{eq:discrete_sets_b}, the equality is met.

	\item 	$\cmtx \vm v = \cmtvb$ $\Longrightarrow$ $\vm v \in \vectorspacecb$.
	
	Since the vector $\vm v$ satisfies \eqref{eq:discrete_conds}, so after some simple operations we can write   
	\begin{align*}
		\gaintt{\skeletonfem} \inbbr{\Gshapev_\porder}\tran\inbbr{\Gshapev_\porder} \vm v \dGa  + 
		\gaintt{\Boundarydirich} \Gshapev_\porder\tran \left(\Gshapev_\porder \vm v  - \boundarysign{\temperature}\right)\dGa =0 \eqdot
	\end{align*}
	Both integrals are independent of each other since $\skeletonfem \cap \Boundarydirich = \emptyset$, so we have
	\begin{align*}
		&\gaintt{\skeletonfem} \inbbr{\Gshapev_\porder}\tran\inbbr{\Gshapev_\porder} \vm v \,\dGa = 0   \quad\Longrightarrow\quad
		\inbbr{\Gshapev_\porder} \vm v = 0 \ \ \text{on } \skeletonfem \eqcomma
		\\ &
		\gaintt{\Boundarydirich} \Gshapev_\porder\tran \left(\Gshapev_\porder \vm v  - \boundarysign{\temperature}\right)\dGa =0
		\quad\Longrightarrow\quad \Gshapev_\porder \vm v = \projectionp\boundarysign{\temperature} \ \ \text{on } \Boundarydirich \eqdot
		\qedhere
	\end{align*}
\end{enumerate}
\end{proof}

\begin{theorem}
There is a vector $\vm v \in \vectorspacecz$ 
if and only if
$\cmtx \vm v = \vm 0$.
\end{theorem}
Now, the definitions of the sets $\vectorspacecz$ and $\vectorspacecb$ can be expressed as follows:
\begin{subequations}\label{eq:discrete_sets2}
\begin{align}
	\label{eq:discrete_sets2_a}
	&\vectorspacecz = \left\{ \vm v:\, \vm v \in  \ker{(\cmtx)} \right\} \eqcomma
	\\&
	\label{eq:discrete_sets2_b}
	\vectorspacecb = \left\{ \vm v:\, \norm{\cmtx\vm v - \cmtvb} = 0  \right\} \eqdot
\end{align}
\end{subequations}

The system of equations in~\eqref{eq:discrete_conds} is square but is singular. The rank of the left-hand side matrix is much smaller in comparison to its dimension. This implies that the system of equations has many linearly dependent equations. The solution of this system of equations is performed as described in Appendix~\ref{sec:dg_solver} and expressed in the form shown in~\eqref{eq:solution}, which is
\begin{align}\label{eq:dg_discr_constr_solve}
\vm v = \transfmatrix\freevar + \vm v_b \ \ \forall \freevar \eqcomma
\qquad \transfmatrix\in\realspace^{m_p} \times \realspace^{m_p-r}\,,\,
\freevar\in \realspace^{m_p-r}\,,\, \vm v_b \in \realspace^{m_p}
\eqcomma
\end{align}
where $r$ is the rank of $\cmtx$ matrix, $\transfmatrix$ is the matrix that consists of basis vectors in $\ker{(\cmtx)}$ and the vector $\vm v_b \in \vectorspacecb$, which is stated via the following theorem:
\begin{theorem}
If there is a vector $\vm v \in \vectorspacecb$ and $\vm v = \transfmatrix\freevar + \vm v_b$  $\forall \freevar$, then $\vm w = \transfmatrix\freevar \in \vectorspacecz$ and $\vm v_b \in \vectorspacecb$.
\end{theorem}

\begin{proof}
When $\vm v \in \vectorspacecb$, then $\vm v$ satisfies \eqref{eq:discrete_conds} and we have
\begin{align*}
	&\cmtx \transfmatrix\freevar + \cmtx\vm v_b = \cmtvb \ \ \forall \freevar
	\quad\Longleftrightarrow\quad \cmtx \transfmatrix\freevar = \vm 0 \ \ \forall \freevar
	\quad\Longrightarrow \quad \vm w \in \vectorspacecz \eqcomma
	\\&
	\cmtx\vm v_b = \cmtvb  \quad\Longleftrightarrow\quad \vm v_b \in \vectorspacecb \eqdot  \qedhere
\end{align*}
\end{proof}

Equation~\eqref{eq:discrete_conds} defines constraints for the vectors, and therefore to enforce continuity on the mesh skeleton and to satisfy the Dirichlet boundary conditions, we need to first solve this singular system of equations. In the calculations, these restrictions prove to be overly restrictive, especially for polygonal meshes, which as a consequence can lead to numerical instabilities that can lead to an inaccurate final solution. Hence, we \emph{loosen} these restrictions to realize a numerically stable approach. To this end, instead of~\eqref{eq:discrete_conds} we propose to use the following relation:
\begin{subequations}\label{eq:discrete_conds_ext}
\begin{align}
	\label{eq:discrete_conds_ext-a}
	{\tcmtx}\vm v &= \tcmtvb \eqcomma 
	\qquad \tcmtx\in\realspace^{m_\bporder}\times \realspace^{m_\bporder}
	\,,\, \tcmtvb \in \realspace^{m_\bporder}\eqcomma
	\\
	\intertext{where}
	\tcmtx &= \gaintt{\skeletonfem} \inbbr{\Gshapev_\bporder}\tran\inbbr{\Gshapev_\porder} \dGa  +  
	\gaintt{\Boundarydirich} \Gshapev_\bporder\tran\Gshapev_\porder \dGa,
	\quad
	\tcmtvb =  \gaintt{\Boundarydirich} \Gshapev_\bporder\tran \boundarysign{\temperature}\dGa \eqcomma \label{eq:discrete_conds_ext-b}
\end{align}
\end{subequations}
with $\bporder\leq\porder$. In the computations, we choose $\bporder = \porder - 1$. The matrix ${\tcmtx}$ is rectangular with the same number of columns as $\cmtx$ matrix, but  with fewer number of rows and so $\rank(\cmtx) > \rank(\tcmtx)\, \Rightarrow\, \ker(\cmtx) \subset \ker(\tcmtx)$. The set of all possible solutions of~\eqref{eq:discrete_conds_ext} is larger in comparison to the one in~\eqref{eq:discrete_conds}, allowing for some level of discontinuities or variation of the boundary conditions. The solution of the system in~\eqref{eq:discrete_conds_ext} has the following form:
\begin{align}\label{eq:dg_discr_constr_ext_solve}
\vm v = \ttransfmatrix\freevar + \ttransfvectorb \ \ \forall \freevar \,,
\qquad \ttransfmatrix\in \realspace^{m_\porder} \times \realspace^{m_\bporder-\overline r}
\,,\, \freevar \in \realspace^{m_\bporder-\overline r}
\,,\, \ttransfvectorb \in \realspace^{m_\porder}
\eqcomma
\end{align}
where $\overline r = \rank(\tcmtx)$.

Consequently, the sets of vectors in~\eqref{eq:discrete_sets2} now become \begin{subequations}\label{eq:discrete_sets2ext}
\begin{align}
	\label{eq:discrete_sets2ext_a}
	&\tvectorspacecz = \left\{ \vm v:\, \vm v \in  \ker{(\tcmtx)} \right\} \eqcomma
	\\&
	\label{eq:discrete_sets2ext_b}
	\tvectorspacecb = \left\{ \vm v:\, \norm{\tcmtx\vm v - \tcmtvb} = 0  \right\} \eqcomma
\end{align}
\end{subequations}
where $\tvectorspacecb$ together with the basis functions define the augmented admissible subset.

When we have the solution as expressed in\eqref{eq:dg_discr_constr_ext_solve}, then the discrete version of the weak problem formulation in~\eqref{eq:heatweak3} leads to the following problem:
For the given $\ttransfmatrix$ and $\ttransfvectorb$, find $\dftempfree$ such that
\begin{align}\label{eq:dg_constr_problem}
\ttransfmatrix\tran \vm K \ttransfmatrix \dftempfree = \ttransfmatrix\tran\left(\vm f - \vm K \ttransfvectorb \right) \,,
\qquad \vm K \in \realspace^{m_\porder}\times\realspace^{m_\porder}\,,\,
\vm f \in  \realspace^{m_\porder}
\end{align}
where $K_{ij}=\bilinearoper(\Gshapesc_j,\Gshapesc_i)$ and $f_i =  \linearoper(\Gshapesc_i)$ and the final solution vector is 
\begin{align}\label{eq:dg_constr}
\dftemp = \ttransfmatrix\dftempfree + \ttransfvectorb \eqdot
\end{align}

The matrix  $\ttransfmatrix\tran \vm K \ttransfmatrix$  is nonsymmetric and well-conditioned and can be solved by many solvers  for finding the $\dftempfree$ vector and later on to obtain the global vector $\dftemp$. However, the matrix $\ttransfmatrix$ and the vector $\ttransfvectorb$ have to be constructed a priori using the procedure presented in the following section. The well-posedness of the problem~\eqref{eq:dg_constr_problem} is shown in the following theorem.

\begin{theorem}
Under the assumptions of Theorem~\ref{th:well_posedness_inf} and Conjecture~\ref{th:augmentet_coerc} for any closed subspace $\approxfunspace\subset\funspace$, there exist a unique solution $\temperaturep = \Gshapev_\porder\dftemp \in \approxfunspace$ of the problem~\eqref{eq:dg_constr_problem}.
\end{theorem}

\begin{proof}
Consider the basis $\Gshapev_\porder$ in $\approxfunspace$, i.e. $\approxfunspace = \vspan{\Gshapev_\porder}$. Because $\ker(\tcmtx) = \vspan{\ttransfmatrix}$, so the function set $\GGshapev_\porder = \Gshapev_\porder \ttransfmatrix$ is  the basis functions of the space $\overline{\mnull}(\cmtroper)\Bigr|_{\approxfunspace} = \vspan{\GGshapev_\porder}$. The following matrix is constructed

\begin{align}
	\overline{ \vm K} = \bilinearoper(\GGshapev_\porder\tran,\GGshapev_\porder) =\ttransfmatrix\tran\bilinearoper(\Gshapev_\porder\tran,\Gshapev) \ttransfmatrix 	=\ttransfmatrix\tran \vm K \ttransfmatrix \eqdot
\end{align}
From the coercivity condition in Conjecture~\ref{th:augmentet_coerc}, we obtain that the matrix $\overline{\vm K}$ is positive-definite, and hence it is invertible.
\end{proof}

\section{Construction of the constraint equations}\label{sec:constriants}
In the \PFDG{} method presented herein, the constraints~\eqref{eq:discrete_conds_ext} that enforce the continuity and boundary conditions are the critical component of the approach. The solution of~\eqref{eq:dg_discr_constr_ext_solve} is nonunique due to the vector of free variables $\freevar$ (which is part of the vector $\vm v$ when using the approach described in Appendix~\ref{sssec:glob}) whose choice is in general arbitrary, and as a consequence, construction of $\ttransfmatrix$ and $\ttransfvectorb$ depend on the choice of the free variables. In the case when we deal with a relatively small problem with few low-order elements, the final result is independent of the choice of $\freevar$. However, for larger problems, choice of $\freevar$ becomes important, since otherwise truncation errors can accumulate in $\ttransfmatrix$ and $\ttransfvectorb$ which can influence the correctness of the solution. Hence, special attention needs to be paid to construct this particular matrix and vector.

The construction of the solution as shown in~\eqref{eq:dg_discr_constr_ext_solve} can be done using two approaches. In the first approach, the system in~\eqref{eq:discrete_conds_ext} is solved globally, as shown in Appendix~\ref{sec:dg_solver}. In the second approach, the solution is obtained by sequential procedure along the mesh skeleton and outer boundary segments, which are presented in this section. The solution in~\eqref{eq:dg_constr_problem} guarantees that the approximate solution is both continuous and meets the boundary conditions. In order to construct the constraint components, i.e.,  $\ttransfmatrix$ and $\ttransfvectorb$, the procedure that involves integration along the mesh skeleton and Dirichlet boundary is utilized. It is performed sequentially for each segment of the mesh using the procedure described in Section~\ref{ssec:iterative}. We first present the part of the algorithm to satisfy continuity on the mesh skeleton and then for the boundary conditions.

This procedure is performed for each skeleton segment to obtain the global final solution. Each skeleton segment has exactly two adjacent elements. We want the jump in the approximation on each skeleton segment to be zero. Consider the $k$-th skeleton segment with $i$ and $j$ adjustment elements, where the $i$-th element is on the '+' side and $j$-th element is on the '--'  side of the segment. The approximation of the discontinuity jump on the $k$-th segment is:
\begin{align}\label{eq:segment_approx}
\inbbr{\temperatureh} = \inbbr{\Gshapev_\porder}\dftemp = 
\begin{bmatrix}	\Gshapev^i_{\porder} & - \Gshapev^j_\porder\end{bmatrix}
\begin{bmatrix} \dftemp^i \\ \dftemp^j \end{bmatrix}=
\Gshapesegment^k_{\porder} \dftempsegment^k
\quad\text{on } \skeletonfem^k \eqcomma
\end{align}
where $\Gshapev^i_{\porder}$ is the vector of basis function for the $i$-th element, $\dftemp^i$ is the dof vector for the $i$-th element, $\Gshapesegment^k$ and $\dftempsegment^k$ are the relevant quantities for the $k$-th segment. The vector $\dftempsegment^k$ has to satisfy the following equation:
\begin{align}\label{eq:segm_jump_eq}
\segmentmatrix^k \dftempsegment^k = \vm 0
,\quad
\segmentmatrix^k = \mint_{\Boundary_\skeletonsign^k} {\Gshapesegment^k_{\bporder}} \tran \Gshapesegment^k_{\porder} \dGa .
\end{align}

The matrix $\segmentmatrix^k$ is rectangular and additionally is singular. The dependent rows are removed from the matrix using the procedure described in Appendix~\ref{sec:dg_solver}, which yields the  $\segmentmatrixred^k$ matrix. The current version of the matrix $\ttransfmatrix$ and the current version of the vector of free variables $\freevar$, obtained from the previous skeleton segment, can be used for expressing the vector $\dftemp$. The vector $\dftempsegment^k$ is the part of the vector $\dftemp$, as shown in~\eqref{eq:segment_approx} and so the vector $\dftempsegment^k$ can be expressed with the same vector $\freevar$ as the $\dftemp$:
\begin{align}\label{eq:segm_jump_constr}
\dftemp = \ttransfmatrix \freevar \quad \longrightarrow \quad \dftempsegment^k = \ttransfmatrix^k\freevar 
\eqdot
\end{align} 
Applying the relation in \eqref{eq:segm_jump_constr} to the \eqref{eq:segm_jump_eq}, we have
\begin{align}\label{eq:seg_redexp_eq}
\segmentmatrixred^k \ttransfmatrix^k \freevar = \vm 0 \eqdot
\end{align}
The solution of~\eqref{eq:seg_redexp_eq} results in the current version of the $\freevar$ vector, which can be written in the  recurrence form 
\begin{align}\label{eq:seg_redexp_eq_sol}
\freevar = \overline{\transfmatrix}^k\freevar \eqdot
\end{align}
In order to obtain this solution, the pivoting variables in the Gauss-Jordan procedure have to be chosen in the solution procedure. For the best choice for the pivots, $\ttransfmatrix$ is analyzed, in such a way that the number of nonzero components in each column of this matrix is counted. Then the pivoting variables are chosen with the minimum nonzero components in the column of the $\ttransfmatrix$ and which also has the leading component in~\eqref{eq:seg_redexp_eq}. When the solution~\eqref{eq:seg_redexp_eq_sol} is obtained the matrix $\ttransfmatrix$ is updated:
\begin{align}
\ttransfmatrix = \ttransfmatrix\overline{\transfmatrix}^k \eqdot
\end{align}

The procedure is repeated for all the skeleton segments, and finally, we obtain the $\ttransfmatrix$ that defines the $\vectorspacecz$ set. The procedure for the boundary conditions is quite similar, remembering that there is only a single element connected to a single outer boundary segment. Finally, at the end of the procedure the matrix $\ttransfmatrix$ and vector $\ttransfvectorb$ are calculated that are applied in the approach presented in~\eqref{eq:dg_constr_problem}.

\section{\PFDG{} for linear elasticity}\label{sec:pfdg_elasticity}
The analysis of the mechanical problem starts with the standard equilibrium equations (momentum balance) in the domain $\Domain$ and the appropriate boundary conditions on the outer boundary $\Boundary$:
\begin{subequations}\label{eq:momentum}
\begin{align}
	\diverg\sig +\smas &= \vm 0  \quad \text{in } \Domain,
	\\
	\sig \cdot \Vnormal &= \boundarysign{\spow} \quad \text{on } \Boundaryneum ,
	\\
	\displu &= \boundarysign{\displu}  \quad \text{on } \Boundarydirich \eqcomma
\end{align}
\end{subequations}
where $\sig$ is the Cauchy stress tensor, $\smas$ is the body force vector per unit volume, $\boundarysign{\spow}$ is the prescribed traction vector, $\displu$ is the displacement field and $\boundarysign{\displu}$ is the prescribed displacement boundary data on $\Boundarydirich$, and $\Boundaryneum$ and $\Boundarydirich$ are the parts of the outer boundary on which the Neumman and essential boundary conditions are prescribed, respectively.

Equation~\myeqref{eq:momentum} is supplemented with the generalized Hooke's law:
\begin{align}\label{eq:hooke_full}
\sigma_{ij} = \youngsk_{ijkl}\,\strainsk_{kl} \eqcomma
\end{align}
where $\strain$ is the small-strain tensor and $\young$ is the fourth-order material moduli tensor. We adopt small-strain kinematics:
\begin{align}
\strain = \gradoper_s \displu 
= \frac{1}{2}\Bigl[\gradoper \displu + (\gradoper \displu)\tran \Bigr] \quad \text{in } \Domain \eqdot
\end{align}

In the PF-DG method, the displacement field is approximated using the same formula as in \eqref{eq:glob_approx}
\begin{align}
\displu \approx \Gshapev \dfu \eqcomma
\end{align}
where now $\Gshapev$ consists of three rows in 3D and two rows in the 2D case.

To find the approximate solution of the elasticity problem the system of equations in the form presented in \eqref{eq:dg_constr_problem} have to be solved, where in this case the definitions of the stiffness matrix and the loading vector are as follows:
\begin{align}
&	\stiffmatrix = \omint \gradoper_{sh}\Gshapev\tran : \young :  \gradoper_{sh}\Gshapev \dOm
+  \gaintt{\skeletonfem} \inbbr{\Gshapev}\tran \scalarprod \young_n : \gradoper_{sh}\Gshapev \dGa
- 	\gaintt{\Boundarydirich} \Gshapev\tran \scalarprod \young_n : \gradoper_{sh}\Gshapev \dGa \eqcomma \label{eq:K}
\\
&\rightvector = \omint \Gshapev\tran \smas \dOm + \gaintt{\Boundaryneum} \Gshapev\tran\boundarysign{\spow} \dGa \eqcomma
\end{align}
where $\young_n = \young \scalarprod \Vnormal$. The stiffness matrix $\stiffmatrix$ is nonsymmetric due to the last two components on the mesh skeleton and Dirichlet boundary that appear in~\eqref{eq:K}.

\section{\PFDG{} for biharmonic problem}\label{sec:pfdg_4order}

In this section, the DG method is developed for the fourth-order biharmonic equation. The strong form of the model problem is:
\begin{subequations}\label{eq:4problem}
\begin{align}
	\label{eq:4problem_a}
	\Delta (\paramone \Delta \temperature) + \nablav\tran (\paramtwo \nablav \temperature) + \paramthree \temperature &= f
	\quad \text{in } \Domain,
	\\
	\label{eq:4problem_b}
	\temperature &= \boundarysign{\temperature}  \quad \text{on } \Boundary,
	\\
	\label{eq:4problem_c}
	\nablav \temperature \scalarprod \Vnormal  &= \boundarysign{q}  \quad \text{on } \Boundary \eqcomma
\end{align}
\end{subequations}
where $\paramone$, $\paramtwo$, $\paramthree$, are parameter functions. The boundary conditions in \myeqref{eq:4problem_b} are of Dirichlet type. Other possible type of boundary conditions are connected with the second and third derivatives on the boundary, which are of Neumann type.  Dirichlet type boundary conditions are particularly challenging in DG methods. 

In the weak form, the Laplace operator $\Delta$ and the nabla operator $\nablav$ are viewed to be element-wise and for the sake of clarity the subscript 'h' is omitted. The construction of the weak form starts by multiplying~\eqref{eq:4problem_a} with a test function and integrating over the domain:
\begin{align}\label{eq:4weak1}
\omint \fvt\, \Delta (\paramone \Delta \temperature) \dOm
+\omint \fvt\, \nablav\tran (\paramtwo \nablav \temperature) \dOm
+ \omint\fvt\ \paramthree \temperature\dOm  - \omint \fvt f = 0 \eqdot
\end{align}
The first integral in~\eqref{eq:4weak1} is integrated by parts to yield
\begin{align}\label{eq:first_integ}
\begin{aligned}
	\omint \fvt\, \Delta (\paramone \Delta \temperature) \dOm &= 
	\omint \nablav \left( \fvt \nablav (\paramone \Delta \temperature) \right) \dOm 
	-\omint \nablav \fvt \nablav (\paramone \Delta \temperature) \dOm
	\\&=
	\gaint{} \fvt \nablav (\paramone \Delta \temperature) \scalarprod \Vnormal \dGa
	- \gaint{\skeletonsign} \inbbr{\fvt \nablav (\paramone \Delta \temperature)} \scalarprod \skeletonnormal \dGa
	-\omint \nablav \fvt \nablav (\paramone \Delta \temperature) \dOm \eqdot
\end{aligned}
\end{align}
The standard decomposition for the jump of two function is applied:
\begin{align}
\inbbr{\fvt \nablav (\paramone \Delta \temperature)} \scalarprod \skeletonnormal
= \inbbr{\fvt  } \iinbbr{\nablav(\paramone \Delta \temperature)} \scalarprod \skeletonnormal
+ \iinbbr{\fvt  } \cancelto{0}{\inbbr{\nablav(\paramone \Delta \temperature)} \scalarprod \skeletonnormal}
\qquad \text{on } \,\skeletonfem \eqdot
\end{align}
It is assumed that the jump of the gradient of Laplacian is zero. Furthermore, the last integral in \myeqref{eq:first_integ} is integrated by parts to yield
\begin{align}
\begin{aligned}
	-\omint \nablav \fvt \nablav (\paramone \Delta \temperature) \dOm &= 
	-\omint \nablav \left(  \nablav \fvt\, \paramone \Delta \temperature \right) \dOm
	+ \omint \paramone\, \Delta \fvt\,  \Delta \temperature \dOm
	\\&=
	-\gaint{} \nablav \fvt\, \paramone \Delta \temperature \scalarprod \Vnormal\dGa 
	+ \gaint{\skeletonsign} \inbbr{\nablav \fvt\, \paramone \Delta \temperature} \scalarprod \skeletonnormal \dGa
	+ \omint \paramone \, \Delta \fvt\, \Delta \temperature \dOm \eqdot
\end{aligned}
\end{align}
After the standard decomposition, it is assumed (as in the finite element method) that the jump of the Laplacian on the mesh skeleton is set to zero:
\begin{align}
\inbbr{\nablav \fvt\, \paramone \Delta \temperature} \scalarprod \skeletonnormal
= \inbbr{\nablav \fvt} \iinbbr{\paramone \Delta \temperature}\scalarprod \skeletonnormal
+ \iinbbr{\nablav \fvt} \cancelto{0}{\inbbr{\paramone \Delta \temperature}\scalarprod \skeletonnormal}
\qquad \text{on } \,\skeletonfem \eqdot
\end{align}
Finally, \eqref{eq:first_integ} is
\begin{align}
\begin{aligned}
	\omint \fvt\, \Delta (\paramone \Delta \temperature) \dOm &= 
	\omint \paramone \, \Delta \fvt\, \Delta \temperature \dOm \dGa
	+ \gaint{} \fvt \nablav (\paramone \Delta \temperature) \scalarprod \Vnormal \dGa
	-\gaint{} \nablav \fvt\, \paramone \Delta \temperature \scalarprod \Vnormal\dGa 
	\\&
	- \gaint{\skeletonsign} \inbbr{\fvt}  \iinbbr{\nablav (\paramone \Delta \temperature)} \scalarprod \skeletonnormal
	+ \gaint{\skeletonsign} \inbbr{\nablav \fvt} \iinbbr{\paramone \Delta \temperature} \scalarprod \skeletonnormal \dGa
	\eqdot
\end{aligned}
\end{align}
The second integral in \myeqref{eq:4weak1} is
\begin{align}
\omint \fvt\, \nablav\tran (\paramtwo \nablav \temperature) \dOm &= 
\omint \nablav \left(\paramtwo\,\fvt \nablav \temperature \right) \dOm
- \omint \paramtwo\, \nablav \fvt \nablav \temperature  \dOm
\\&=
\gaint{} \paramtwo\,\fvt \nablav \temperature \scalarprod \Vnormal \dGa
- \gaint{\skeletonsign}  \paramtwo\, \inbbr{\fvt \nablav \temperature } \scalarprod \skeletonnormal \dGa
- \omint \paramtwo\, \nablav \fvt \nablav \temperature  \dOm  \eqdot
\end{align}
It is assumed, similar to arriving at~\eqref{eq:heatweak2}, that the normal derivative  on the mesh skeleton is zero:
\begin{align}
\inbbr{\fvt \nablav \temperature } \scalarprod \skeletonnormal = 
\inbbr{\fvt} \iinbbr{\nablav \temperature } \scalarprod \skeletonnormal
+\iinbbr{\fvt} \cancelto{0}{\inbbr{\nablav \temperature } \scalarprod \skeletonnormal}
\qquad \text{on } \,\skeletonfem  \eqdot
\end{align}

The weak form in~\myeqref{eq:4weak1} can now be written as:
\begin{align}\label{eq:4weak2}
\begin{aligned}
	&\omint \paramone \, \Delta \fvt\, \Delta \temperature \dOm \dGa
	- \omint \paramtwo\, \nablav \fvt \nablav \temperature  \dOm
	+ \omint\fvt\ \paramthree \temperature\dOm  - \omint \fvt f 
	\\&
	- \gaint{\skeletonsign} \inbbr{\fvt}  \iinbbr{\nablav (\paramone \Delta \temperature)} \scalarprod \skeletonnormal\dGa
	+ \gaint{\skeletonsign} \inbbr{\nablav \fvt} \iinbbr{\paramone \Delta \temperature} \scalarprod \skeletonnormal \dGa
	- \gaint{\skeletonsign}  \paramtwo\, \inbbr{\fvt} \iinbbr{\nablav \temperature } \scalarprod \skeletonnormal \dGa
	\\&
	+ \gaint{} \fvt \nablav (\paramone \Delta \temperature) \scalarprod \Vnormal \dGa
	-\gaint{} \nablav \fvt\, \paramone \Delta \temperature \scalarprod \Vnormal\dGa 
	+ \gaint{} \paramtwo\,\fvt \nablav \temperature \scalarprod \Vnormal \dGa  - \omint \fvt f \dOm
	= 0  \eqdot
\end{aligned}
\end{align}

In this case, the definitions of the bilinear and linear forms from \eqref{eq:heatweak2} are: 
\begin{subequations}\label{eq:4order_bil_oper}
\begin{align}
	\label{eq:4order_bil_oper_a}
	&\begin{aligned}
		\bilinearoper(\fvt,\temperature) &= 
		\omint \paramone \, \Delta \fvt\, \Delta \temperature \dOm \dGa
		- \omint \paramtwo\, \nablav \fvt \nablav \temperature  \dOm
		+ \omint\fvt\ \paramthree \temperature\dOm  - \omint \fvt f 
		\\&
		- \gaint{\skeletonsign} \inbbr{\fvt}  \iinbbr{\nablav (\paramone \Delta \temperature)} \scalarprod \skeletonnormal\dGa
		+ \gaint{\skeletonsign} \inbbr{\nablav \fvt} \iinbbr{\paramone \Delta \temperature} \scalarprod \skeletonnormal \dGa
		- \gaint{\skeletonsign}  \paramtwo\, \inbbr{\fvt} \iinbbr{\nablav \temperature } \scalarprod \skeletonnormal \dGa
		\\&
		+ \gaint{} \fvt \nablav (\paramone \Delta \temperature) \scalarprod \Vnormal \dGa
		-\gaint{} \nablav \fvt\, \paramone \Delta \temperature \scalarprod \Vnormal\dGa 
		+ \gaint{} \paramtwo\,\fvt \nablav \temperature \scalarprod \Vnormal \dGa		
	\end{aligned}
	\\&
	\label{eq:4order_bil_oper_b}
	\linearoper(\fvt) = \omint \fvt f \dOm  \eqdot
\end{align}
\end{subequations}

The definitions of matrix $\tcmtx$ and the vector $\tcmtvb$ from the constraints \eqref{eq:discrete_conds_ext} are now as follows:
\begin{subequations}
\begin{align}
	&\begin{aligned}
		\tcmtx &= \gaint{\skeletonsign} \inbbr{\Gshapev_\bporder}\tran\inbbr{\Gshapev_\porder} \dGa  +  
		\gaint{} \Gshapev_\bporder\tran\Gshapev_\porder \dGa
		\\&+
		\gaint{\skeletonsign}\elemsize^2 \inbbr{\gradoper\Gshapev_\bporder  }\tran
		\skeletonnormal \otimes \skeletonnormal\,
		\inbbr{\gradoper\Gshapev_\porder}\dGa
		+
		\gaint{}\elemsize^2\gradoper\Gshapev_\bporder  \tran
		\Vnormal \otimes \Vnormal
		\gradoper\Gshapev_\porder \dGa  \eqcomma
	\end{aligned}
	\\&
	\tcmtvb = \gaint{} \Gshapev_\bporder\tran \boundarysign{\temperature} \dGa 
	+ \gaint{} \elemsize^2\gradoper\Gshapev_\bporder  \tran \scalarprod \Vnormal\, \boundarysign{q} \dGa \eqdot
\end{align}
\end{subequations}

In the definition of the bilinear form in~\eqref{eq:4order_bil_oper_a}, the third-order component $\nablav (\paramone \Delta \temperature)$ on the mesh skeleton and the outer boundary have to be evaluated. In this application of the DG method, Chebyshev basis functions are used, for which the second- or third-order derivatives are readily computed. The component  expands to
\begin{align}
\nablav (\paramone \Delta \temperature) = 
\begin{bmatrix}
	{\paramone}_{,x} \cdot \left( u_{,xx} + u_{,yy} \right) + \paramone \cdot \left( u_{,xxx} + u_{,xyy} \right)
	\\
	{\paramone}_{,y} \cdot \left( u_{,xx} + u_{,yy} \right) + \paramone \cdot \left( u_{,xxy} + u_{,yyy} \right)
\end{bmatrix} \eqdot
\end{align}

\section{Numerical examples}\label{sec:benchmarks}
In this section, seven numerical examples affirming the accuracy and robustness of the \PFDG{} method are presented. Six problems are chosen in 2D and one example in 3D is solved over the biunit cube. Three Poisson problems are first solved. In the fourth example, a 2D elasticity problem on the L-shaped domain is considered. In the fifth example, the hyperelastic Cook's membrane problem is analyzed. The sixth example is a fourth-order biharmonic problem on the square.

In the 2D problems, polygonal meshes are used for the analyses, and in the first example, additionally, triangular and quadrilateral meshes are used. For the  3D problem,  regular hexahedral meshes are used. Numerical integration has to be performed in the numerical computations. In the case when standard finite elements are used, i.e., triangular, quadrilateral, or hexagonal shapes, the symmetric cubatures presented in the paper by  Witherden and Vincent\cite{WITHERDEN20151232} are applied. The integration over polygonal elements is performed using the scaled boundary cubature scheme proposed by Chin and Sukumar.\cite{CHIN2021113796}

In these examples, the results are mainly shown in the form of convergence plots using the relative error measures over the whole domain, along the mesh skeleton and the Dirichlet boundary: 
\begin{align}
&
\eta = \frac{\norm{\temperaturep - \temperature}_{L^2(\Domain)}}{\norm{\temperature}_{L^2(\Domain)}}\eqcomma
\quad
\eta_E = \frac{\norm{\temperaturep - \temperature}_{E(\Domain)}}{\norm{\temperature}_{E(\Domain)}}\eqcomma
\quad
\eta_s = \frac{\norm{\inbbr{\temperaturep}}_{L^2(\skeletonfem)}}{\abs{\skeletonfem}}\eqcomma
\quad
\eta_D = \frac{\norm{\temperaturep - \boundarysign{\temperature}}_{L^2(\Boundarydirich)}} {\abs{\Boundarydirich}} \eqcomma
\end{align}
where $\norm{\temperature}_{E(\Domain)} = \bilinearoper(\temperature,\temperature)$ is the  energy norm. The measures show the global error, the discontinuity level on the mesh skeleton of the approximate solution, and how well the Dirichlet boundary conditions are met, respectively.

The theoretical rates of the convergence for the approximate solution for the $L^2$ norm and energy seminorm are:
\begin{align}
\norm{\temperaturep - \temperature}_{L^2(\Domain)} \leq C(\Domain,\mu)\elemsize^{p+1} \abs{\temperature}_{H^{\porder+1}}
\,, \quad
\norm{\temperaturep - \temperature}_{E(\Domain)} \leq C(\Domain,\mu)\elemsize^{p} \abs{\temperature}_{H^{\porder+1}}
\end{align}
where the constant $C=C(\Domain,\mu)$ depends on the domain shape of the problem and the mesh regularity parameter $\mu \in (0,1]$, which is defined as
\begin{align*}
\forall \Domainelem \in \Domainfem \quad 
\mu \elemsize_e \leq 2\rho_e \leq \elemsize_e
\end{align*}
where $\elemsize_e$ is the diameter of the finite element and $\rho_e$ is the  radius of the largest sphere contained in $\Domainelem$.

\subsection{Example 1: Poisson problem with trigonometric solution} \label{ssec:poissons_dg}
Consider the benchmark elliptic problem, $- \nabla^2 u = f$ in $\Domain = (-1,1)^2$, with $f$ chosen so that the exact solution is given by the trigonometrical function:
\begin{align} \label{eq:bench_sin_exact}
u(x,y) = x^2y + \sin \left( \frac{11 \pi x}{2} \right)
\sin \left( \frac{11 \pi y}{2} \right)  \eqdot
\end{align}

\begin{figure}
\centering
\begin{subfigure}{0.55\textwidth}
	\centering
	\includegraphics[width=\textwidth]{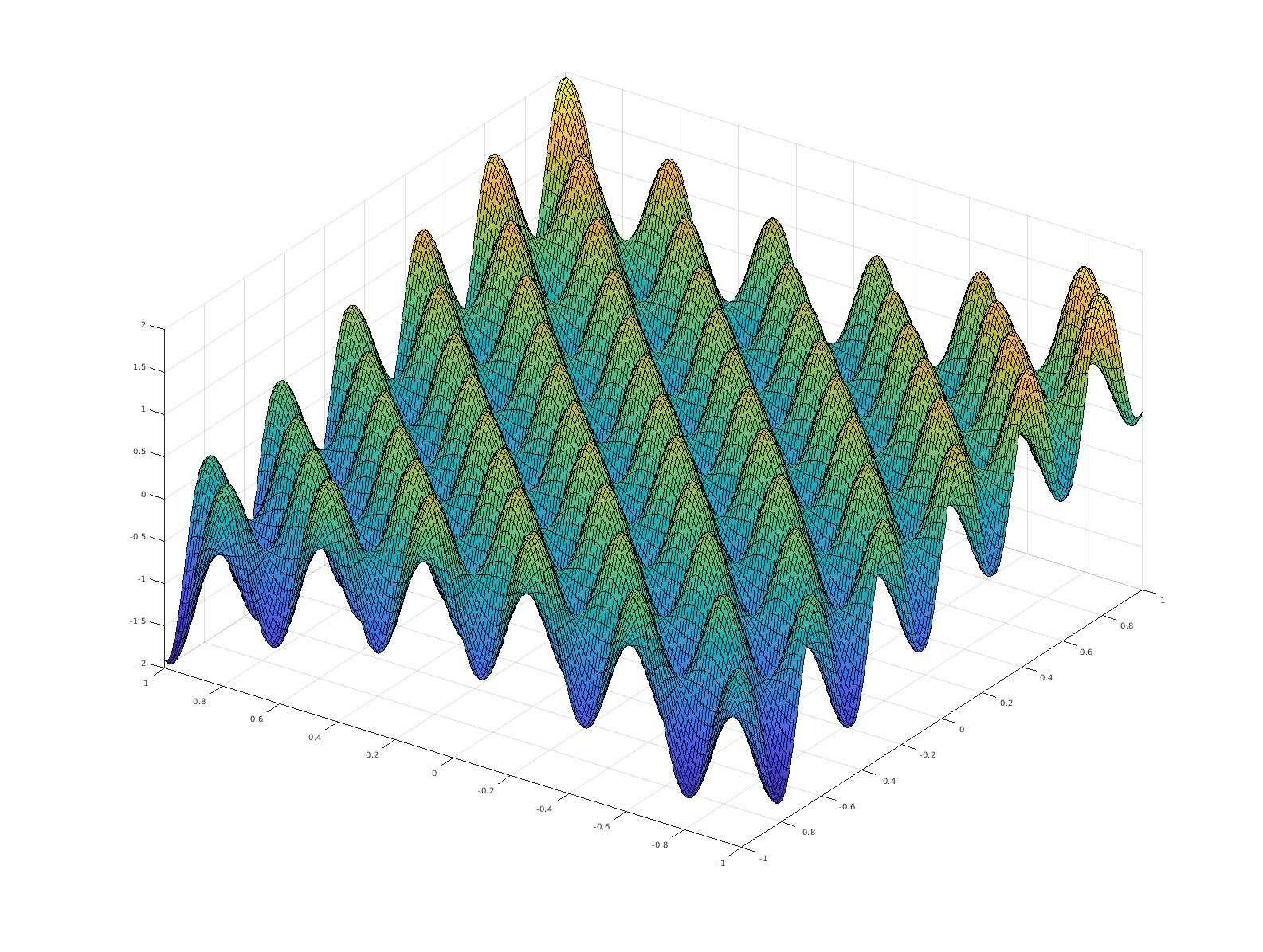}
	\subcaption{}
\end{subfigure}
\hfill
\begin{subfigure}{0.4\textwidth}
\centering
\includegraphics[width=\textwidth]{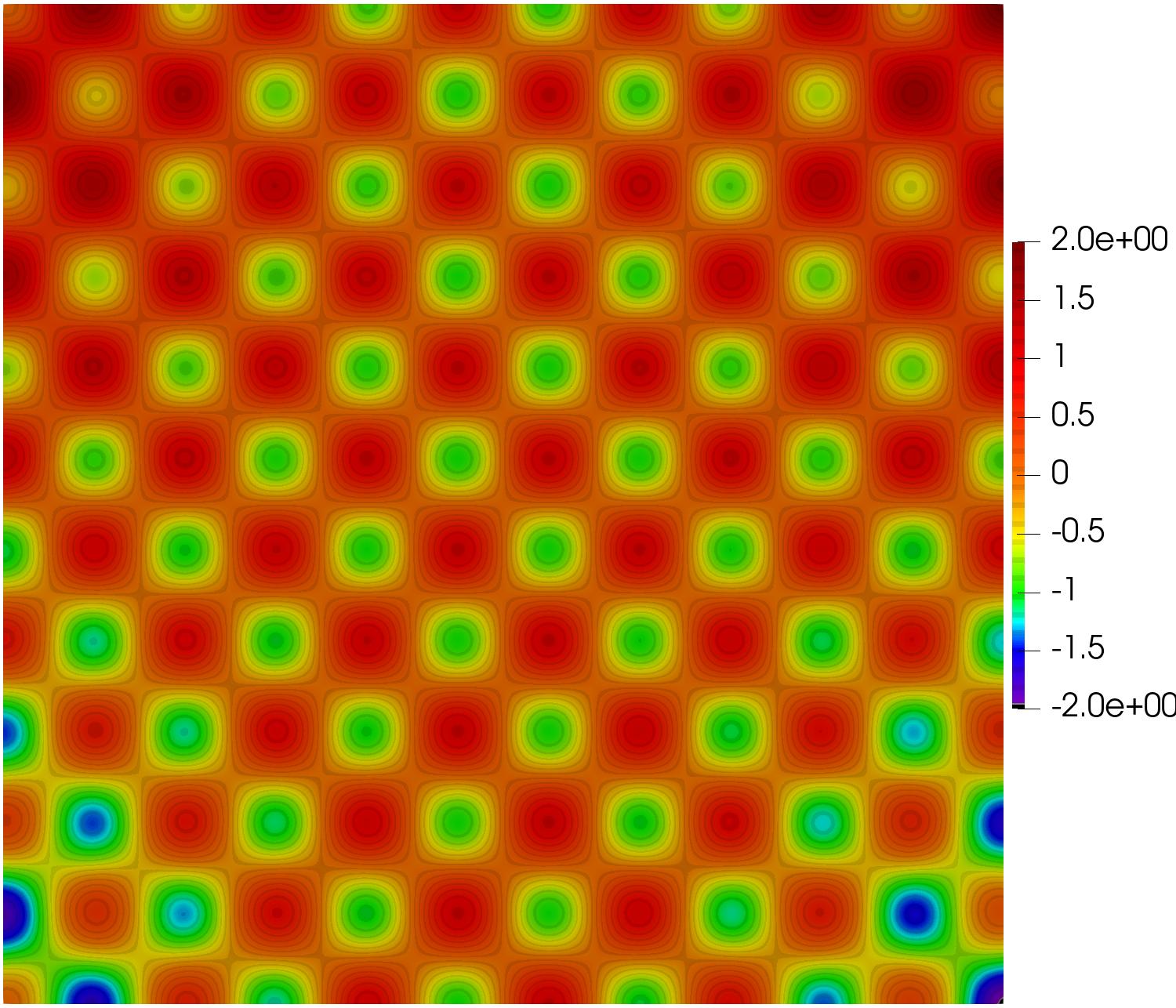}
\subcaption{}
\end{subfigure}
\caption{Surface and contour plots of the exact solution given in~\myeqref{eq:bench_sin_exact}
for Example 1.}
\label{fig:dg_benchmark_function}
\end{figure}
\begin{figure}
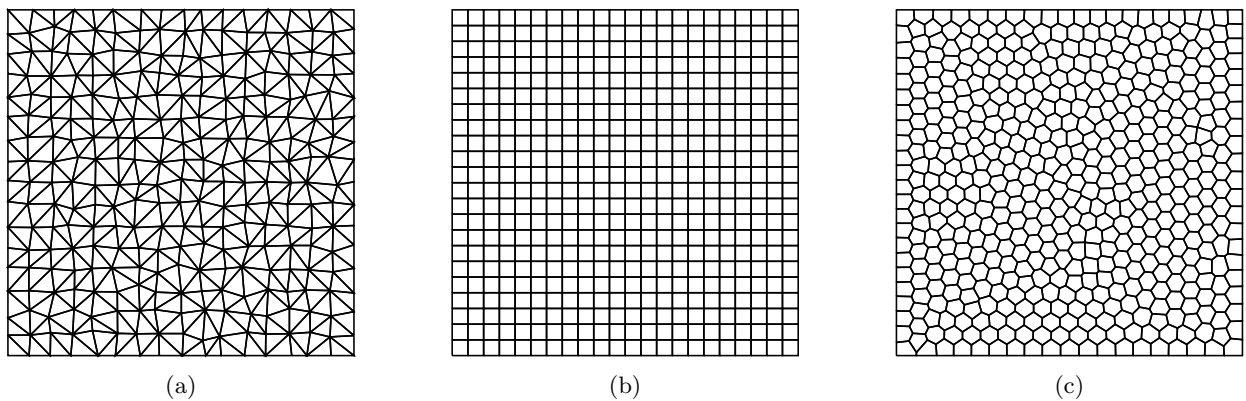
\centering
\begin{subfigure}{0.32\textwidth}\centering
\begin{tikzpicture}[scale=2.3]
	\input{figures/example1/saved_triangle_mesh_500}
\end{tikzpicture}
\subcaption{}
\end{subfigure}
\begin{subfigure}{0.32\textwidth}\centering
\begin{tikzpicture}[scale=2.3]
	\input{figures/example1/saved_square_mesh_500}
\end{tikzpicture}
\subcaption{}
\end{subfigure}
\begin{subfigure}{0.32\textwidth}\centering
\begin{tikzpicture}[scale=2.3]
	\input{figures/example1/saved_polygon_mesh_500}
\end{tikzpicture}
\subcaption{}
\end{subfigure}

\caption{Representative meshes in Example 1.
(a) Triangular mesh,
(a) quadtree mesh and
(b) polygonal mesh.}
\label{fig:ex1_repres_meshes}
\end{figure}
\begin{figure}\centering

\begin{subfigure}{0.32\textwidth} 
	\begin{subfigure}{\textwidth}
		\begin{tikzpicture}[scale=0.9]
			\begin{loglogaxis}[width=14cm, height=14cm,
				ylabel=$\eta$,
				scale=0.43,
				legend style={at={(0.03,0.03)},anchor=south west},
				xmin=10^(1.3), xmax=10^(2.3),
				ymin=1e-9, ymax=10,
				xtick = {24.372115, 39.166312, 62.928531, 101.118742, 162.480768},
				xticklabels = {594, 1534, 3960, 10225, 26400},
				ytickten={-9,-7,-5,...,1},
				]
				]
				\addplot [solid, color=blue, thick, mark=*] 
				table [x index=3, y index=4] {figures/example1/ex1_poly_p3.txt};
				\addlegendentry{\small $p=3$}
				
				\addplot [dashed, color=red, thick, mark=square*]
				table [x index=3, y index=4] {figures/example1/ex1_poly_p5.txt};
				\addlegendentry{\small $p=5$}
				
				\addplot [dashdotted, color=brown, thick, mark=triangle*, mark size=3pt]  
				table [x index=3, y index=4] {figures/example1/ex1_poly_p7.txt};
				\addlegendentry{\small $p=7$}
				
				\addplot [dotted, color=black, thick, mark=diamond*, mark size=3pt]  
				table [x index=3, y index=4] {figures/example1/ex1_poly_p10.txt};
				\addlegendentry{\small $p=10$}
				
				\logLogSlopeTTriangle{0.56}{0.15}{0.2}{11}{black};
				\logLogSlopeTriangle{0.91}{0.2}{0.69}{4}{blue};
				
			\end{loglogaxis}
		\end{tikzpicture}
		
		\begin{tikzpicture}[scale=0.9]
			\begin{loglogaxis}[width=14cm, height=14cm,
				xlabel=${\#dof}$,
				ylabel=$\eta_E$,
				scale=0.43,
				legend style={at={(0.03,0.03)},anchor=south west},
				xmin=10^(1.3), xmax=10^(2.3),
				ymin=1e-9, ymax=10,
				xtick = {24.372115, 39.166312, 62.928531, 101.118742, 162.480768},
				xticklabels = {594, 1534, 3960, 10225, 26400},
				ytickten={-9,-7,-5,...,1},
				]
				\addplot [solid, color=blue, thick, mark=*] 
				table [x index=3, y index=5] {figures/example1/ex1_poly_p3.txt};
				\addlegendentry{\small $p=3$}
				
				\addplot [dashed, color=red, thick, mark=square*]
				table [x index=3, y index=5] {figures/example1/ex1_poly_p5.txt};
				\addlegendentry{\small $p=5$}
				
				\addplot [dashdotted, color=brown, thick, mark=triangle*, mark size=3pt]  
				table [x index=3, y index=5] {figures/example1/ex1_poly_p7.txt};
				\addlegendentry{\small $p=7$}
				
				\addplot [dotted, color=black, thick, mark=diamond*, mark size=3pt]  
				table [x index=3, y index=5] {figures/example1/ex1_poly_p10.txt};
				\addlegendentry{\small $p=10$}
				
				\logLogSlopeTTriangle{0.65}{0.15}{0.2}{10}{black};
				\logLogSlopeTriangle{0.91}{0.2}{0.78}{3}{blue};
				
			\end{loglogaxis}
		\end{tikzpicture}
	\end{subfigure}
	\subcaption{}
\end{subfigure}
\hspace{0.5cm}
\begin{subfigure}{0.32\textwidth} 
		\begin{subfigure}{\textwidth}
			\begin{tikzpicture}[scale=0.9]
				\begin{loglogaxis}[width=14cm, height=14cm,
					scale=0.43,
					legend style={at={(0.03,0.03)},anchor=south west},
					xmin=10^(1.3), xmax=10^(2.3),
					ymin=1e-9, ymax=10,
					xtick = {24.372115, 39.166312, 62.928531, 101.118742, 162.480768},
					xticklabels = {594, 1534, 3960, 10225, 26400},
					ytickten={-9,-7,-5,...,1},
					]
					\addplot [solid, color=blue, thick, mark=*] 
					table [x index=3, y index=4] {figures/example1/ex1_square_p3.txt};
					\addlegendentry{\small $p=3$}
					
					\addplot [dashed, color=red, thick, mark=square*]
					table [x index=3, y index=4] {figures/example1/ex1_square_p5.txt};
					\addlegendentry{\small $p=5$}
					
					\addplot [dashdotted, color=brown, thick, mark=triangle*, mark size=3pt]  
					table [x index=3, y index=4] {figures/example1/ex1_square_p7.txt};
					\addlegendentry{\small $p=7$}
					
					\addplot [dotted, color=black, thick, mark=diamond*, mark size=3pt]  
					table [x index=3, y index=4] {figures/example1/ex1_square_p10.txt};
					\addlegendentry{\small $p=10$}
					
					\logLogSlopeTTriangle{0.61}{0.15}{0.2}{11}{black};
					\logLogSlopeTriangle{0.91}{0.2}{0.70}{4}{blue};
					
				\end{loglogaxis}
			\end{tikzpicture}
			
			\begin{tikzpicture}[scale=0.9]
				\begin{loglogaxis}[width=14cm, height=14cm,
					xlabel=${\#dof}$,
					scale=0.43,
					legend style={at={(0.03,0.03)},anchor=south west},
					xmin=10^(1.3), xmax=10^(2.3),
					ymin=1e-9, ymax=10,
					xtick = {24.372115, 39.166312, 62.928531, 101.118742, 162.480768},
					xticklabels = {594, 1534, 3960, 10225, 26400},
					ytickten={-9,-7,-5,...,1},
					]
					\addplot [solid, color=blue, thick, mark=*] 
					table [x index=3, y index=5] {figures/example1/ex1_square_p3.txt};
					\addlegendentry{\small $p=3$}
					
					\addplot [dashed, color=red, thick, mark=square*]
					table [x index=3, y index=5] {figures/example1/ex1_square_p5.txt};
					\addlegendentry{\small $p=5$}
					
					\addplot [dashdotted, color=brown, thick, mark=triangle*, mark size=3pt]  
					table [x index=3, y index=5] {figures/example1/ex1_square_p7.txt};
					\addlegendentry{\small $p=7$}
					
					\addplot [dotted, color=black, thick, mark=diamond*, mark size=3pt]  
					table [x index=3, y index=5] {figures/example1/ex1_square_p10.txt};
					\addlegendentry{\small $p=10$}
					
					\logLogSlopeTTriangle{0.70}{0.15}{0.2}{10}{black};
					\logLogSlopeTriangle{0.91}{0.2}{0.77}{3}{blue};
					
				\end{loglogaxis}
			\end{tikzpicture}
		\end{subfigure}
		\subcaption{}
	\end{subfigure}
	\begin{subfigure}{0.32\textwidth} 
			\begin{subfigure}{\textwidth}
				\begin{tikzpicture}[scale=0.9]
					\begin{loglogaxis}[width=14cm, height=14cm,
						scale=0.43,
						legend style={at={(0.03,0.03)},anchor=south west},
						xmin=10^(1.3), xmax=10^(2.3),
						ymin=1e-9, ymax=10,
						xtick = {24.372115, 39.166312, 62.928531, 101.118742, 162.480768},
						xticklabels = {594, 1534, 3960, 10225, 26400},
						ytickten={-9,-7,-5,...,1},
						]
						
						\addplot [solid, color=blue, thick, mark=*] 
						table [x index=3, y index=4] {figures/example1/ex1_triangle_p3.txt};
						\addlegendentry{\small $p=3$}
						
						\addplot [dashed, color=red, thick, mark=square*]
						table [x index=3, y index=4] {figures/example1/ex1_triangle_p5.txt};
						\addlegendentry{\small $p=5$}
						
						\addplot [dashdotted, color=brown, thick, mark=triangle*, mark size=3pt]  
						table [x index=3, y index=4] {figures/example1/ex1_triangle_p7.txt};
						\addlegendentry{\small $p=7$}
						
						\addplot [dotted, color=black, thick, mark=diamond*, mark size=3pt]  
						table [x index=3, y index=4] {figures/example1/ex1_triangle_p10.txt};
						\addlegendentry{\small $p=10$}
						
						\logLogSlopeTTriangle{0.72}{0.15}{0.2}{11}{black};
						\logLogSlopeTriangle{0.91}{0.2}{0.72}{4}{blue};
						
					\end{loglogaxis}
				\end{tikzpicture}
				
				\begin{tikzpicture}[scale=0.9]
					\begin{loglogaxis}[width=14cm, height=14cm,
						xlabel=${\#dof}$,
						scale=0.43,
						legend style={at={(0.03,0.03)},anchor=south west},
						xmin=10^(1.3), xmax=10^(2.3),
						ymin=1e-9, ymax=10,
						xtick = {24.372115, 39.166312, 62.928531, 101.118742, 162.480768},
						xticklabels = {594, 1534, 3960, 10225, 26400},
						ytickten={-9,-7,-5,...,1},
						]
						\addplot [solid, color=blue, thick, mark=*] 
						table [x index=3, y index=5] {figures/example1/ex1_triangle_p3.txt};
						\addlegendentry{\small $p=3$}
						
						\addplot [dashed, color=red, thick, mark=square*]
						table [x index=3, y index=5] {figures/example1/ex1_triangle_p5.txt};
						\addlegendentry{\small $p=5$}
						
						\addplot [dashdotted, color=brown, thick, mark=triangle*, mark size=3pt]  
						table [x index=3, y index=5] {figures/example1/ex1_triangle_p7.txt};
						\addlegendentry{\small $p=7$}
						
						\addplot [dotted, color=black, thick, mark=diamond*, mark size=3pt]  
						table [x index=3, y index=5] {figures/example1/ex1_triangle_p10.txt};
						\addlegendentry{\small $p=10$}
						
						\logLogSlopeTTriangle{0.71}{0.15}{0.29}{10}{black};
						\logLogSlopeTriangle{0.91}{0.2}{0.79}{3}{blue};
						
					\end{loglogaxis}
				\end{tikzpicture}
			\end{subfigure}
			\subcaption{}
		\end{subfigure}
		\caption{Convergence of the \PFDG{} method in the $L^2$ 
			norm (top) and the energy seminorm (bottom) for Example 1. 
			Computations are shown for
			(a) polygonal, (b) quadrilateral and (c) triangular meshes. Scale is
			$\frac{1}{2}\log$--$\log$.}
		\label{fig:example1_conv}
	\end{figure}
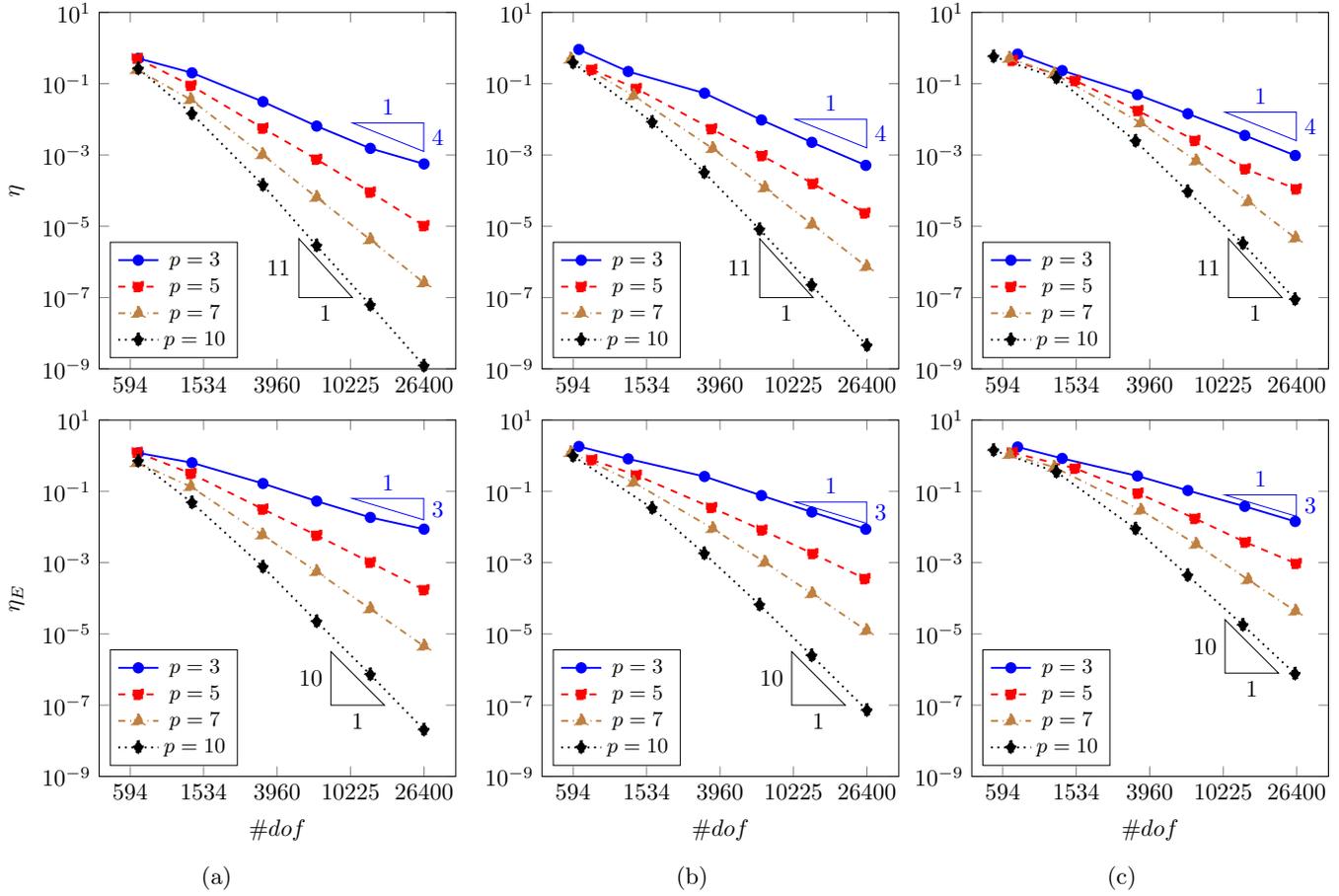

The exact solution has many hills and valleys that are regularly distributed in the domain (see Fig.~\ref{fig:dg_benchmark_function}). In this example, Dirichlet boundary conditions are applied on the whole outer boundary. This problem is solved using triangular, quadrilateral, and polygonal meshes, and representative meshes are shown in Fig.~\ref{fig:ex1_repres_meshes}. 	The calculations are performed for orders $p=3$ to $p=10$ on a sequence of refined meshes. The convergence curves in the $L^2$ and the energy seminorm are presented in Fig.~\ref{fig:example1_conv}. Optimal theoretical convergence rates are obtained by the PF-DG method on all meshes, with polygonal meshes yielding the smallest errors.

\begin{figure}\centering
		\begin{subfigure}{0.32\textwidth}
			\begin{tikzpicture}[scale=0.9]
				\begin{loglogaxis}[width=14cm, height=14cm,
					xlabel=${\#dof}$,
					ylabel=$\eta$,
					scale=0.43,
					legend style={at={(0.03,0.03)},anchor=south west},
					xtick = {24.372115, 39.166312, 62.928531, 101.118742, 162.480768},
					xticklabels = {594, 1534, 3960, 10225, 26400},
					]
					
					\addplot [solid, color=blue, thick, mark=*] 
					table [x index=3, y index=4] {figures/example1/ex1_poly_p3.txt};
					\addlegendentry{\small PF-DG}
					
					\addplot [dashed, color=red, thick, mark=square*] 
					table [x index=3, y index=10] {figures/example1/ex1_poly_p3.txt};
					\addlegendentry{\small DGFD}
					
					\addplot [dashdotted, color=brown, thick, mark=triangle*, mark size=3pt]  
					table [x index=3, y index=8] {figures/example1/ex1_poly_p3.txt};
					\addlegendentry{\small DG\cite{ODEN1998491}}
					
					
				\end{loglogaxis}
			\end{tikzpicture}
			\subcaption{}
		\end{subfigure}
		\hspace{0.4cm}
		\begin{subfigure}{0.32\textwidth}\centering
			\begin{tikzpicture}[scale=0.9]
				\begin{loglogaxis}[width=14cm, height=14cm,
					xlabel=${\#dof}$,
					scale=0.43,
					legend style={at={(0.03,0.03)},anchor=south west},
					xtick = {24.372115, 39.166312, 62.928531, 101.118742, 162.480768},
					xticklabels = {594, 1534, 3960, 10225, 26400},
					]
					
					\addplot [solid, color=blue, thick, mark=*] 
					table [x index=3, y index=4] {figures/example1/ex1_poly_p5.txt};
					\addlegendentry{\small PF-DG}
					
					\addplot [dashed, color=red, thick, mark=square*] 
					table [x index=3, y index=10] {figures/example1/ex1_poly_p5.txt};
					\addlegendentry{\small DGFD}
					
					\addplot [dashdotted, color=brown, thick, mark=triangle*, mark size=3pt]  
					table [x index=3, y index=8] {figures/example1/ex1_poly_p5.txt};
					\addlegendentry{\small DG\cite{ODEN1998491}}

				\end{loglogaxis}
			\end{tikzpicture}
			\subcaption{}
		\end{subfigure}
		\begin{subfigure}{0.32\textwidth}
			\begin{tikzpicture}[scale=0.9]
				\begin{loglogaxis}[width=14cm, height=14cm,
					xlabel=${\#dof}$,
					scale=0.43,
					legend style={at={(0.03,0.03)},anchor=south west},
					xtick = {24.372115, 39.166312, 62.928531, 101.118742, 162.480768},
					xticklabels = {594, 1534, 3960, 10225, 26400},
					]

					\addplot [solid, color=blue, thick, mark=*] 
					table [x index=3, y index=4] {figures/example1/ex1_poly_p10.txt};
					\addlegendentry{\small PF-DG}
					
					\addplot [dashed, color=red, thick, mark=square*] 
					table [x index=3, y index=10] {figures/example1/ex1_poly_p10.txt};
					\addlegendentry{\small DGFD}
					
					\addplot [dashdotted, color=brown, thick, mark=triangle*, mark size=3pt]  
					table [x index=3, y index=8] {figures/example1/ex1_poly_p10.txt};
					\addlegendentry{\small DG\cite{ODEN1998491}}
					
				\end{loglogaxis}
			\end{tikzpicture}
			\subcaption{}
		\end{subfigure}
		
		\caption{Convergence study for Example 1 on the polygonal meshes using
			the three DG methods 
			(DG,\cite{ODEN1998491} DGFD and \PFDG). 
			The plots are presented for (a) $\porder=3$, 
			(b) $\porder=5$ and (c) $\porder=10$.
			Scale is $\frac{1}{2}\log$--$\log$.
		}
		\label{fig:ex1_paths_methods}
	\end{figure}
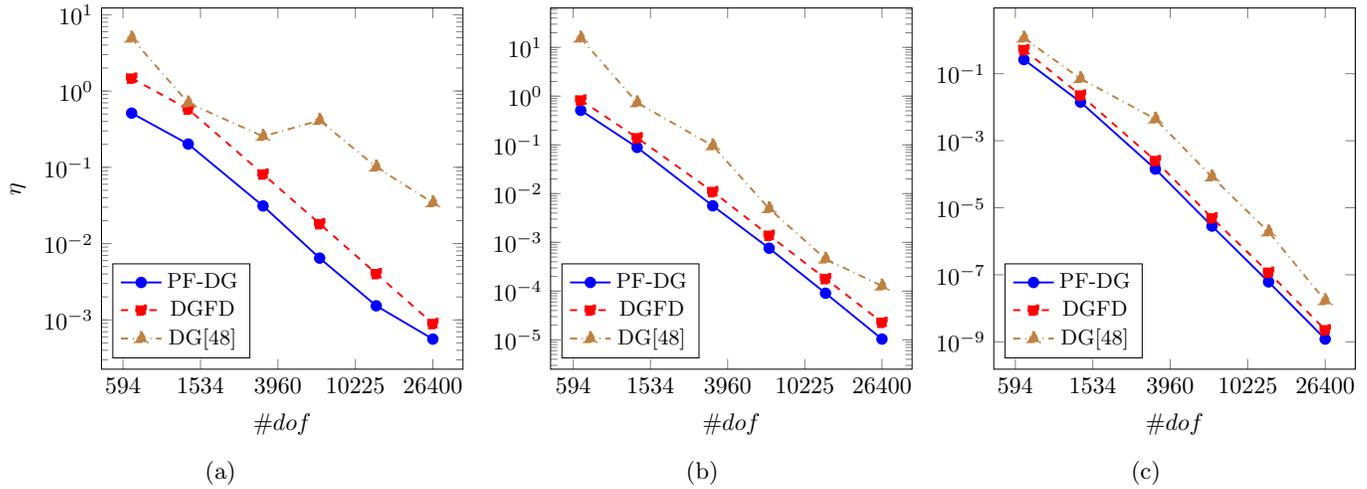

	This example is solved using the \pfdg{} method, and for comparison, using two other DG methods---DGFD method by Ja{\'s}kowiec\cite{jaskowiec_cm_2017} and the version of the penalty-free DG method proposed by Oden at al.\cite{ODEN1998491} In Fig.~\ref{fig:ex1_paths_methods}, the convergence curves of three DG methods are compared to each other: DG method of Oden et al.,~\cite{ODEN1998491} DGFD and the \pfdg{} presented in this work. The DGFD and the \pfdg{} have the same convergence rate and the level of accuracy for these two methods is similar. Significantly	worse results in each case are obtained using the DG\cite{ODEN1998491}	method, which is appreciable for $\porder=3$, where additionally, the convergence is nonmonotonic.
	
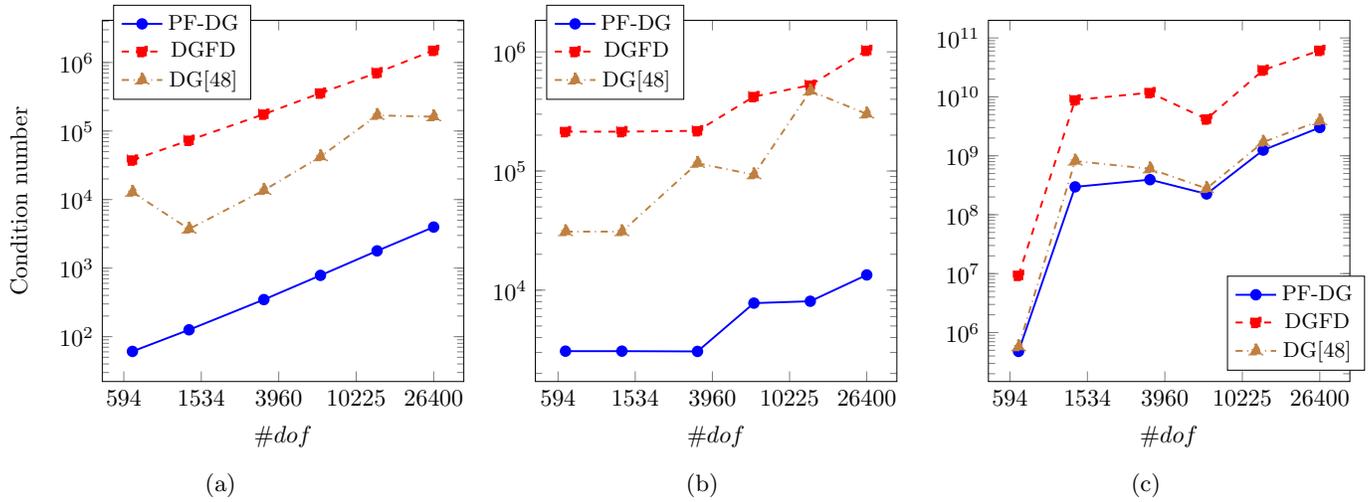
\begin{figure}
		\begin{subfigure}{0.32\textwidth}
			\begin{tikzpicture}[scale=0.9]
				\begin{loglogaxis}[width=14cm, height=14cm,
					xlabel=${\#dof}$,
					ylabel=Condition number,
					scale=0.43,
					legend style={at={(0.03,0.78)},anchor=south west},
					xtick = {24.372115, 39.166312, 62.928531, 101.118742, 162.480768},
					xticklabels = {594, 1534, 3960, 10225, 26400},
					]
					
					\addplot [solid, color=blue, thick, mark=*] 
					table [x index=3, y index=4] {figures/example1/ex1_sin_poly_conds_p3.txt};
					\addlegendentry{\small PF-DG}
					
					\addplot [dashed, color=red, thick, mark=square*]
					table [x index=3, y index=5] {figures/example1/ex1_sin_poly_conds_p3.txt};
					\addlegendentry{\small DGFD}	
					
					\addplot [dashdotted, color=brown, thick, mark=triangle*, mark size=3pt] 
					table [x index=3, y index=6] {figures/example1/ex1_sin_poly_conds_p3.txt};
					\addlegendentry{\small DG\cite{ODEN1998491}}			
				\end{loglogaxis}
			\end{tikzpicture}
			\subcaption{}
		\end{subfigure}
		\hspace{0.4cm}
		\begin{subfigure}{0.32\textwidth}
			\begin{tikzpicture}[scale=0.9]
				\begin{loglogaxis}[width=14cm, height=14cm,
					xlabel=${\#dof}$,
					scale=0.43,
					legend style={at={(0.03,0.78)},anchor=south west},
					xtick = {24.372115, 39.166312, 62.928531, 101.118742, 162.480768},
					xticklabels = {594, 1534, 3960, 10225, 26400},
					]
					
					\addplot [solid, color=blue, thick, mark=*] 
					table [x index=3, y index=4] {figures/example1/ex1_sin_poly_conds_p5.txt};
					\addlegendentry{\small PF-DG}
					
					\addplot [dashed, color=red, thick, mark=square*]
					table [x index=3, y index=5] {figures/example1/ex1_sin_poly_conds_p5.txt};
					\addlegendentry{\small DGFD}
					
					\addplot [dashdotted, color=brown, thick, mark=triangle*, mark size=3pt] 
					table [x index=3, y index=6] {figures/example1/ex1_sin_poly_conds_p5.txt};
					\addlegendentry{\small DG\cite{ODEN1998491}}				
				\end{loglogaxis}
			\end{tikzpicture}
			\subcaption{}
		\end{subfigure}
		\begin{subfigure}{0.32\textwidth}
			\begin{tikzpicture}[scale=0.9]
				\begin{loglogaxis}[width=14cm, height=14cm,
					xlabel=${\#dof}$,
					scale=0.43,
					legend style={at={(0.66,0.03)},anchor=south west},
					xtick = {24.372115, 39.166312, 62.928531, 101.118742, 162.480768},
					xticklabels = {594, 1534, 3960, 10225, 26400},
					]
					
					\addplot [solid, color=blue, thick, mark=*] 
					table [x index=3, y index=4] {figures/example1/ex1_sin_poly_conds_p10.txt};
					\addlegendentry{\small PF-DG}
					
					\addplot [dashed, color=red, thick, mark=square*]
					table [x index=3, y index=5] {figures/example1/ex1_sin_poly_conds_p10.txt};
					\addlegendentry{\small DGFD}		
					
					\addplot [dashdotted, color=brown, thick, mark=triangle*, mark size=3pt] 
					table [x index=3, y index=6] {figures/example1/ex1_sin_poly_conds_p10.txt};
					\addlegendentry{\small DG\cite{ODEN1998491}}		
					
				\end{loglogaxis}
			\end{tikzpicture}
			\subcaption{}
		\end{subfigure}

		\caption{
			Condition numbers for the matrices used in Example 1, for 
			the three DG methods 
			(DG,\cite{ODEN1998491} DGFD and \PFDG). 
			The plots are presented for (a) $\porder=3$, 
			(b) $\porder=5$ and (c) $\porder=10$.
		}
		\label{fig:ex1_conds}
	\end{figure}

	The condition number, which is the ratio of the maximum to minimum singular value of a matrix, is presented for the stiffness matrix in Fig.~\ref{fig:ex1_conds}. The condition numbers are presented for the three methods \pfdg, DGFD and DG\cite{ODEN1998491} and for approximation orders $p=3,\, 5,\, 10$. The stiffness matrices for \pfdg{} are much better conditioned in comparison to the other two methods. It confirms that the \pfdg{} method is numerically stable. 
	
\begin{figure}\centering
		\begin{subfigure}{0.32\textwidth}
			\centering
			\begin{tikzpicture}[scale=2.3]
				\input{figures/example1/saved_rand_mesh}
			\end{tikzpicture}
			\subcaption{}
		\end{subfigure}
		\begin{subfigure}{0.32\textwidth}\centering
			\includegraphics[height=0.8\textwidth]{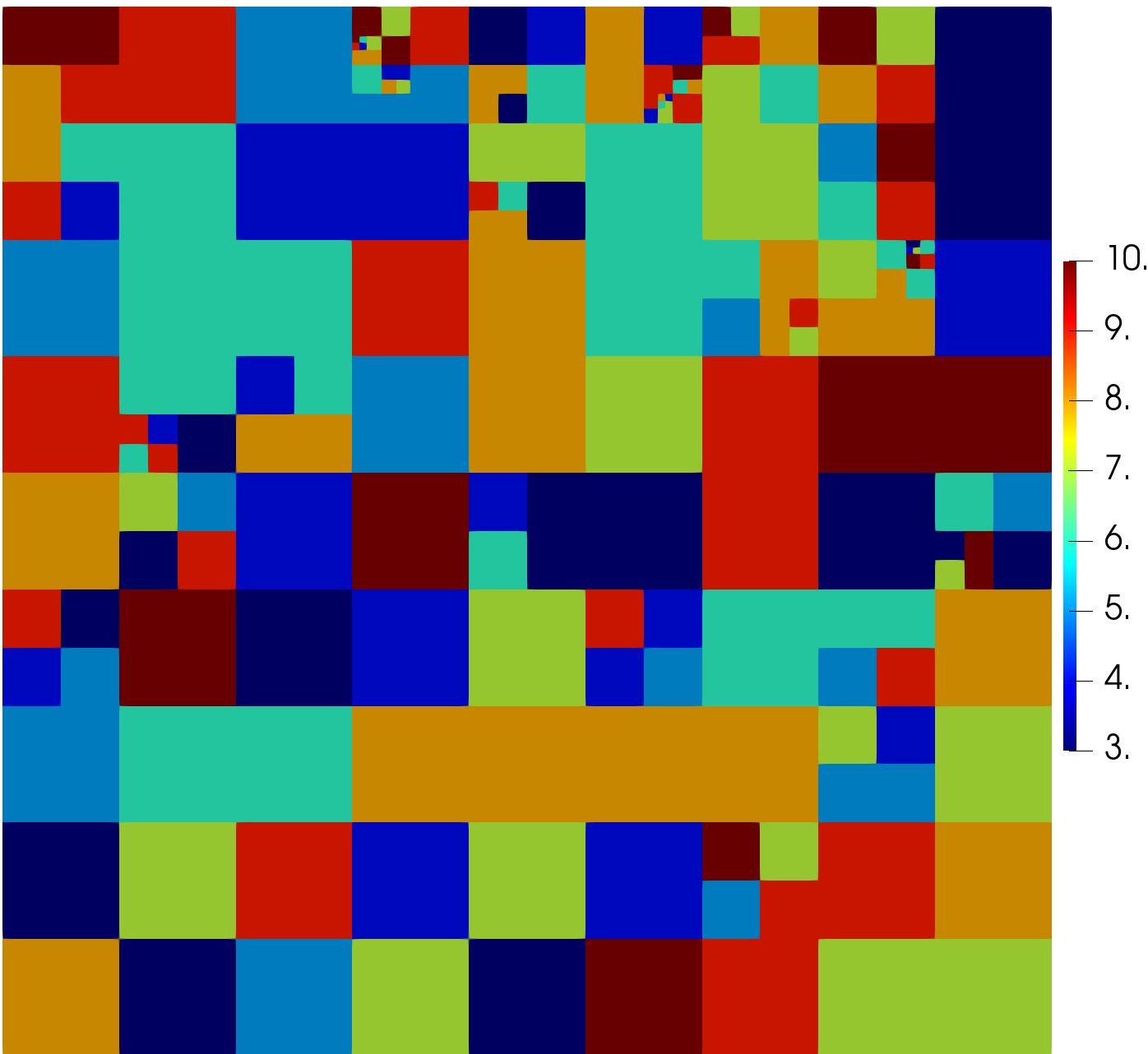}
			\subcaption{}
		\end{subfigure}
		\caption{Randomly refined quadrilateral mesh for Example 1.
			(a) Quadtree mesh with hanging nodes and
			(b) approximation order in each element.}
		\label{fig:ex1_random_mesh}
	\end{figure}

\begin{figure}\centering
		\begin{subfigure}{0.32\textwidth}
			\includegraphics[height=0.8\textwidth]{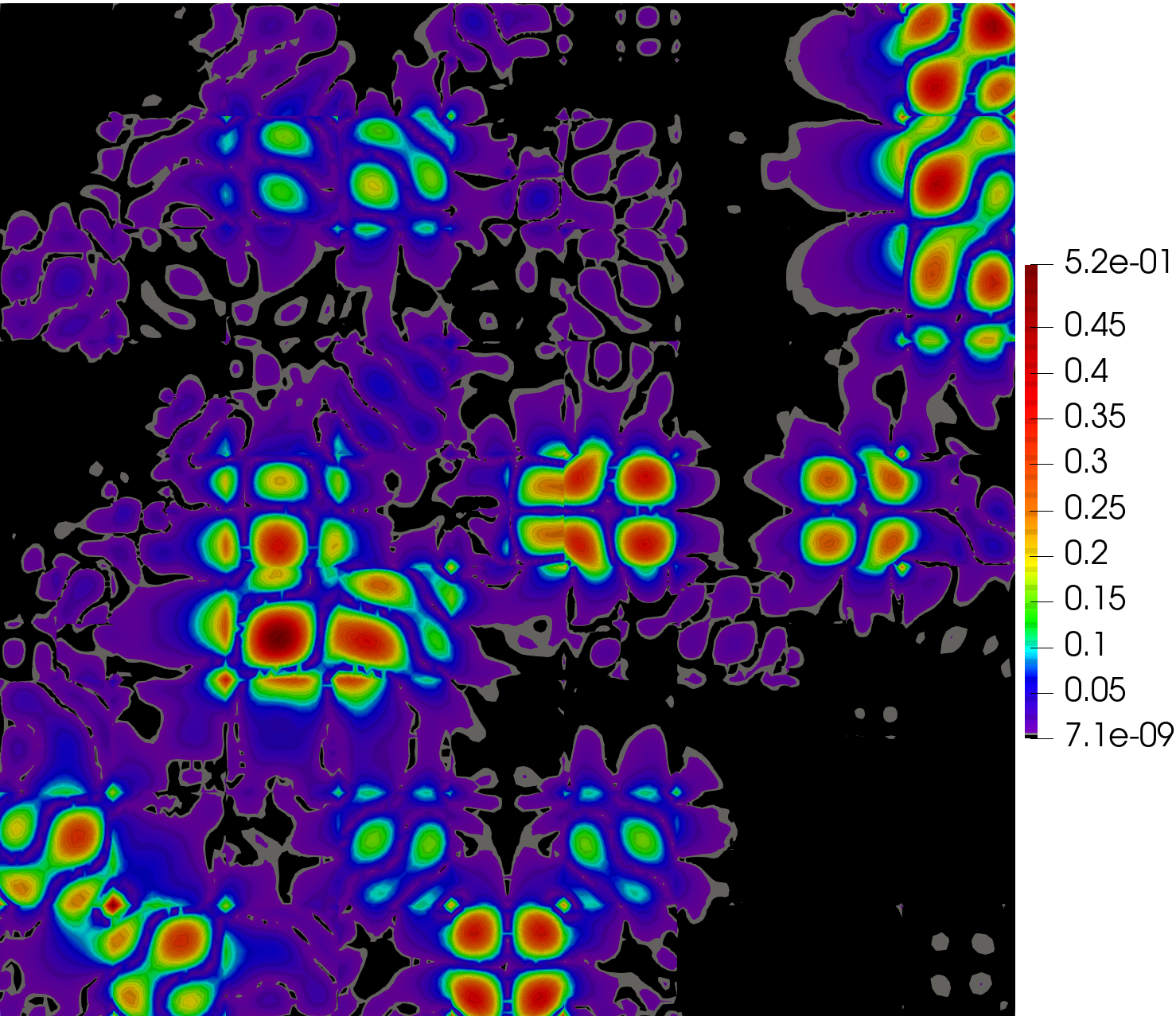}
			\subcaption{}
		\end{subfigure}
		\begin{subfigure}{0.32\textwidth}
			\includegraphics[height=0.8\textwidth]{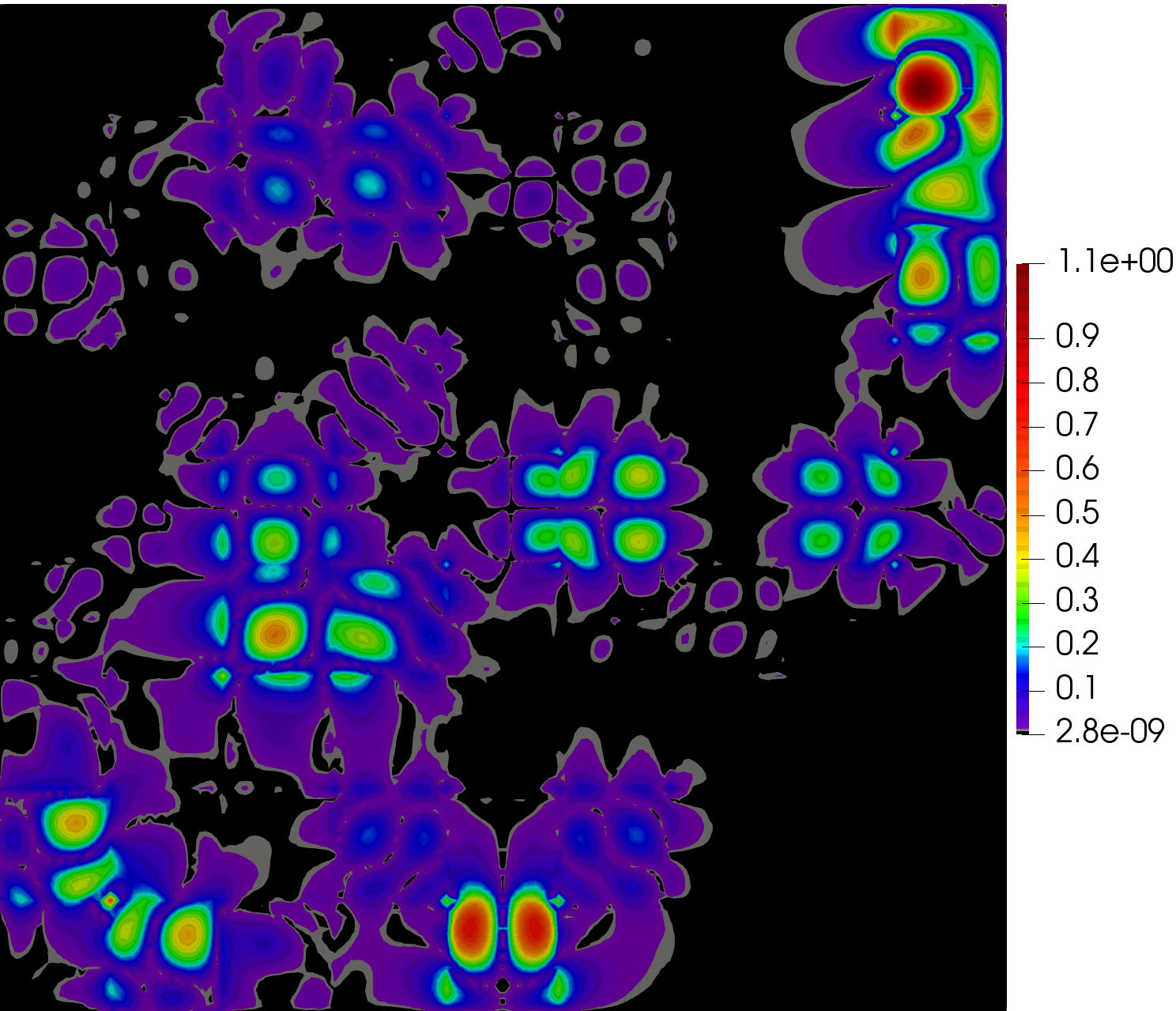}
			\subcaption{}
		\end{subfigure}
		\begin{subfigure}{0.32\textwidth}
			\includegraphics[height=0.8\textwidth]{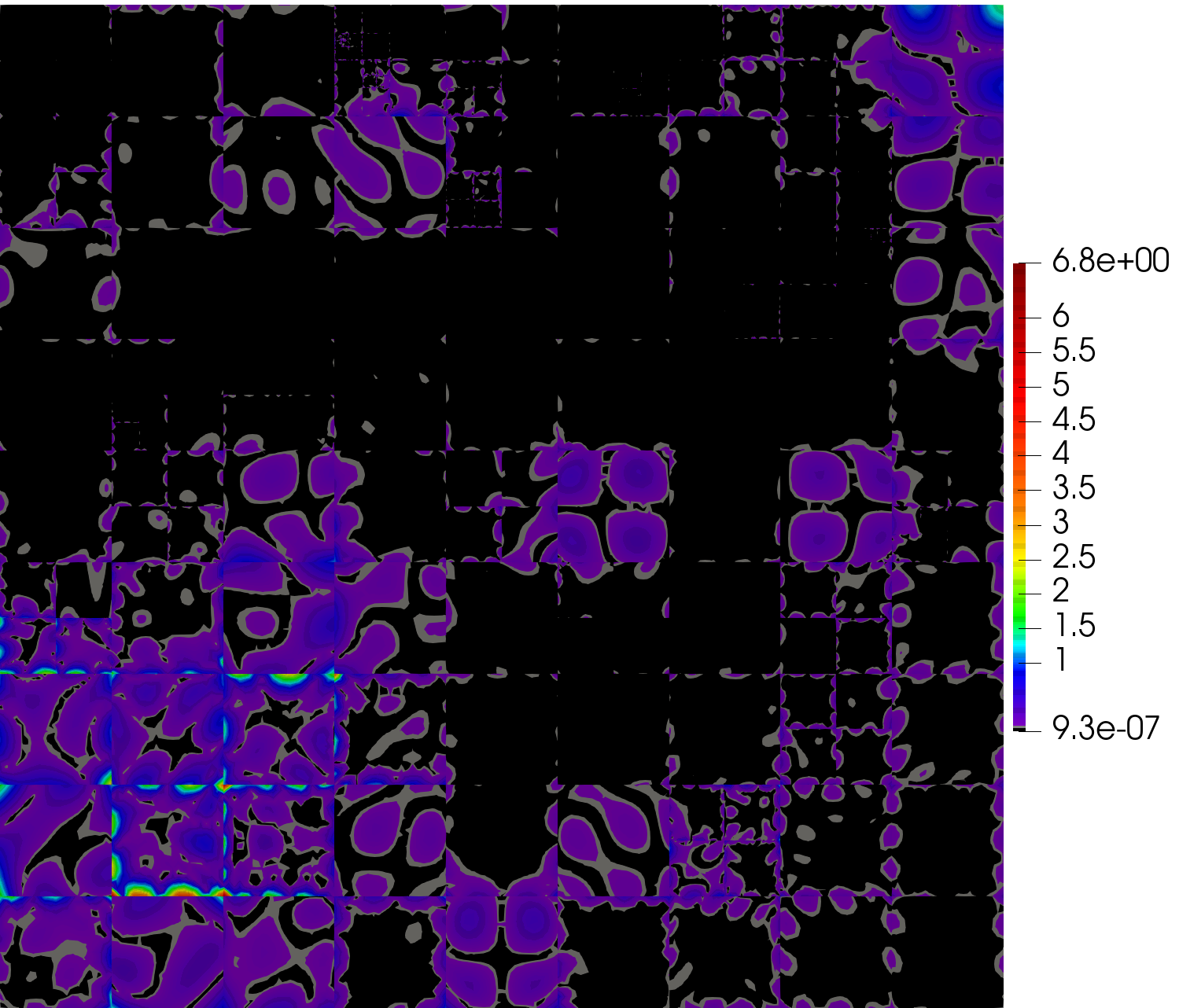}
			\subcaption{}
		\end{subfigure}
		
		\caption{Error distribution on the randomly refined quadrilateral
			mesh for Example 1. (a) \PFDG, (b) DGFD 
			and (c) DG method.\cite{ODEN1998491}
		}
		\label{fig:ex1_random_errors}
	\end{figure}

	The same example is solved on the quadrilateral mesh that is randomly refined and the approximation orders in the elements $\porder\in[3,10]$ are set randomly (see Fig.~\ref{fig:ex1_random_mesh}). The error distributions for the three versions of the DG method are shown in Fig.~\ref{fig:ex1_random_errors}. The errors for the \PFDG{} and the DGFD are comparable to each other. Much larger errors appear in the DG method of Oden et al.,~\cite{ODEN1998491} which reveals that it is not well-suited for such meshes.

	\subsection{Example 2: Poisson problem with exponential solution} \label{ssec:poissons_dg_exp}
\begin{figure}\centering
		\begin{subfigure}{0.48\textwidth}
			\centering
			\includegraphics[width=0.96\textwidth]{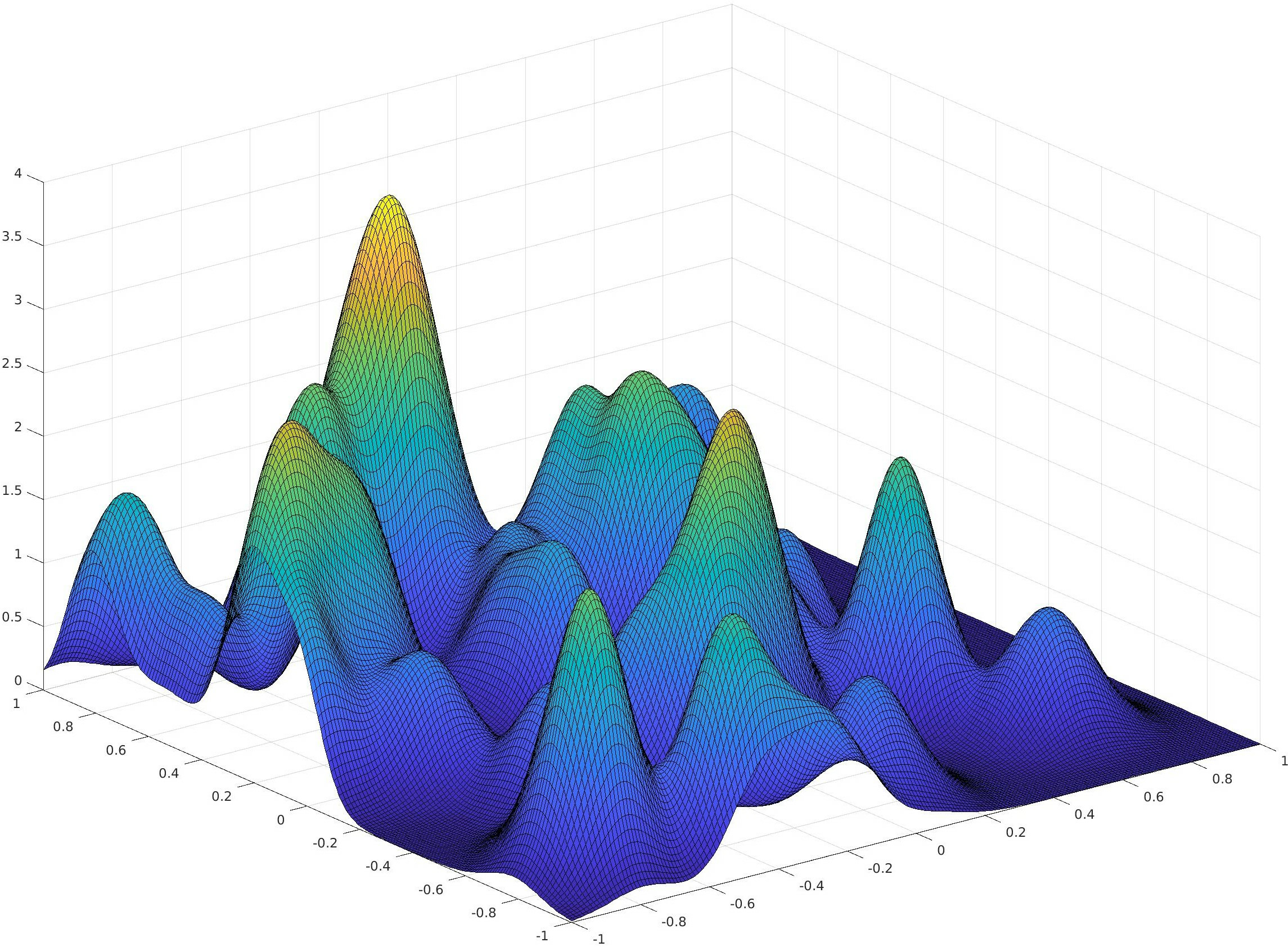}
			\subcaption{}
		\end{subfigure}
		\hfill
		\begin{subfigure}{0.48\textwidth}
			\centering
			\includegraphics[width=0.74\textwidth]{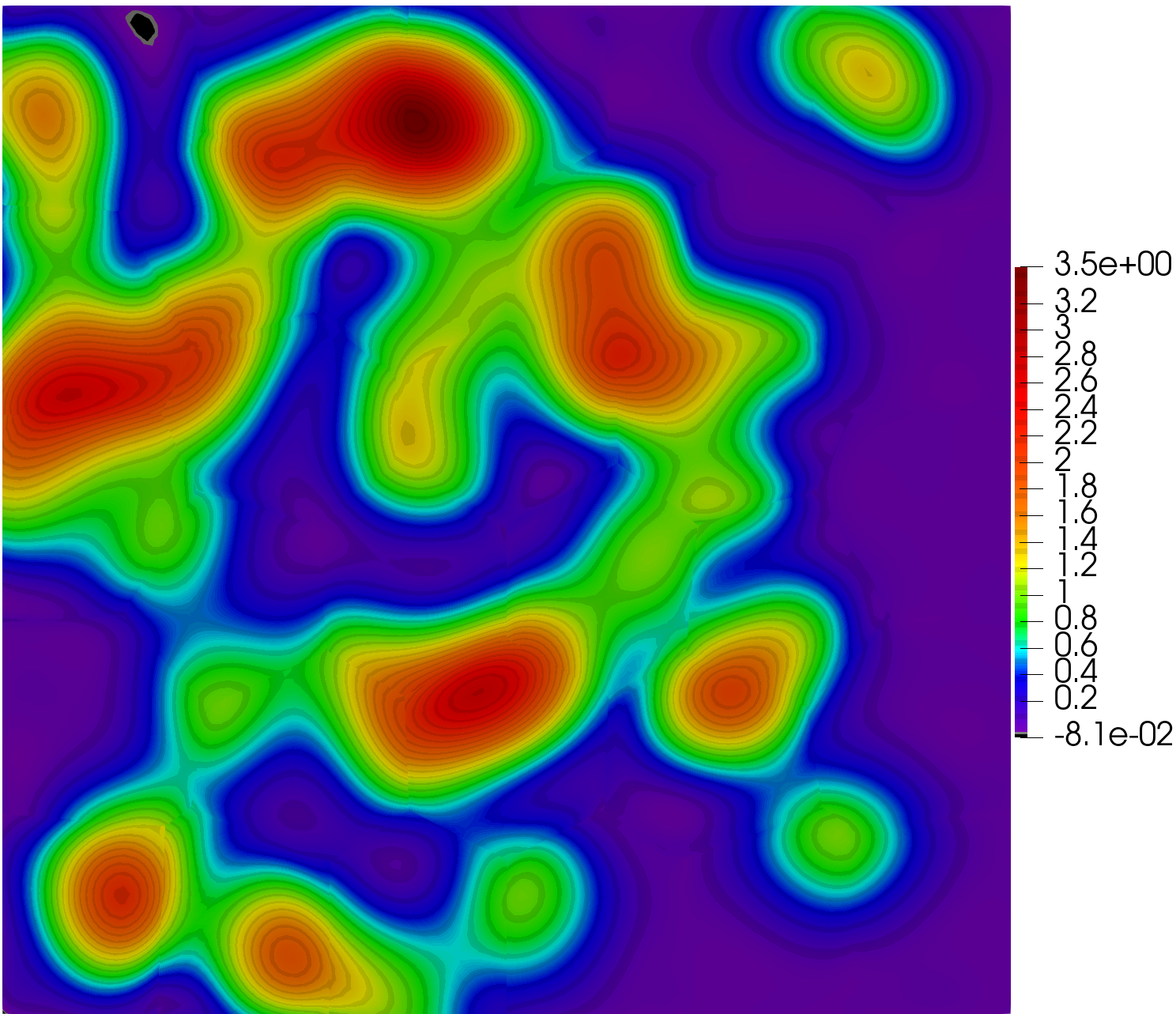}
			\subcaption{}
		\end{subfigure}
		\caption{Surface and contour plots of the
			exact solution of the exact solution~\eqref{eq:exp2d_example} for Example 2.}
		\label{fig:example1b_plots}
	\end{figure}

	In this example the Poisson problem in $\Omega = (-1,1)^2$ is considered with the exact solution of the following exponential form:
	\begin{align}\label{eq:exp2d_example}
		u(x,y) = \sum_{i=1}^{50}  \exp - \left(a_i (x - x_i)^2 + b_i (y-y_i)^2\right) \eqcomma
	\end{align}
	where $a_i$ and $b_i$ are randomly generated positive values, $(x_i,y_i)$ are randomly generated points in the domain. The plots showing the function $u$ 			in~\eqref{eq:exp2d_example} are depicted in Fig.~\ref{fig:example1b_plots}. This function is irregular with many hills randomly distributed in the domain. This example has been solved using polygonal meshes for uniform approximation orders in the elements, $p=$ 3, 5, 7, and 10.
	
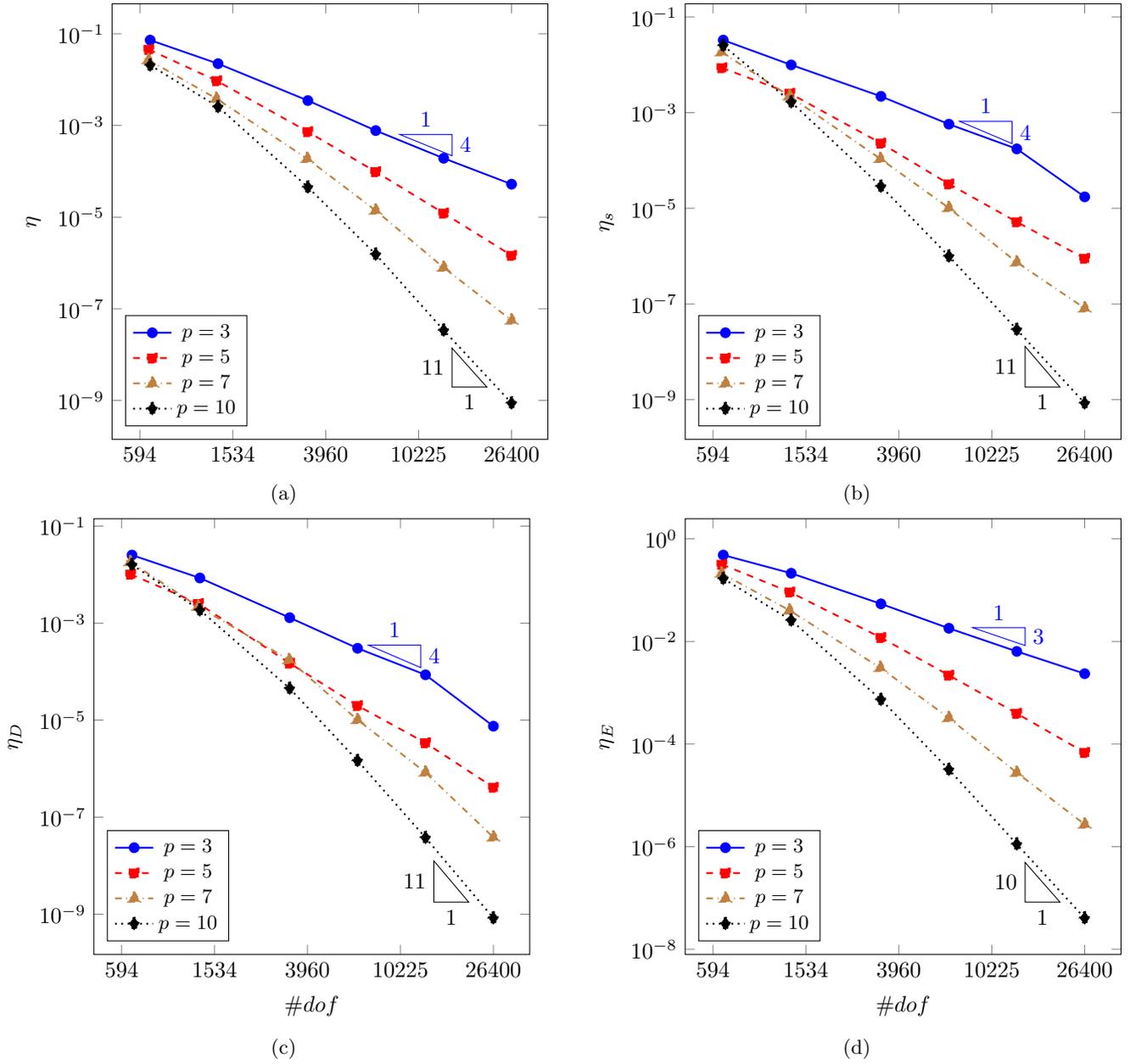
\begin{figure}\centering
		\begin{subfigure}{0.49\textwidth}\centering
			\begin{tikzpicture}[scale=1]
				\begin{axis}[width=14cm, height=14cm,
					ymode=log,
					xmode=log,
					ylabel=$\eta$,
					scale=0.55,
					legend style={at={(0.03,0.03)},anchor=south west},
					xtick = {24.372115, 39.166312, 62.928531, 101.118742, 162.480768},
					xticklabels = {594, 1534, 3960, 10225, 26400},
					ytickten={-9,-7,-5,...,0},
					]

					\addplot [solid, color=blue, thick, mark=*]
					table [x index=3, y index=4] {figures/example2/ex1b_poly_p3.txt};
					\addlegendentry{\small $p=3$}
					
					\addplot [dashed, color=red, thick, mark=square*]
					table [x index=3, y index=4] {figures/example2/ex1b_poly_p5.txt};
					\addlegendentry{\small $p=5$}
					
					\addplot [dashdotted, color=brown, thick, mark=triangle*, mark size=3pt] 
					table [x index=3, y index=4] {figures/example2/ex1b_poly_p7.txt};
					\addlegendentry{\small $p=7$}

					\addplot [dotted, color=black, thick, mark=diamond*, mark size=3pt] 
					table [x index=3, y index=4] {figures/example2/ex1b_poly_p10.txt};
					\addlegendentry{\small $p=10$}

					\logLogSlopeTTriangle{0.78}{0.08}{0.12}{11}{black};
					\logLogSlopeTriangle{0.78}{0.12}{0.70}{4}{blue};
					
				\end{axis}
			\end{tikzpicture}
			\subcaption{}
		\end{subfigure}
		\begin{subfigure}{0.49\textwidth}\centering
			\begin{tikzpicture}[scale=1]
				\begin{axis}[width=14cm, height=14cm,
					ymode=log,
					xmode=log,
					ylabel=$\eta_s$,
					scale=0.55,
					legend style={at={(0.03,0.03)},anchor=south west},
					xtick = {24.372115, 39.166312, 62.928531, 101.118742, 162.480768},
					xticklabels = {594, 1534, 3960, 10225, 26400},
					ytickten={-9,-7,...,-1},
					]
					
					\addplot [solid, color=blue, thick, mark=*]
					table [x index=3, y index=6] {figures/example2/ex1b_poly_p3.txt};
					\addlegendentry{\small $p=3$}
					
					\addplot [dashed, color=red, thick, mark=square*]
					table [x index=3, y index=6] {figures/example2/ex1b_poly_p5.txt};
					\addlegendentry{\small $p=5$}
					
					\addplot [dashdotted, color=brown, thick, mark=triangle*, mark size=3pt] 
					table [x index=3, y index=6] {figures/example2/ex1b_poly_p7.txt};
					\addlegendentry{\small $p=7$}
					
					\addplot [dotted, color=black, thick, mark=diamond*, mark size=3pt] 
					table [x index=3, y index=6] {figures/example2/ex1b_poly_p10.txt};
					\addlegendentry{\small $p=10$}

					\logLogSlopeTTriangle{0.78}{0.08}{0.12}{11}{black};
					\logLogSlopeTriangle{0.75}{0.12}{0.73}{4}{blue};
					
				\end{axis}
			\end{tikzpicture}
			\subcaption{}
		\end{subfigure}

		\begin{subfigure}{0.49\textwidth}
			\begin{tikzpicture}[scale=1]
				\begin{axis}[width=14cm, height=14cm,
					ymode=log,
					xmode=log,
					xlabel=${\#dof}$,
					ylabel=$\eta_D$,
					scale=0.55,
					legend style={at={(0.03,0.03)},anchor=south west},
					xtick = {24.372115, 39.166312, 62.928531, 101.118742, 162.480768},
					xticklabels = {594, 1534, 3960, 10225, 26400},
					ytickten={-9,-7,...,-1},
					]
					
					\addplot [solid, color=blue, thick, mark=*]
					table [x index=3, y index=7] {figures/example2/ex1b_poly_p3.txt};
					\addlegendentry{\small $p=3$}
					
					\addplot [dashed, color=red, thick, mark=square*]
					table [x index=3, y index=7] {figures/example2/ex1b_poly_p5.txt};
					\addlegendentry{\small $p=5$}
					
					\addplot [dashdotted, color=brown, thick, mark=triangle*, mark size=3pt] 
					table [x index=3, y index=7] {figures/example2/ex1b_poly_p7.txt};
					\addlegendentry{\small $p=7$}
					
					\addplot [dotted, color=black, thick, mark=diamond*, mark size=3pt] 
					table [x index=3, y index=7] {figures/example2/ex1b_poly_p10.txt};
					\addlegendentry{\small $p=10$}

					\logLogSlopeTTriangle{0.78}{0.08}{0.12}{11}{black};
					\logLogSlopeTriangle{0.75}{0.12}{0.71}{4}{blue};
					
				\end{axis}
			\end{tikzpicture}
			\subcaption{}
		\end{subfigure}
		\begin{subfigure}{0.49\textwidth}\centering
			\begin{tikzpicture}[scale=1]
				\begin{axis}[width=14cm, height=14cm,
					ymode=log,
					xmode=log,
					xlabel=${\#dof}$,
					ylabel=$\eta_E$,
					scale=0.55,
					legend style={at={(0.03,0.03)},anchor=south west},
					xtick = {24.372115, 39.166312, 62.928531, 101.118742, 162.480768},
					xticklabels = {594, 1534, 3960, 10225, 26400},
					ytickten={-8,-6,...,0},
					]
					
					\addplot [solid, color=blue, thick, mark=*]
					table [x index=3, y index=5] {figures/example2/ex1b_poly_p3.txt};
					\addlegendentry{\small $p=3$}
					
					\addplot [dashed, color=red, thick, mark=square*]
					table [x index=3, y index=5] {figures/example2/ex1b_poly_p5.txt};
					\addlegendentry{\small $p=5$}
					
					\addplot [dashdotted, color=brown, thick, mark=triangle*, mark size=3pt] 
					table [x index=3, y index=5] {figures/example2/ex1b_poly_p7.txt};
					\addlegendentry{\small $p=7$}
					
					\addplot [dotted, color=black, thick, mark=diamond*, mark size=3pt] 
					table [x index=3, y index=5] {figures/example2/ex1b_poly_p10.txt};
					\addlegendentry{\small $p=10$}

					\logLogSlopeTTriangle{0.78}{0.08}{0.12}{10}{black};
					\logLogSlopeTriangle{0.78}{0.12}{0.75}{3}{blue};
					
				\end{axis}
			\end{tikzpicture}
			\subcaption{}
		\end{subfigure}

		\caption{Convergence study of the \PFDG{} method 
			for Example 2 on polygonal meshes.
			Study is conducted in the $L^2$ norm
			on the (a) whole domain, 
			(b) mesh skeleton and (c) outer boundary, and
			(d) energy seminorm over the whole domain.
			Scale is $ \frac{1}{2}\log $--$\log$.} 
		\label{fig:example1b_convergence}
	\end{figure}

	Convergence plots for the \PFDG{} method are presented in Fig.~\ref{fig:example1b_convergence} on the whole domain and also along the mesh skeleton and outer boundary. The rates of the convergences are as expected for the particular orders. The computed residuals confirm  the correctness of the  \PFDG{} method.
	
	\subsection{Example 3: Poisson problem on L-shaped domain} \label{ssec:poissons_dg_lshape}
\begin{figure}
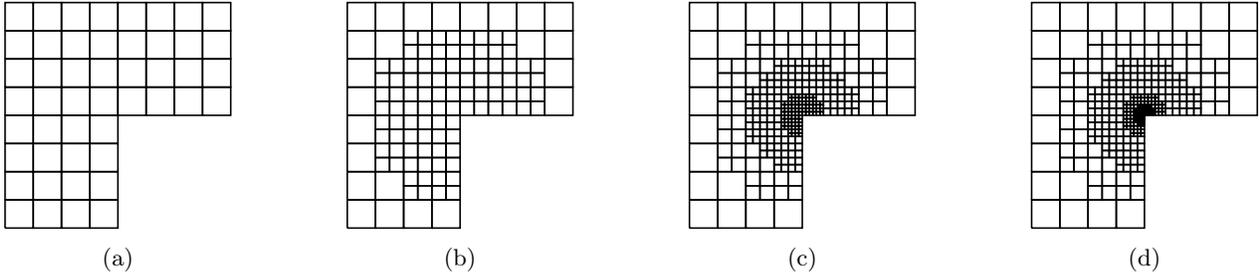
\centering
		\begin{subfigure}{0.245\textwidth}
			\centering
			\begin{tikzpicture}[scale=1.5]
				\draw[line width=0.6pt] (-1.000000, -1.000000) -- (-0.750000, -1.000000) -- (-0.750000, -0.750000) -- (-1.000000, -0.750000) -- (-1.000000, -1.000000);
\draw[line width=0.6pt] (-1.000000, -0.750000) -- (-0.750000, -0.750000) -- (-0.750000, -0.500000) -- (-1.000000, -0.500000) -- (-1.000000, -0.750000);
\draw[line width=0.6pt] (-1.000000, -0.500000) -- (-0.750000, -0.500000) -- (-0.750000, -0.250000) -- (-1.000000, -0.250000) -- (-1.000000, -0.500000);
\draw[line width=0.6pt] (-1.000000, -0.250000) -- (-0.750000, -0.250000) -- (-0.750000, 0.000000) -- (-1.000000, 0.000000) -- (-1.000000, -0.250000);
\draw[line width=0.6pt] (-1.000000, 0.000000) -- (-0.750000, 0.000000) -- (-0.750000, 0.250000) -- (-1.000000, 0.250000) -- (-1.000000, 0.000000);
\draw[line width=0.6pt] (-1.000000, 0.250000) -- (-0.750000, 0.250000) -- (-0.750000, 0.500000) -- (-1.000000, 0.500000) -- (-1.000000, 0.250000);
\draw[line width=0.6pt] (-1.000000, 0.500000) -- (-0.750000, 0.500000) -- (-0.750000, 0.750000) -- (-1.000000, 0.750000) -- (-1.000000, 0.500000);
\draw[line width=0.6pt] (-1.000000, 0.750000) -- (-0.750000, 0.750000) -- (-0.750000, 1.000000) -- (-1.000000, 1.000000) -- (-1.000000, 0.750000);
\draw[line width=0.6pt] (-0.750000, -1.000000) -- (-0.500000, -1.000000) -- (-0.500000, -0.750000) -- (-0.750000, -0.750000) -- (-0.750000, -1.000000);
\draw[line width=0.6pt] (-0.750000, -0.750000) -- (-0.500000, -0.750000) -- (-0.500000, -0.500000) -- (-0.750000, -0.500000) -- (-0.750000, -0.750000);
\draw[line width=0.6pt] (-0.750000, -0.500000) -- (-0.500000, -0.500000) -- (-0.500000, -0.250000) -- (-0.750000, -0.250000) -- (-0.750000, -0.500000);
\draw[line width=0.6pt] (-0.750000, -0.250000) -- (-0.500000, -0.250000) -- (-0.500000, 0.000000) -- (-0.750000, 0.000000) -- (-0.750000, -0.250000);
\draw[line width=0.6pt] (-0.750000, 0.000000) -- (-0.500000, 0.000000) -- (-0.500000, 0.250000) -- (-0.750000, 0.250000) -- (-0.750000, 0.000000);
\draw[line width=0.6pt] (-0.750000, 0.250000) -- (-0.500000, 0.250000) -- (-0.500000, 0.500000) -- (-0.750000, 0.500000) -- (-0.750000, 0.250000);
\draw[line width=0.6pt] (-0.750000, 0.500000) -- (-0.500000, 0.500000) -- (-0.500000, 0.750000) -- (-0.750000, 0.750000) -- (-0.750000, 0.500000);
\draw[line width=0.6pt] (-0.750000, 0.750000) -- (-0.500000, 0.750000) -- (-0.500000, 1.000000) -- (-0.750000, 1.000000) -- (-0.750000, 0.750000);
\draw[line width=0.6pt] (-0.500000, -1.000000) -- (-0.250000, -1.000000) -- (-0.250000, -0.750000) -- (-0.500000, -0.750000) -- (-0.500000, -1.000000);
\draw[line width=0.6pt] (-0.500000, -0.750000) -- (-0.250000, -0.750000) -- (-0.250000, -0.500000) -- (-0.500000, -0.500000) -- (-0.500000, -0.750000);
\draw[line width=0.6pt] (-0.500000, -0.500000) -- (-0.250000, -0.500000) -- (-0.250000, -0.250000) -- (-0.500000, -0.250000) -- (-0.500000, -0.500000);
\draw[line width=0.6pt] (-0.500000, -0.250000) -- (-0.250000, -0.250000) -- (-0.250000, 0.000000) -- (-0.500000, 0.000000) -- (-0.500000, -0.250000);
\draw[line width=0.6pt] (-0.500000, 0.000000) -- (-0.250000, 0.000000) -- (-0.250000, 0.250000) -- (-0.500000, 0.250000) -- (-0.500000, 0.000000);
\draw[line width=0.6pt] (-0.500000, 0.250000) -- (-0.250000, 0.250000) -- (-0.250000, 0.500000) -- (-0.500000, 0.500000) -- (-0.500000, 0.250000);
\draw[line width=0.6pt] (-0.500000, 0.500000) -- (-0.250000, 0.500000) -- (-0.250000, 0.750000) -- (-0.500000, 0.750000) -- (-0.500000, 0.500000);
\draw[line width=0.6pt] (-0.500000, 0.750000) -- (-0.250000, 0.750000) -- (-0.250000, 1.000000) -- (-0.500000, 1.000000) -- (-0.500000, 0.750000);
\draw[line width=0.6pt] (-0.250000, -1.000000) -- (0.000000, -1.000000) -- (0.000000, -0.750000) -- (-0.250000, -0.750000) -- (-0.250000, -1.000000);
\draw[line width=0.6pt] (-0.250000, -0.750000) -- (0.000000, -0.750000) -- (0.000000, -0.500000) -- (-0.250000, -0.500000) -- (-0.250000, -0.750000);
\draw[line width=0.6pt] (-0.250000, -0.500000) -- (0.000000, -0.500000) -- (0.000000, -0.250000) -- (-0.250000, -0.250000) -- (-0.250000, -0.500000);
\draw[line width=0.6pt] (-0.250000, -0.250000) -- (0.000000, -0.250000) -- (0.000000, 0.000000) -- (-0.250000, 0.000000) -- (-0.250000, -0.250000);
\draw[line width=0.6pt] (-0.250000, 0.000000) -- (0.000000, 0.000000) -- (0.000000, 0.250000) -- (-0.250000, 0.250000) -- (-0.250000, 0.000000);
\draw[line width=0.6pt] (-0.250000, 0.250000) -- (0.000000, 0.250000) -- (0.000000, 0.500000) -- (-0.250000, 0.500000) -- (-0.250000, 0.250000);
\draw[line width=0.6pt] (-0.250000, 0.500000) -- (0.000000, 0.500000) -- (0.000000, 0.750000) -- (-0.250000, 0.750000) -- (-0.250000, 0.500000);
\draw[line width=0.6pt] (-0.250000, 0.750000) -- (0.000000, 0.750000) -- (0.000000, 1.000000) -- (-0.250000, 1.000000) -- (-0.250000, 0.750000);
\draw[line width=0.6pt] (0.000000, 0.000000) -- (0.250000, 0.000000) -- (0.250000, 0.250000) -- (0.000000, 0.250000) -- (0.000000, 0.000000);
\draw[line width=0.6pt] (0.000000, 0.250000) -- (0.250000, 0.250000) -- (0.250000, 0.500000) -- (0.000000, 0.500000) -- (0.000000, 0.250000);
\draw[line width=0.6pt] (0.000000, 0.500000) -- (0.250000, 0.500000) -- (0.250000, 0.750000) -- (0.000000, 0.750000) -- (0.000000, 0.500000);
\draw[line width=0.6pt] (0.000000, 0.750000) -- (0.250000, 0.750000) -- (0.250000, 1.000000) -- (0.000000, 1.000000) -- (0.000000, 0.750000);
\draw[line width=0.6pt] (0.250000, 0.000000) -- (0.500000, 0.000000) -- (0.500000, 0.250000) -- (0.250000, 0.250000) -- (0.250000, 0.000000);
\draw[line width=0.6pt] (0.250000, 0.250000) -- (0.500000, 0.250000) -- (0.500000, 0.500000) -- (0.250000, 0.500000) -- (0.250000, 0.250000);
\draw[line width=0.6pt] (0.250000, 0.500000) -- (0.500000, 0.500000) -- (0.500000, 0.750000) -- (0.250000, 0.750000) -- (0.250000, 0.500000);
\draw[line width=0.6pt] (0.250000, 0.750000) -- (0.500000, 0.750000) -- (0.500000, 1.000000) -- (0.250000, 1.000000) -- (0.250000, 0.750000);
\draw[line width=0.6pt] (0.500000, 0.000000) -- (0.750000, 0.000000) -- (0.750000, 0.250000) -- (0.500000, 0.250000) -- (0.500000, 0.000000);
\draw[line width=0.6pt] (0.500000, 0.250000) -- (0.750000, 0.250000) -- (0.750000, 0.500000) -- (0.500000, 0.500000) -- (0.500000, 0.250000);
\draw[line width=0.6pt] (0.500000, 0.500000) -- (0.750000, 0.500000) -- (0.750000, 0.750000) -- (0.500000, 0.750000) -- (0.500000, 0.500000);
\draw[line width=0.6pt] (0.500000, 0.750000) -- (0.750000, 0.750000) -- (0.750000, 1.000000) -- (0.500000, 1.000000) -- (0.500000, 0.750000);
\draw[line width=0.6pt] (0.750000, 0.000000) -- (1.000000, 0.000000) -- (1.000000, 0.250000) -- (0.750000, 0.250000) -- (0.750000, 0.000000);
\draw[line width=0.6pt] (0.750000, 0.250000) -- (1.000000, 0.250000) -- (1.000000, 0.500000) -- (0.750000, 0.500000) -- (0.750000, 0.250000);
\draw[line width=0.6pt] (0.750000, 0.500000) -- (1.000000, 0.500000) -- (1.000000, 0.750000) -- (0.750000, 0.750000) -- (0.750000, 0.500000);
\draw[line width=0.6pt] (0.750000, 0.750000) -- (1.000000, 0.750000) -- (1.000000, 1.000000) -- (0.750000, 1.000000) -- (0.750000, 0.750000);
			\end{tikzpicture}
			\subcaption{}
		\end{subfigure}
		\begin{subfigure}{0.245\textwidth}
			\centering
			\begin{tikzpicture}[scale=1.5]
				\input{figures/example_lshape_pois/lshape_mesh_120}
			\end{tikzpicture}
			\subcaption{}
		\end{subfigure}
		\begin{subfigure}{0.245\textwidth}
			\centering
			\begin{tikzpicture}[scale=1.5]
				\input{figures/example_lshape_pois/lshape_mesh_309}
			\end{tikzpicture}
			\subcaption{}
		\end{subfigure}
		\begin{subfigure}{0.245\textwidth}
			\centering
			\begin{tikzpicture}[scale=1.5]
				\input{figures/example_lshape_pois/lshape_mesh_597}
			\end{tikzpicture}
			\subcaption{}
		\end{subfigure}
		\caption{Quadrilateral meshes for the L-shaped domain in Example 3. Mesh is refined in the vicinity of the reentrant corner.}
		\label{fig:lshape_qmeshes}
	\end{figure}

	We consider the standard  benchmark Poisson problem on the L-shaped domain, with the exact solution
	\begin{align}
		\temperature(r,\theta) = r^{\frac{2}{3}} \sin\left( \frac{2}{3}\theta \right) .
	\end{align}
	This function is harmonic, so the right-hand side of~\eqref{eq:heatcond1} is zero. The exact solution has  a derivative singularity at the origin, and hence the error in the numerical solution concentrates at the origin. In this example, the vertices of the L-shaped domain are chosen as $(-1, -1)$, $(0, -1)$, $(0, 0)$, $(1, 0)$, $(1, 1)$, $(-1, 1)$ and the domain is meshed with structured quadrilaterals and then refined several times in the vicinity of the origin (see Fig.~\ref{fig:lshape_qmeshes}). The approximation orders in the elements are homogeneous in the entire mesh and calculations are performed on elements with orders $p=$ 3, 5, 7, and 10.

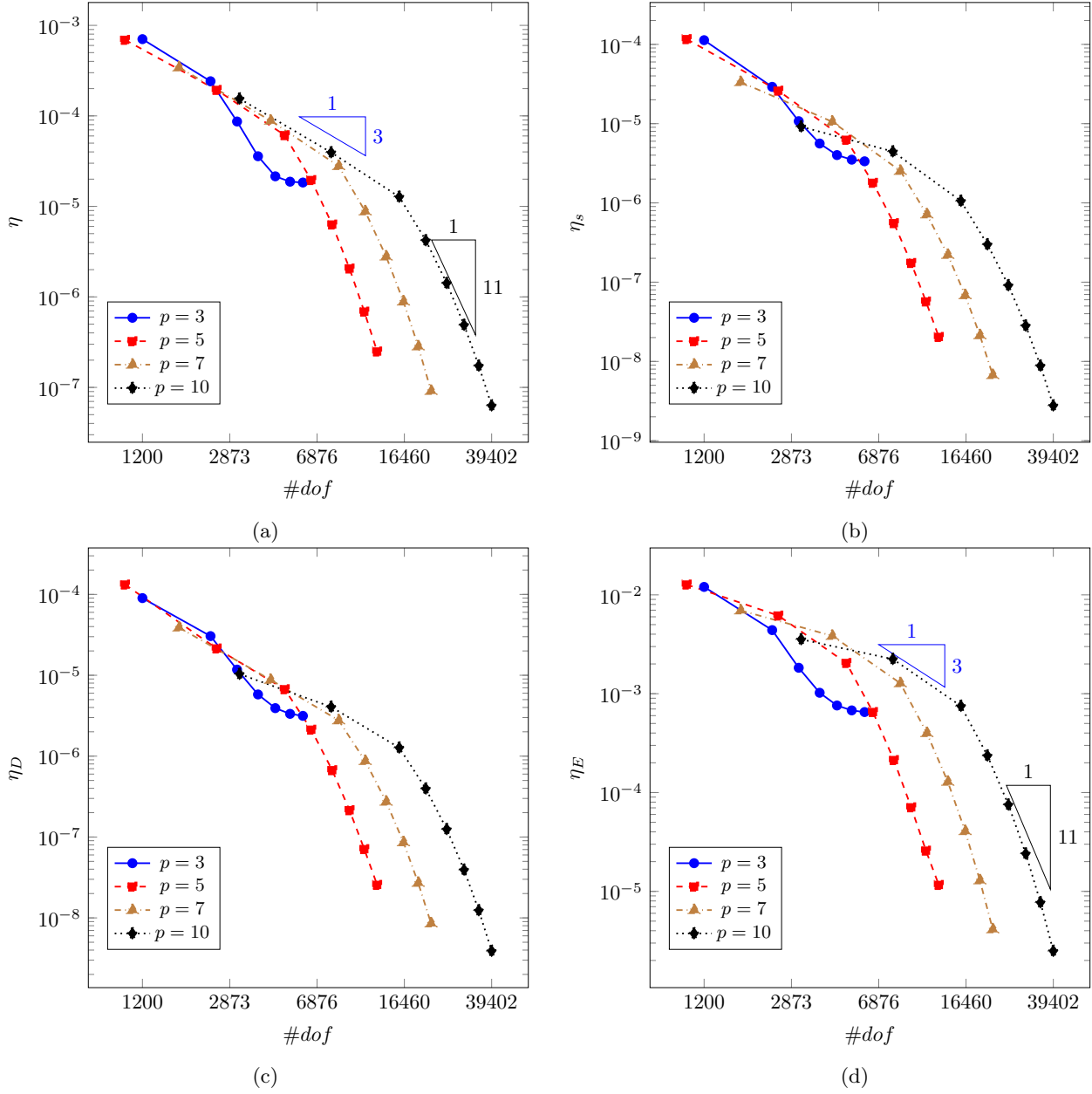
\begin{figure}\centering
		\begin{subfigure}{0.49\textwidth}
			\centering
			\begin{tikzpicture}[scale=0.9]
				\begin{axis}[width=14cm, height=14cm,
					xlabel=${\#dof}$,
					ylabel=$\eta$,
					scale=0.6,
					ymode=log,
					xmode=log,
					legend style={at={(0.3,0.32)}},
					xtick = {34.641016, 53.600373, 82.921650, 128.296532, 198.499370},
					xticklabels = {1200, 2873, 6876, 16460, 39402},				
					]

					\addplot [solid, color=blue, thick, mark=*]
					table [x index=3, y index=4] {figures/example_lshape_pois/lshape_results_p3.txt};
					\addlegendentry{\small $p=3$}
					
					\addplot [dashed, color=red, thick, mark=square*] 
					table [x index=3, y index=4] {figures/example_lshape_pois/lshape_results_p5.txt};
					\addlegendentry{\small $p=5$}
					
					\addplot [dashdotted, color=brown, thick, mark=triangle*, mark size=3pt]  
					table [x index=3, y index=4] {figures/example_lshape_pois/lshape_results_p7.txt};
					\addlegendentry{\small $p=7$}

					\addplot[dotted, color=black, thick, mark=diamond*, mark size=3pt]
					table [x index=3, y index=4] {figures/example_lshape_pois/lshape_results_p10.txt};
					\addlegendentry{\small $p=10$}

					\logLogSlopeTriangle{0.63}{0.15}{0.74}{3}{blue};
					\logLogSlopeTriangle{0.88}{0.1}{0.46}{11}{black};			
				\end{axis}
			\end{tikzpicture}
			\subcaption{}
		\end{subfigure}
		\begin{subfigure}{0.49\textwidth}
			\begin{tikzpicture}[scale=0.9]
				\begin{axis}[width=14cm, height=14cm,
					xlabel=${\#dof}$,
					ylabel=$\eta_s$,
					scale=0.6,
					ymode=log,
					xmode=log,
					legend style={at={(0.3,0.32)}},
					xtick = {34.641016, 53.600373, 82.921650, 128.296532, 198.499370},
					xticklabels = {1200, 2873, 6876, 16460, 39402},	
					]
					\addplot [solid, color=blue, thick, mark=*]
					table [x index=3, y index=6] {figures/example_lshape_pois/lshape_results_p3.txt};
					\addlegendentry{\small $p=3$}
					
					\addplot [dashed, color=red, thick, mark=square*] 
					table [x index=3, y index=6] {figures/example_lshape_pois/lshape_results_p5.txt};
					\addlegendentry{\small $p=5$}
					
					\addplot [dashdotted, color=brown, thick, mark=triangle*, mark size=3pt]  
					table [x index=3, y index=6] {figures/example_lshape_pois/lshape_results_p7.txt};
					\addlegendentry{\small $p=7$}

					\addplot[dotted, color=black, thick, mark=diamond*, mark size=3pt]
					table [x index=3, y index=6] {figures/example_lshape_pois/lshape_results_p10.txt};
					\addlegendentry{\small $p=10$}

				\end{axis}
			\end{tikzpicture}
			\subcaption{}
		\end{subfigure}

		\begin{subfigure}{0.49\textwidth}
			\centering
			\begin{tikzpicture}[scale=0.9]
				\begin{axis}[width=14cm, height=14cm,
					xlabel=${\#dof}$,
					ylabel=$\eta_D$,
					scale=0.6,
					ymode=log,
					xmode=log,
					legend style={at={(0.3,0.32)}},
					xtick = {34.641016, 53.600373, 82.921650, 128.296532, 198.499370},
					xticklabels = {1200, 2873, 6876, 16460, 39402},	
					]
					\addplot [solid, color=blue, thick, mark=*]
					table [x index=3, y index=7] {figures/example_lshape_pois/lshape_results_p3.txt};
					\addlegendentry{\small $p=3$}
					
					\addplot [dashed, color=red, thick, mark=square*] 
					table [x index=3, y index=7] {figures/example_lshape_pois/lshape_results_p5.txt};
					\addlegendentry{\small $p=5$}
					
					\addplot [dashdotted, color=brown, thick, mark=triangle*, mark size=3pt]  
					table [x index=3, y index=7] {figures/example_lshape_pois/lshape_results_p7.txt};
					\addlegendentry{\small $p=7$}

					\addplot[dotted, color=black, thick, mark=diamond*, mark size=3pt]
					table [x index=3, y index=7] {figures/example_lshape_pois/lshape_results_p10.txt};
					\addlegendentry{\small $p=10$}
					
				\end{axis}
			\end{tikzpicture}
			\subcaption{}
		\end{subfigure}
		\begin{subfigure}{0.49\textwidth}
			\begin{tikzpicture}[scale=0.9]
				\begin{axis}[width=14cm, height=14cm,
					xlabel=${\#dof}$,
					ylabel=$\eta_E$,
					scale=0.6,
					ymode=log,
					xmode=log,
					legend style={at={(0.3,0.32)}},
					xtick = {34.641016, 53.600373, 82.921650, 128.296532, 198.499370},
					xticklabels = {1200, 2873, 6876, 16460, 39402},	
					]
					\addplot [solid, color=blue, thick, mark=*]
					table [x index=3, y index=5] {figures/example_lshape_pois/lshape_results_p3.txt};
					\addlegendentry{\small $p=3$}
					
					\addplot [dashed, color=red, thick, mark=square*] 
					table [x index=3, y index=5] {figures/example_lshape_pois/lshape_results_p5.txt};
					\addlegendentry{\small $p=5$}
					
					\addplot [dashdotted, color=brown, thick, mark=triangle*, mark size=3pt]  
					table [x index=3, y index=5] {figures/example_lshape_pois/lshape_results_p7.txt};
					\addlegendentry{\small $p=7$}

					\addplot[dotted, color=black, thick, mark=diamond*, mark size=3pt]
					table [x index=3, y index=5] {figures/example_lshape_pois/lshape_results_p10.txt};
					\addlegendentry{\small $p=10$}

					\logLogSlopeTriangle{0.67}{0.15}{0.78}{3}{blue};
					\logLogSlopeTriangle{0.91}{0.1}{0.46}{11}{black};			
				\end{axis}
			\end{tikzpicture}
			\subcaption{}
		\end{subfigure}

		\caption{Convergence study of the \PFDG{} method 
			for Example 3 on quadrilateral meshes.
			Study is conducted in the $L^2$ norm
			on the (a) whole domain, 
			(cb) mesh skeleton and (c) outer boundary, and
			(d) energy seminorm over the whole domain.
			Scale is $ \frac{1}{2}\log $--$\log$.} 
		\label{fig:lshape_poiss_convergences}
	\end{figure}

	The results of this analysis are presented in the form of the convergence curves over the domain and along the mesh skeleton and Dirichlet boundary in Fig.~\ref{fig:lshape_poiss_convergences}. Due to the successive mesh $h$-refinement, steep convergence rates are obtained for all approximation orders. Also in this case the convergence rates over the domain, mesh skeleton, and Dirichlet boundary are about the same.
	
	\subsection{Example 4: Linear elasticity problem on L-shaped domain}\label{ssec:elasticity_dg}
	In this example, a benchmark elasticity problem on an L-shaped domain is considered. Plane strain conditions are assumed. The governing equations are stated in~\myeqref{eq:momentum}. The exact solution of this problem in polar coordinates is:
	\begin{align}
		\label{eq:ex3_exact}
		\begin{aligned}
			\displusk_r(r,\theta) &= \frac{1}{2\lamemu} r^{\beta} \Bigl[
			-(\beta+1)  \cos\left( (\beta+1) \theta \right)
			+ (C_2- \beta-1 ) C_1  \cos\left( (\beta-1) \theta \right)
			\Bigr] \eqcomma
			\\
			\displusk_\theta(r,\theta) &= \frac{1}{2\lamemu} r^{\beta} \Bigl[
			(\beta+1)  \sin\left( (\beta+1) \theta \right)
			+ (C_2 +\beta-1) C_1  \sin\left( (\beta-1) \theta \right) 
			\Bigr] \eqdot
		\end{aligned}
	\end{align}
	where
	$C_1 = -\cos\left((\beta+1) \omega \right)/\cos\left((\beta-1) \omega \right)$,
	$C_2 = 2 \left( \lamelambda + 2\lamemu \right) / \left(\lamelambda + \lamemu \right)$,
	$\omega=3\pi/4$ and the critical exponent $\beta$ is the positive solution of the equation
	$\beta\sin(2\omega) + \sin(2\omega\beta) = 0$, and
	so $\beta\approx 0.544483737$.
	The calculations have been performed for Young's modulus $\youngmodul = 1$ and shear modulus
	$\lamemu=0.2$.

\begin{figure}\centering
		\begin{subfigure}[b]{0.45\textwidth}
			\centering
			\includegraphics[width=0.9\textwidth]{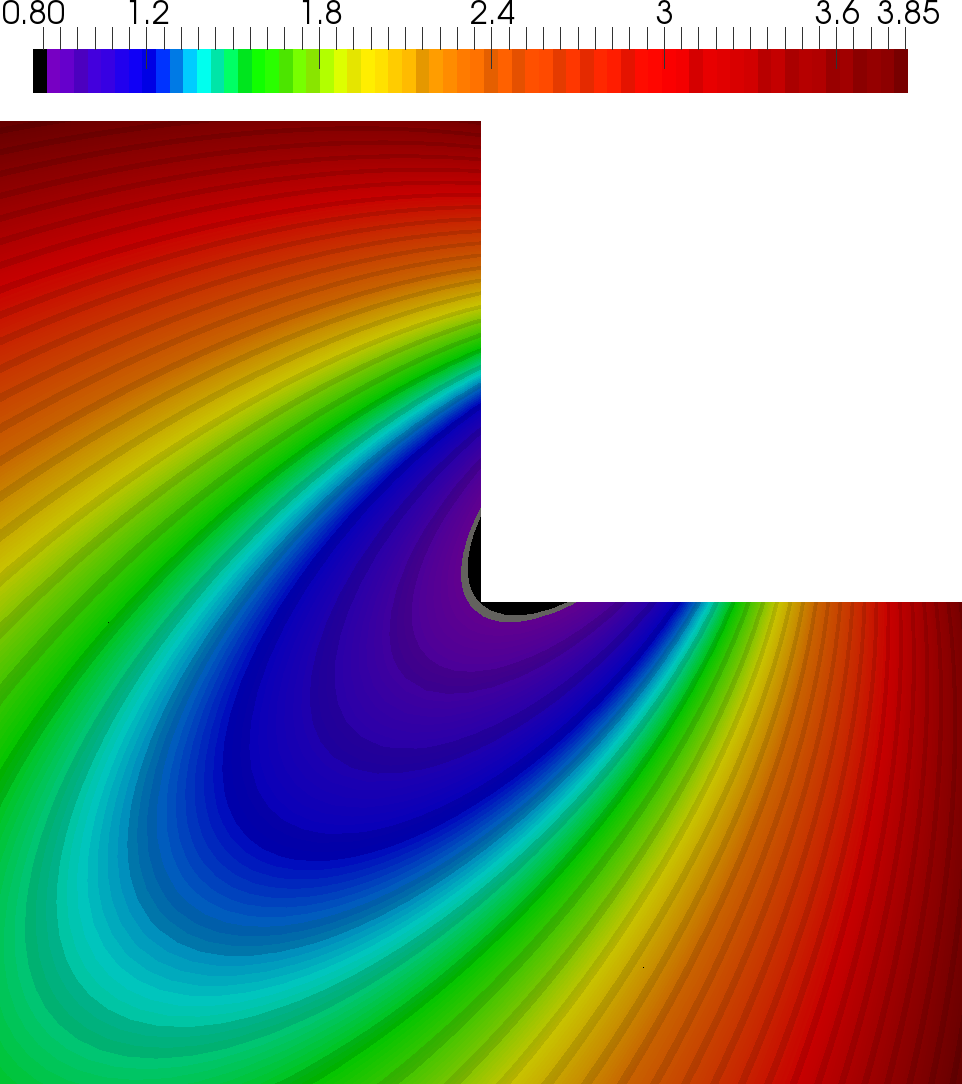}
			\subcaption{}
			\label{fig:ex_bench_exact_u}
		\end{subfigure}
		\hfill
		\begin{subfigure}[b]{0.45\textwidth}
			\centering
			\includegraphics[width=0.9\textwidth]{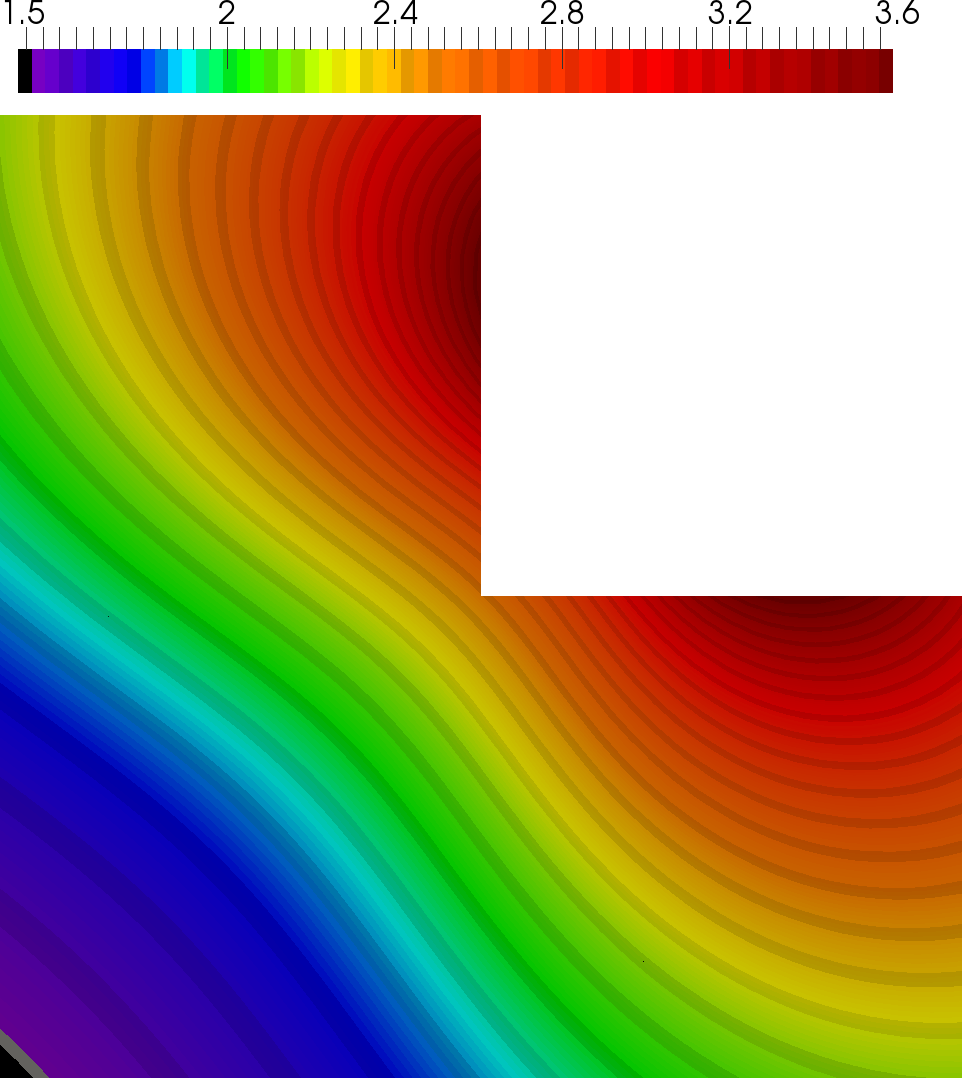}
			\subcaption{}
			\label{fig:ex_bench_exact_sv}
		\end{subfigure}
		\caption{Contour plot of (a) $\norm{\displu}$ and (b) von Mises stress on the L-shaped domain in Example 4.}
		\label{fig:ex_bench_exact}
	\end{figure}

	Originally the benchmark example has been solved on the so-called rotated L-shaped domain and with singularity at the origin, see Carstensen and Gedicke.\cite{CARSTENSEN2016245} In this paper, the  L-shaped domain is considered with vertices at points: $(-1, -1),\ (1, -1),\ (1, 0), \ (0, 0), \ (0, 1), \ (-1, 1)$ and the singularity is shifted out of the domain. To achieve the solution in such a domain, both the global coordinates as well as the displacements need to be appropriately transformed, as shown below. The polar coordinates are constructed using the Cartesian auxiliary coordinates $\xi$ and $\eta$:
	\begin{align}
		r = \sqrt{\xi^2 + \eta^2 } 
		\,,\quad
		\theta = \tan^{-1} \left( \frac{\eta}{\xi} \right)
	\end{align}
	where the auxiliary coordinates $(\xi,\eta)$ come from the following transformation of the global coordinates using the transformation angle $\phi$:
	\begin{align}
		\begin{aligned}
			&
			\phi=-3\pi/4 
			\,,\quad
			\xi = \cos(\phi) ({\cargone} - \cargone_s) + \sin(\phi) ({\cargtwo - \cargtwo_s}) 
			\,,\quad
			\eta = -\sin(\phi) ({\cargone} - \cargone_s) + \cos(\phi) ({\cargtwo - \cargtwo_s}) \eqcomma
		\end{aligned}
	\end{align}
	where $\cargone_s=\frac{1}{2}$ and $\cargtwo_s=\frac{1}{2}$ are additional shifting parameters that are set to avoid a singular solution at the origin of the global coordinates. When the displacements in  polar coordinates are calculated the displacements in  global coordinates are obtained by the same transformation angle $\phi$:
	\begin{align}
		\begin{aligned}
			&
			\displusk_x(x,y) = \cos(\phi) \displusk_r - \sin(\phi) \displusk_\theta
			\,,\quad
			\displusk_y(x,y) = \sin(\phi) \displusk_r + \cos(\phi) \displusk_\theta \eqdot
		\end{aligned}
	\end{align}
	The maps of exact solution in the displacements norm as well as von Mises stresses are shown in Fig.~\ref{fig:ex_bench_exact}.

	\begin{figure}
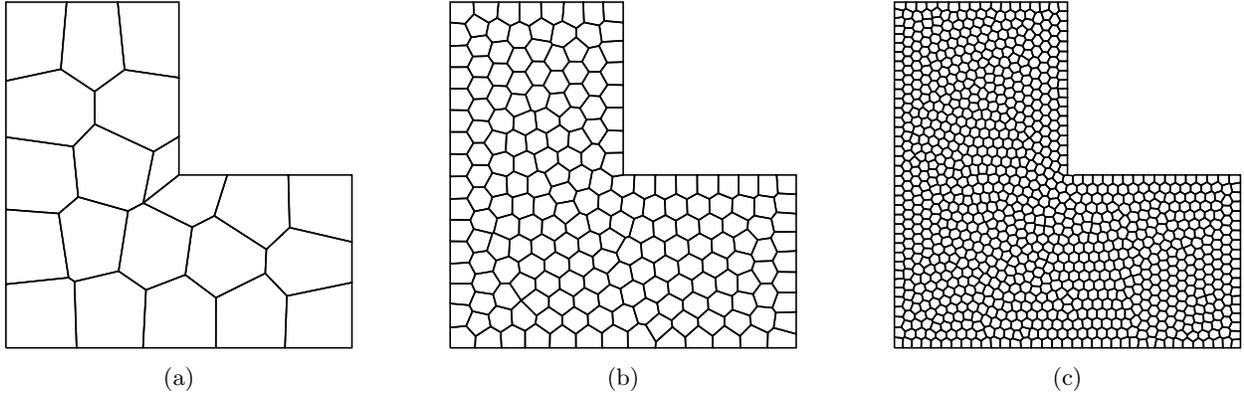
\centering	
		\begin{subfigure}[b]{0.32\textwidth}
			\centering
			\begin{tikzpicture}[scale=2.3]
				\draw[line width=0.6pt] (-0.348847, -0.556412) -- (-0.190340, -0.664818) -- (0.037132, -0.577685) -- (0.076479, -0.302053) -- (-0.206416, -0.162196) -- (-0.293892, -0.208599) -- (-0.348847, -0.556412);
\draw[line width=0.6pt] (-0.682733, 0.611971) -- (-0.487341, 0.482546) -- (-0.313939, 0.608732) -- (-0.353988, 1.000000) -- (-0.649089, 1.000000) -- (-0.682733, 0.611971);
\draw[line width=0.6pt] (-0.591467, -1.000000) -- (-0.207749, -1.000000) -- (-0.190340, -0.664818) -- (-0.348847, -0.556412) -- (-0.603810, -0.616460) -- (-0.591467, -1.000000);
\draw[line width=0.6pt] (-0.486936, 0.293968) -- (-0.145997, 0.134946) -- (-0.000000, 0.222160) -- (-0.000000, 0.562853) -- (-0.313939, 0.608732) -- (-0.487341, 0.482546) -- (-0.486936, 0.293968);
\draw[line width=0.6pt] (-1.000000, 0.544006) -- (-0.682733, 0.611971) -- (-0.649089, 1.000000) -- (-1.000000, 1.000000) -- (-1.000000, 0.544006);
\draw[line width=0.6pt] (-0.313939, 0.608732) -- (-0.000000, 0.562853) -- (-0.000000, 1.000000) -- (-0.353988, 1.000000) -- (-0.313939, 0.608732);
\draw[line width=0.6pt] (0.037132, -0.577685) -- (0.213740, -0.709968) -- (0.498828, -0.568862) -- (0.505447, -0.430067) -- (0.204828, -0.236335) -- (0.076479, -0.302053) -- (0.037132, -0.577685);
\draw[line width=0.6pt] (0.204828, -0.236335) -- (0.505447, -0.430067) -- (0.639678, -0.304075) -- (0.632188, -0.000000) -- (0.279688, -0.000000) -- (0.204828, -0.236335);
\draw[line width=0.6pt] (-0.582398, -0.127215) -- (-0.293892, -0.208599) -- (-0.206416, -0.162196) -- (-0.145997, 0.134946) -- (-0.486936, 0.293968) -- (-0.611827, 0.164670) -- (-0.582398, -0.127215);
\draw[line width=0.6pt] (-1.000000, 0.219188) -- (-0.611827, 0.164670) -- (-0.486936, 0.293968) -- (-0.487341, 0.482546) -- (-0.682733, 0.611971) -- (-1.000000, 0.544006) -- (-1.000000, 0.219188);
\draw[line width=0.6pt] (-1.000000, -1.000000) -- (-0.591467, -1.000000) -- (-0.603810, -0.616460) -- (-0.637739, -0.593752) -- (-1.000000, -0.631978) -- (-1.000000, -1.000000);
\draw[line width=0.6pt] (-1.000000, -0.631978) -- (-0.637739, -0.593752) -- (-0.694910, -0.222803) -- (-1.000000, -0.199310) -- (-1.000000, -0.631978);
\draw[line width=0.6pt] (0.612877, -1.000000) -- (1.000000, -1.000000) -- (1.000000, -0.622084) -- (0.624137, -0.699773) -- (0.612877, -1.000000);
\draw[line width=0.6pt] (0.498828, -0.568862) -- (0.624137, -0.699773) -- (1.000000, -0.622084) -- (1.000000, -0.387509) -- (0.639678, -0.304075) -- (0.505447, -0.430067) -- (0.498828, -0.568862);
\draw[line width=0.6pt] (-1.000000, -0.199310) -- (-0.694910, -0.222803) -- (-0.582398, -0.127215) -- (-0.611827, 0.164670) -- (-1.000000, 0.219188) -- (-1.000000, -0.199310);
\draw[line width=0.6pt] (-0.207749, -1.000000) -- (0.213368, -1.000000) -- (0.213740, -0.709968) -- (0.037132, -0.577685) -- (-0.190340, -0.664818) -- (-0.207749, -1.000000);
\draw[line width=0.6pt] (-0.637739, -0.593752) -- (-0.603810, -0.616460) -- (-0.348847, -0.556412) -- (-0.293892, -0.208599) -- (-0.582398, -0.127215) -- (-0.694910, -0.222803) -- (-0.637739, -0.593752);
\draw[line width=0.6pt] (0.000000, 0.000000) -- (-0.000000, 0.222160) -- (-0.145997, 0.134946) -- (-0.206416, -0.162196) -- (0.000000, 0.000000);
\draw[line width=0.6pt] (0.213368, -1.000000) -- (0.612877, -1.000000) -- (0.624137, -0.699773) -- (0.498828, -0.568862) -- (0.213740, -0.709968) -- (0.213368, -1.000000);
\draw[line width=0.6pt] (0.639678, -0.304075) -- (1.000000, -0.387509) -- (1.000000, -0.000000) -- (0.632188, -0.000000) -- (0.639678, -0.304075);
\draw[line width=0.6pt] (-0.206416, -0.162196) -- (0.076479, -0.302053) -- (0.204828, -0.236335) -- (0.279688, -0.000000) -- (0.000000, 0.000000) -- (-0.206416, -0.162196);		  
			\end{tikzpicture}
			\subcaption{}
		\end{subfigure}
		\begin{subfigure}[b]{0.32\textwidth}
			\centering
			\begin{tikzpicture}[scale=2.3]
				\input{figures/example3/saved_mesh_p200}
			\end{tikzpicture}
			\subcaption{}
		\end{subfigure}
		\begin{subfigure}[b]{0.32\textwidth}
			\centering
			\begin{tikzpicture}[scale=2.3]
				\input{figures/example3/saved_mesh_p1000}
			\end{tikzpicture}
			\subcaption{}
		\end{subfigure}
		\caption{Sample polygonal meshes used
			for Example 4 with (a) 21, (b) 200
			and (c) 1000 elements.}
		\label{fig:ex_bench_meshes}
	\end{figure}

\begin{figure}[!htb]
		\centering
		\begin{subfigure}{0.49\textwidth}
			\begin{tikzpicture}[scale=0.9]
				\begin{axis}[width=14cm, height=14cm,
					xlabel=${\#dof}$,
					ylabel=$\eta_D$,
					scale=0.6,
					ymode=log,
					xmode=log,
					xmin=10^(1.2), 	xmax=10^(2.3),
					legend style={at={(0.3,0.32)}},
					xtick = {20.493902, 35.749126, 62.353829, 108.770400, 189.736660},
					xticklabels = {420, 1278, 3888, 11831, 36000},
					]
					\addplot [solid, color=blue, thick, mark=*]
					table [x index=3, y index=6] {figures/example3/conv2_pss_polyg_p3.txt};
					\addlegendentry{\small $p=3$}
					
					\addplot [dashed, color=red, thick, mark=square*]
					table [x index=3, y index=6] {figures/example3/conv2_pss_polyg_p5.txt};
					\addlegendentry{\small $p=5$}
					
					\addplot [dashdotted, color=brown, thick, mark=triangle*, mark size=3pt]
					table [x index=3, y index=6] {figures/example3/conv2_pss_polyg_p7.txt};
					\addlegendentry{\small $p=7$}
					
					\addplot [dotted, color=black, thick, mark=diamond*, mark size=3pt]
					table [x index=3, y index=6] {figures/example3/conv2_pss_polyg_p10.txt};
					\addlegendentry{\small $p=10$}
					
					\logLogSlopeTriangle{0.70}{0.15}{0.73}{4}{blue};
					\logLogSlopeTTriangle{0.48}{0.1}{0.31}{10}{black};	
					\logLogSlopeTTriangle{0.75}{0.1}{0.25}{8}{brown};			
				\end{axis}
			\end{tikzpicture}
			\subcaption{}
		\end{subfigure}
		\begin{subfigure}{0.49\textwidth}
			\centering
			\begin{tikzpicture}[scale=0.9]
				\begin{axis}[width=14cm, height=14cm,
					xlabel=${\#dof}$,
					ylabel=$\eta_D$,
					scale=0.6,
					ymode=log,
					xmode=log,
					xmin=10^(1.2), 	xmax=10^(2.3),
					legend style={at={(0.3,0.32)}},
					xtick = {20.493902, 35.749126, 62.353829, 108.770400, 189.736660},
					xticklabels = {420, 1278, 3888, 11831, 36000},
					]
					\addplot [solid, color=blue, thick, mark=*] 
					table [x index=3, y index=6] {figures/example3/conv2_pss_square_p3.txt};
					\addlegendentry{\small $p=3$}
					
					\addplot [dashed, color=red, thick, mark=square*] 
					table [x index=3, y index=6] {figures/example3/conv2_pss_square_p5.txt};
					\addlegendentry{\small $p=5$}

					\addplot [dashdotted, color=brown, thick, mark=triangle*, mark size=3pt] 
					table [x index=3, y index=6]{figures/example3/conv2_pss_square_p7.txt};
					\addlegendentry{\small $p=7$}
					
					\addplot[dotted, color=black, thick, mark=diamond*, mark size=3pt]
					table [x index=3, y index=6]{figures/example3/conv2_pss_square_p10.txt};
					\addlegendentry{\small $p=10$}
					
					\logLogSlopeTriangle{0.76}{0.18}{0.82}{4}{blue};
					\logLogSlopeTTriangle{0.60}{0.12}{0.18}{10}{blue};
					\logLogSlopeTTriangle{0.72}{0.1}{0.38}{8}{brown};	
					
				\end{axis}
			\end{tikzpicture}
			\subcaption{}
		\end{subfigure}
		
		\caption{Convergence study of the \PFDG{} method 
			along the Dirichlet boundary
			for Example 4 on (a) polygonal meshes
			and (b) quadrilateral meshes. Scale is 
			$\frac{1}{2}\log$--$\log$.}
		\label{fig:dg_convergence_elast_boundary}
	\end{figure}
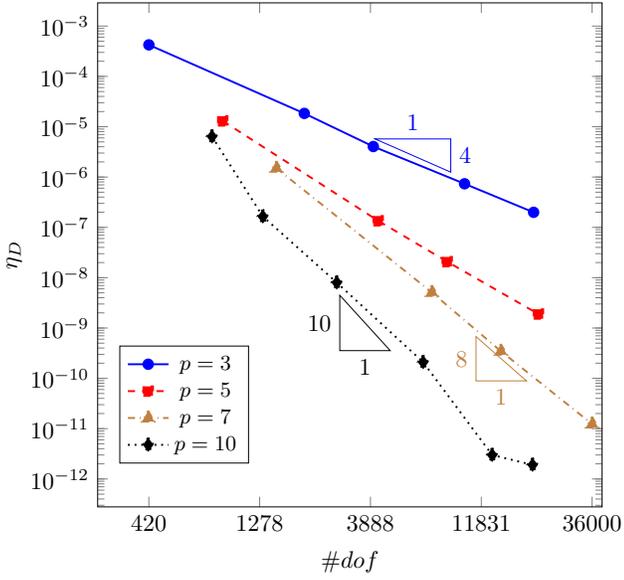
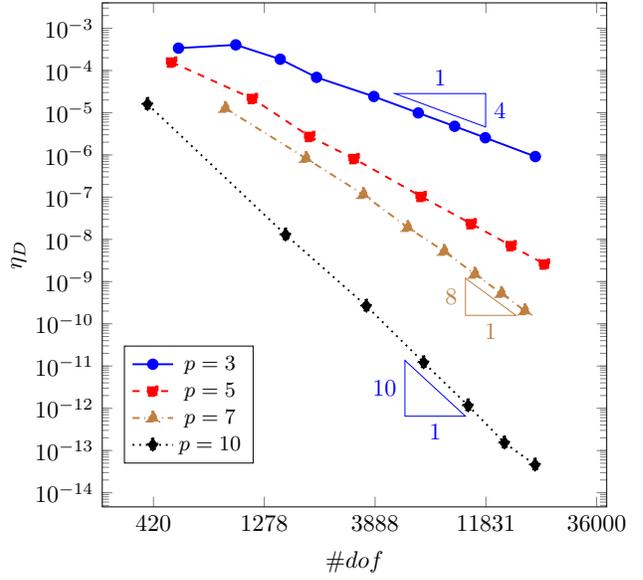

	This example has been solved on	meshes with polygonal and quadrilateral elements for various mesh densities and approximation orders $p=3$ up to $p=10$ on each mesh. Three polygonal meshes are shown in Fig.~\ref{fig:ex_bench_meshes}.	The results are presented in the form of convergence rates for various values of $p$ versus the square root of the number of degrees of freedom (see Fig.~\ref{fig:dg_convergence_elast_domain}). Numerically computed convergence rates are consistent with the theoretical estimates. Note that the rate of convergence in the $L^2$ norm is $\porder+1$ except for the case when $\porder=10$ on polygonal meshes. The convergence curve for $\porder = 10$ on the polygonal case tilts for the last mesh, which is likely caused due to round-off errors when the results are close to machine precision. Such a situation is not observed for the quadrilateral mesh, which is because the Chebyshev basis functions are weight-orthogonal on square elements and so rounding errors are much smaller in this case. In  Fig.~\ref{fig:dg_convergence_elast_skeleton}, the integrated values of the jumps of the approximate solution along the mesh skeleton are presented. Similar curves but those that pertain to the boundary conditions are depicted in Fig.~\ref{fig:dg_convergence_elast_boundary}. The rates of convergence of the jumps of the approximate solution and violation of boundary conditions are similar to the global error rates. It shows that the final approximate solution is not strictly continuous but the level of the discontinuity is small and tends to zero with mesh refinement.
	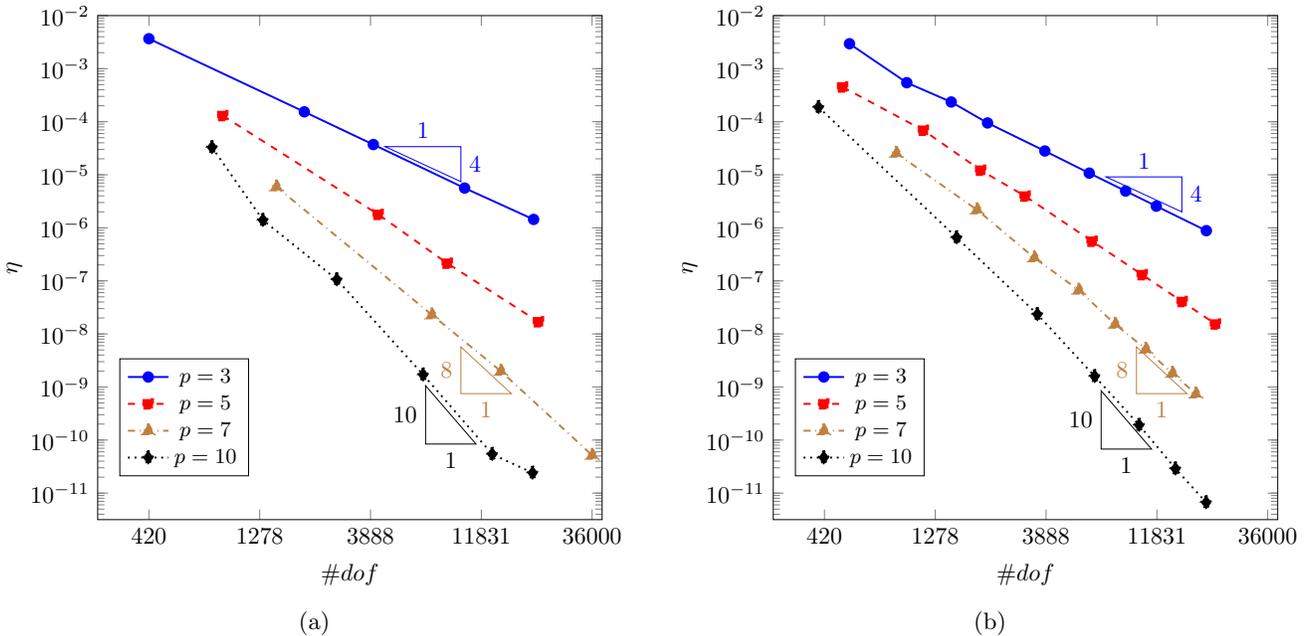
\begin{figure}[!htb]
		\centering
		\begin{subfigure}{0.49\textwidth}
			\centering
			\begin{tikzpicture}[scale=0.9]
				\begin{axis}[width=14cm, height=14cm,
					xlabel=${\#dof}$,
					ylabel=$\eta$,
					scale=0.6,
					ymode=log,
					xmode=log,
					xmin=10^(1.2), 	xmax=10^(2.3),
					ymin=10^(-11.5), ymax=10^(-2),
					legend style={at={(0.3,0.32)}},
					xtick = {20.493902, 35.749126, 62.353829, 108.770400, 189.736660},
					xticklabels = {420, 1278, 3888, 11831, 36000},
					]
					
					\addplot [solid, color=blue, thick, mark=*] 
					table [x index=3, y index=4] {figures/example3/conv2_pss_polyg_p3.txt};
					\addlegendentry{\small $p=3$}
					
					\addplot [dashed, color=red, thick, mark=square*] 
					table [x index=3, y index=4] {figures/example3/conv2_pss_polyg_p5.txt};
					\addlegendentry{\small $p=5$}
					
					\addplot [dashdotted, color=brown, thick, mark=triangle*, mark size=3pt]
					table [x index=3, y index=4] {figures/example3/conv2_pss_polyg_p7.txt};
					\addlegendentry{\small $p=7$}
					
					\addplot [dotted, color=black, thick, mark=diamond*, mark size=3pt]
					table [x index=3, y index=4] {figures/example3/conv2_pss_polyg_p10.txt};
					\addlegendentry{\small $p=10$}
					
					\logLogSlopeTriangle{0.72}{0.15}{0.74}{4}{blue};
					\logLogSlopeTTriangle{0.65}{0.1}{0.15}{10}{black};
					\logLogSlopeTTriangle{0.72}{0.1}{0.25}{8}{brown};			
				\end{axis}
			\end{tikzpicture}
			\subcaption{}
		\end{subfigure}
		\begin{subfigure}{0.49\textwidth}
			\centering
			\begin{tikzpicture}[scale=0.9]
				\begin{axis}[width=14cm, height=14cm,
					xlabel=${\#dof}$,
					ylabel=$\eta$,
					scale=0.6,
					ymode=log,
					xmode=log,
					xmin=10^(1.2), 	xmax=10^(2.3),
					ymin=10^(-11.5), ymax=10^(-2),
					legend style={at={(0.3,0.32)}},
					xtick = {20.493902, 35.749126, 62.353829, 108.770400, 189.736660},
					xticklabels = {420, 1278, 3888, 11831, 36000},
					]
					\addplot [solid, color=blue, thick, mark=*]  
					table [x index=3, y index=4] {figures/example3/conv2_pss_square_p3.txt};
					\addlegendentry{\small $p=3$}
					
					\addplot [dashed, color=red, thick, mark=square*] 
					table [x index=3, y index=4] {figures/example3/conv2_pss_square_p5.txt};
					\addlegendentry{\small $p=5$}

					\addplot [dashdotted, color=brown, thick, mark=triangle*, mark size=3pt]
					table [x index=3, y index=4]{figures/example3/conv2_pss_square_p7.txt};
					\addlegendentry{\small $p=7$}
					
					\addplot [dotted, color=black, thick, mark=diamond*, mark size=3pt]
					table [x index=3, y index=4]{figures/example3/conv2_pss_square_p10.txt};
					\addlegendentry{\small $p=10$}
					
					\logLogSlopeTriangle{0.81}{0.15}{0.68}{4}{blue};
					\logLogSlopeTTriangle{0.65}{0.1}{0.14}{10}{black};
					\logLogSlopeTTriangle{0.72}{0.1}{0.25}{8}{brown};	
					
				\end{axis}
			\end{tikzpicture}
			\subcaption{}
		\end{subfigure}
		
		\caption{Convergence study of the \PFDG{} method 
			over the domain for Example 4 on (a) polygonal meshes
			and (b) quadrilateral meshes. Scale is 
			$\frac{1}{2}\log$--$\log$.
		}
		\label{fig:dg_convergence_elast_domain}
	\end{figure}
	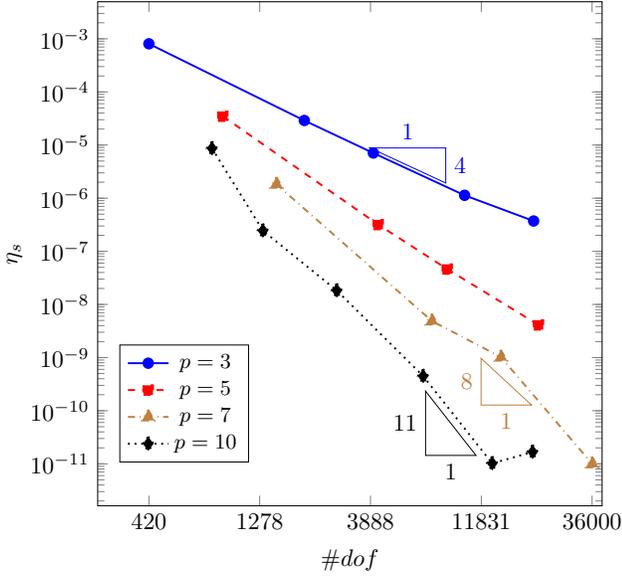
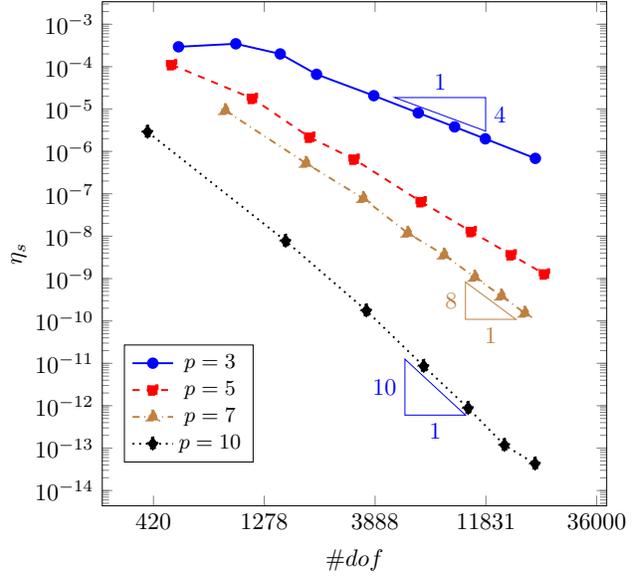
\begin{figure}[!htb]
		\centering
		\begin{subfigure}{0.49\textwidth}
			\centering
			\begin{tikzpicture}[scale=0.9]
				\begin{axis}[width=14cm, height=14cm,
					xlabel=${\#dof}$,
					ylabel=$\eta_s$,
					scale=0.6,
					ymode=log,
					xmode=log,
					xmin=10^(1.2), 	xmax=10^(2.3),
					legend style={at={(0.3,0.32)}},
					xtick = {20.493902, 35.749126, 62.353829, 108.770400, 189.736660},
					xticklabels = {420, 1278, 3888, 11831, 36000},
					]
					
					\addplot [solid, color=blue, thick, mark=*] 
					table [x index=3, y index=5] {figures/example3/conv2_pss_polyg_p3.txt};
					\addlegendentry{\small $p=3$}
					
					\addplot [dashed, color=red, thick, mark=square*] 
					table [x index=3, y index=5] {figures/example3/conv2_pss_polyg_p5.txt};
					\addlegendentry{\small $p=5$}
					
					\addplot [dashdotted, color=brown, thick, mark=triangle*, mark size=3pt]
					table [x index=3, y index=5] {figures/example3/conv2_pss_polyg_p7.txt};
					\addlegendentry{\small $p=7$}
					
					\addplot [dotted, color=black, thick, mark=diamond*, mark size=3pt]
					table [x index=3, y index=5] {figures/example3/conv2_pss_polyg_p10.txt};
					\addlegendentry{\small $p=10$}
					
					\logLogSlopeTriangle{0.69}{0.15}{0.71}{4}{blue};
					\logLogSlopeTTriangle{0.65}{0.1}{0.10}{11}{black};	
					\logLogSlopeTTriangle{0.76}{0.1}{0.20}{8}{brown};			
				\end{axis}
			\end{tikzpicture}
			\subcaption{}
		\end{subfigure}
		\begin{subfigure}{0.49\textwidth}
			\centering
			\begin{tikzpicture}[scale=0.9]
				\begin{axis}[width=14cm, height=14cm,
					xlabel=${\#dof}$,
					ylabel=$\eta_s$,
					scale=0.6,
					ymode=log,
					xmode=log,
					xmin=10^(1.2), 	xmax=10^(2.3),
					legend style={at={(0.3,0.32)}},
					xtick = {20.493902, 35.749126, 62.353829, 108.770400, 189.736660},
					xticklabels = {420, 1278, 3888, 11831, 36000},
					]
					\addplot [solid, color=blue, thick, mark=*]
					table [x index=3, y index=5] {figures/example3/conv2_pss_square_p3.txt};
					\addlegendentry{\small $p=3$}
					
					\addplot [dashed, color=red, thick, mark=square*]
					table [x index=3, y index=5] {figures/example3/conv2_pss_square_p5.txt};
					\addlegendentry{\small $p=5$}

					\addplot [dashdotted, color=brown, thick, mark=triangle*, mark size=3pt]
					table [x index=3, y index=5]{figures/example3/conv2_pss_square_p7.txt};
					\addlegendentry{\small $p=7$}
					
					\addplot [dotted, color=black, thick, mark=diamond*, mark size=3pt]
					table [x index=3, y index=5]{figures/example3/conv2_pss_square_p10.txt};
					\addlegendentry{\small $p=10$}
					
					\logLogSlopeTriangle{0.76}{0.18}{0.81}{4}{blue};
					\logLogSlopeTTriangle{0.60}{0.12}{0.18}{10}{blue};
					\logLogSlopeTTriangle{0.72}{0.1}{0.37}{8}{brown};	
					
				\end{axis}
			\end{tikzpicture}
			\subcaption{}
		\end{subfigure}
		
		\caption{Convergence study of the \PFDG{} method 
			along the mesh skeleton
			(discontinuity measure) for Example 4 on (a) polygonal meshes
			and (b) quadrilateral meshes. Scale is 
			$\frac{1}{2}\log$--$\log$.}
		\label{fig:dg_convergence_elast_skeleton}
	\end{figure}
	
	\subsection{Example 5: Cook’s membrane problem for hyperelasticity}\label{ssec:hyperelasticity}
	
\begin{figure}
		\centering
		\begin{subfigure}[t]{0.49\textwidth}
			\centering
			\begin{tikzpicture}[scale=0.1]
				\filldraw[line width=1.5pt,fill=gray!10] (0,0) -- ++(0,44) -- ++(48,16) -- ++(0, -16) -- cycle;
				\draw[-latex] (48+1.5,44) -- ++(0,16/5);
				\draw[-latex] (48+1.5,44+16/5) -- ++(0,16/5);
				\draw[-latex] (48+1.5,44+2*16/5) -- ++(0,16/5);
				\draw[-latex] (48+1.5,44+3*16/5) -- ++(0,16/5);
				\draw[-latex] (48+1.5,44+4*16/5) -- ++(0,16/5);
				\draw(48+1.5,44+16/2) node[right]{$p_0$};
				\filldraw[black] (48,44+16) circle (15pt) node[anchor=south west]{$A$};
				\draw[] (0,0)--++(-2,-4) ++(2,4+4) -- ++(-2,-4)
				++(2,4+4) -- ++(-2,-4)  ++(2,4+4) -- ++(-2,-4)  ++(2,4+4) -- ++(-2,-4)
				++(2,4+4) -- ++(-2,-4)  ++(2,4+4) -- ++(-2,-4)  ++(2,4+4) -- ++(-2,-4)
				++(2,4+4) -- ++(-2,-4)  ++(2,4+4) -- ++(-2,-4)  ++(2,4+4) -- ++(-2,-4)
				++(2,4+4) -- ++(-2,-4)  ++(2,4+4)
				; 
				
				\draw[] (0,0-5)
				-- ++(-1,-1) -- ++(2,2) ++(-1,-1)
				-- ++(48,0)
				-- ++(-1,-1) -- ++(2,2) ++(-1,-1)
				++(-48/2,0) node[anchor = north]{48 mm}
				;
				
				\draw[] (44+12,0)
				-- ++(-1,-1) -- ++(2,2) ++(-1,-1)
				-- ++(0,44)
				-- ++(-1,-1) -- ++(2,2) ++(-1,-1)
				-- ++(0,16)
				-- ++(-1,-1) -- ++(2,2) ++(-1,-1)
				++ (0,-16/2) node[anchor = north, rotate=90]{16 mm}
				++ (0,-16/2-44/2) node[anchor = north, rotate=90]{44 mm}
				;    
				\draw[](10,40) node[]{$\times$~1 mm};
				\draw[] (21,17) node[anchor=west]{$E = 500$ MPa}
				(21,17-5) node[anchor=west]{$\nu = 0.35$}
				(21,17-10) node[anchor=west]{$p_0 = 20$ MPa}
				;
			\end{tikzpicture}
			\subcaption{}\label{fig:ex6:cooksdomain-a}
		\end{subfigure}
		\begin{subfigure}[t]{0.49\textwidth}
			\centering
			\begin{tikzpicture}[scale=0.1] 	\input{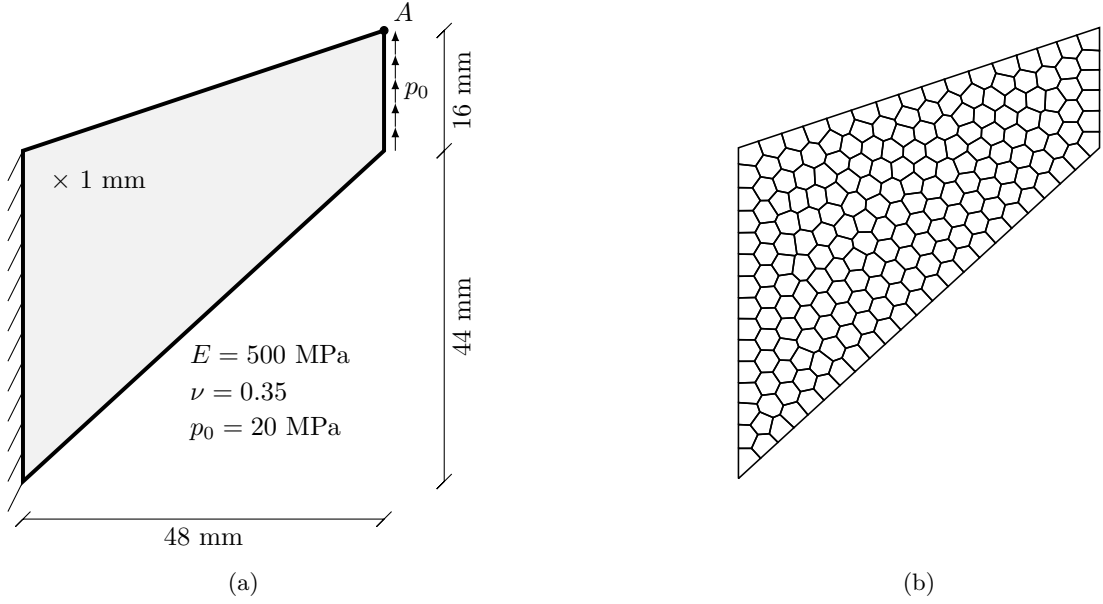}
				\draw[](0,-10);
			\end{tikzpicture}
			\subcaption{}\label{fig:ex6:cooksdomain-b}
		\end{subfigure}
		\caption{Cook's membrane problem (Example 5)
			with (a) geometry, material parameters and boundary conditions, and (b) polygonal mesh
			discretization of the domain.}
		\label{fig:ex6:cooksdomain}
	\end{figure}

	In this example, we solve the well-known Cook's membrane problem,~\cite{Cook1974} which combines bending and shear deformation. The cantilever is clamped on the left-hand side and is subjected to a constant shear loading on the right-hand side. The geometry, boundary conditions, and material parameters for this problem are shown in Fig.~\ref{fig:ex6:cooksdomain-a}.
	
	The Cook's membrane is analyzed under plane strain conditions. The hyperelastic model for the material is defined by the strain energy density function\cite{schroder2021selection}
	\begin{align}
		\psi = \frac{\mu}{2} \left(\Tr \RightCG -3 \right) + \frac{\lambda}{4} \left(J^2 -1 \right)
		- \left( \frac{\lambda}{2} + \mu \right) \ln J  \eqcomma
	\end{align}
	where $\lambda$ and $\mu$ are the Lam{\'e} parameters, $J$ is the determinant of the deformation gradient $\Jacob=\gradoper \displu + \unitmatrix$, and $\bm C = \Jacob\tran \bm \Jacob$ is the right Cauchy-Green tensor. The second Piola Kirchhoff stress tensor and the material tangent tensor are
	\begin{align}
		&
		\spk = 2\frac{\partial \psi}{\partial \RightCG} = 
		\frac{\lambda}{2} \left(J^2-1 \right) \bm C^{-1} + \mu \left( \unitmatrix - \RightCG^{-1} \right)  \eqcomma
		\\&
		\mattangent = 2\frac{\partial \spk}{\partial \RightCG} \eqcomma
	\end{align}
	where symbolic computations are used to compute the tensor $\mattangent$.
	
	The momentum balance, in the absence of the body forces, in the initial configuration can be expressed in terms of the first Piola-Kirchhoff stress tensor as:
	\begin{subequations}\label{eq:momentum_pk}
		\begin{align}
			&
			\diverg \fpk = \bm 0 \quad \text{in } \Domain \eqcomma
			\\&
			\fpk\scalarprod \VVnormal = \boundarysign{\spow} \quad \text{on } \Boundaryneum \eqcomma
			\\&
			\displu= \bm 0  \quad \text{on } \Boundarydirich \eqcomma
		\end{align}
	\end{subequations}
	where $\VVnormal$ denotes the unit outward normal. 
	
	The weak form of \eqref{eq:momentum_pk} is: Find $\displu \in \funspacecz$ such that
	\begin{align}\label{eq:momentum_pk_weak}
		\omint \gradoper \fvu   : \fpk  \dOm  - \gaintt{\Boundaryneum} \fvu\tran \boundarysign{\spow} \dGa = 0 \ \
		\forall \fvu \in \funspacecz 
		\eqdot
	\end{align}
	
	The Green-Lagrange strain tensor and its variation are
	\begin{align}
		&
		\GLstrain = \frac{1}{2} \left( \Jacob\tran\Jacob  - \unitmatrix \right) \eqcomma
		\\&
		\delta\GLstrain = \frac{1}{2} \left(\delta\Jacob\tran\Jacob + \Jacob\tran\delta\Jacob \right)
		,\quad
		\delta \Jacob = \frac{\partial \delta \displu}{\partial \bm X} = \frac{\partial \fvu }{\partial \bm X} = \gradoper \fvu
	\end{align}
	
	Using the relation between the first and the second Piola-Kirchhoff stress tensors, $\fpk = \Jacob \spk$, the inner product~\eqref{eq:momentum_pk_weak} is rewritten as follows:
	\begin{align}
		&
		\gradoper \fvu   : \fpk =   
		\gradoper \fvu  : \left(\Jacob\cdot\spk\right) = \left( \gradoper\fvu\tran \cdot \Jacob \right) : \spk
		\\&= 
		\frac{1}{2}  \left( \gradoper\fvu\tran \cdot \Jacob  + \Jacob\tran \gradoper\fvu \right) : \spk 
		+ \cancelto{0}{\frac{1}{2}  \left( \gradoper\fvu\tran \cdot \Jacob  - \Jacob\tran \gradoper\fvu \right): \spk}
		= 
		\delta\GLstrain(\fvu) : \spk \eqdot
	\end{align}
	
	The tensor $\spk$ is an unknown tensor and is updated using its increments, which can be expressed in the recurrence form:
	\begin{align}
		\spk = \spk  + \Delta \spk \eqdot
	\end{align}
	The tensor increment $\Delta\spk$ is expressed using a linearized relation.  There is no unique technique for linearization and in this paper, the following one is used:
	\begin{align}
		\Delta\spk \approx \mattangent \strain(\Delta\displu),
		\quad 
		\strain(\Delta\displu) = \frac{1}{2} \left( \gradoper \Delta\displu + 
		\gradoper \Delta\displu\tran \right) \eqdot
	\end{align}
	
	The nonlinear problem in \eqref{eq:momentum_pk_weak} is solved using an incremental iterative procedure, which is written in the following form:
	\begin{align}
		\omint 	\delta\GLstrain(\fvu) : \mattangent \strain(\Delta\displu) \dOm
		= \Lambda \gaintt{\Boundaryneum} \fvu\tran \boundarysign{\spow} \dGa - \omint \gradoper \fvu   : \fpk  \dOm 
		\eqdot
	\end{align}
	In this procedure, the external force multiplier $0< \Lambda <1$ is increased at every step. At every incremental step, equilibrium is found by updating the tensor $\fpk$ at every iterative step. The tangent tensor $\mattangent$ is updated at every incremental step.
	
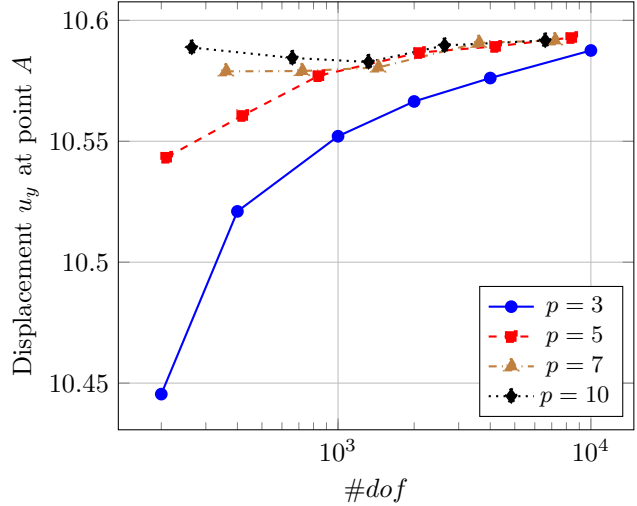
\begin{figure}\centering
		\begin{tikzpicture}
			\begin{axis}[
				grid,
				xlabel=$\#dof$,
				ylabel=Displacement $u_y$ at point $A$,
				legend style={at={(0.7,0.03)},anchor=south west},
				xmode = log,
				]

				\addplot [solid, color=blue, thick, mark=*]  
				table [x index=3, y index=4] {figures/example6_cook/cook2_p3.txt};
				\addlegendentry{\small $p=3$}
				
				\addplot [dashed, color=red, thick, mark=square*] 
				table [x index=3, y index=4] {figures/example6_cook/cook2_p5.txt};
				\addlegendentry{\small $p=5$}
				
				\addplot [dashdotted, color=brown, thick, mark=triangle*, mark size=3pt]
				table [x index=3, y index=4] {figures/example6_cook/cook2_p7.txt};
				\addlegendentry{\small $p=7$}
				
				\addplot [dotted, color=black, thick, mark=diamond*, mark size=3pt] 
				table [x index=3, y index=4] {figures/example6_cook/cook2_p10.txt};
				\addlegendentry{\small $p=10$}
				
			\end{axis}
		\end{tikzpicture}
		\caption{Convergence curves of 
			displacement $u_y$ at point $A$ for Example 5.}
		\label{fig:ex6:convergence}
	\end{figure}

\begin{figure}\centering
		\begin{subfigure}{0.32\textwidth}
			\includegraphics[height=0.8\textwidth]{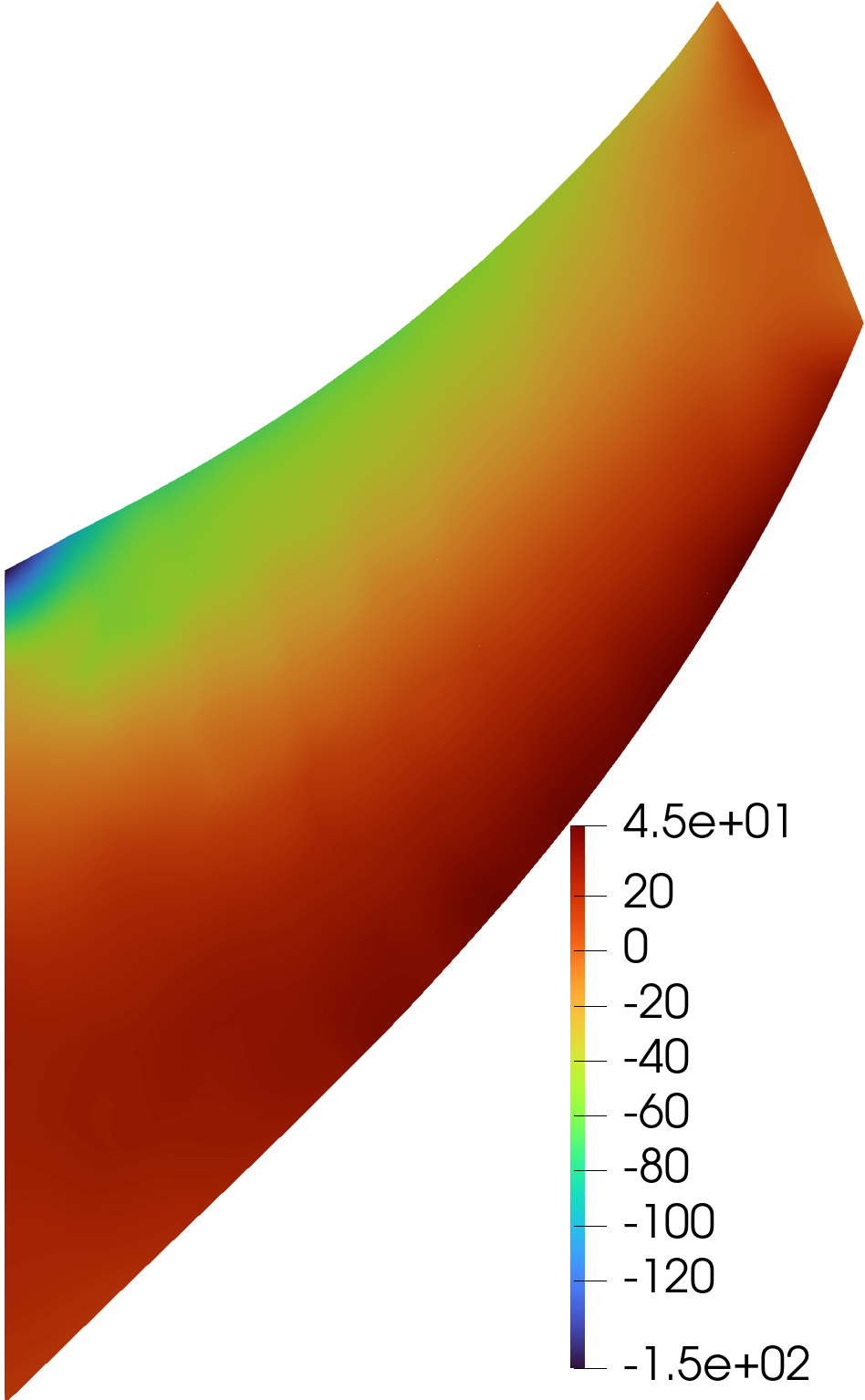}
			\subcaption{}
		\end{subfigure}
		\begin{subfigure}{0.32\textwidth}
			\includegraphics[height=0.8\textwidth]{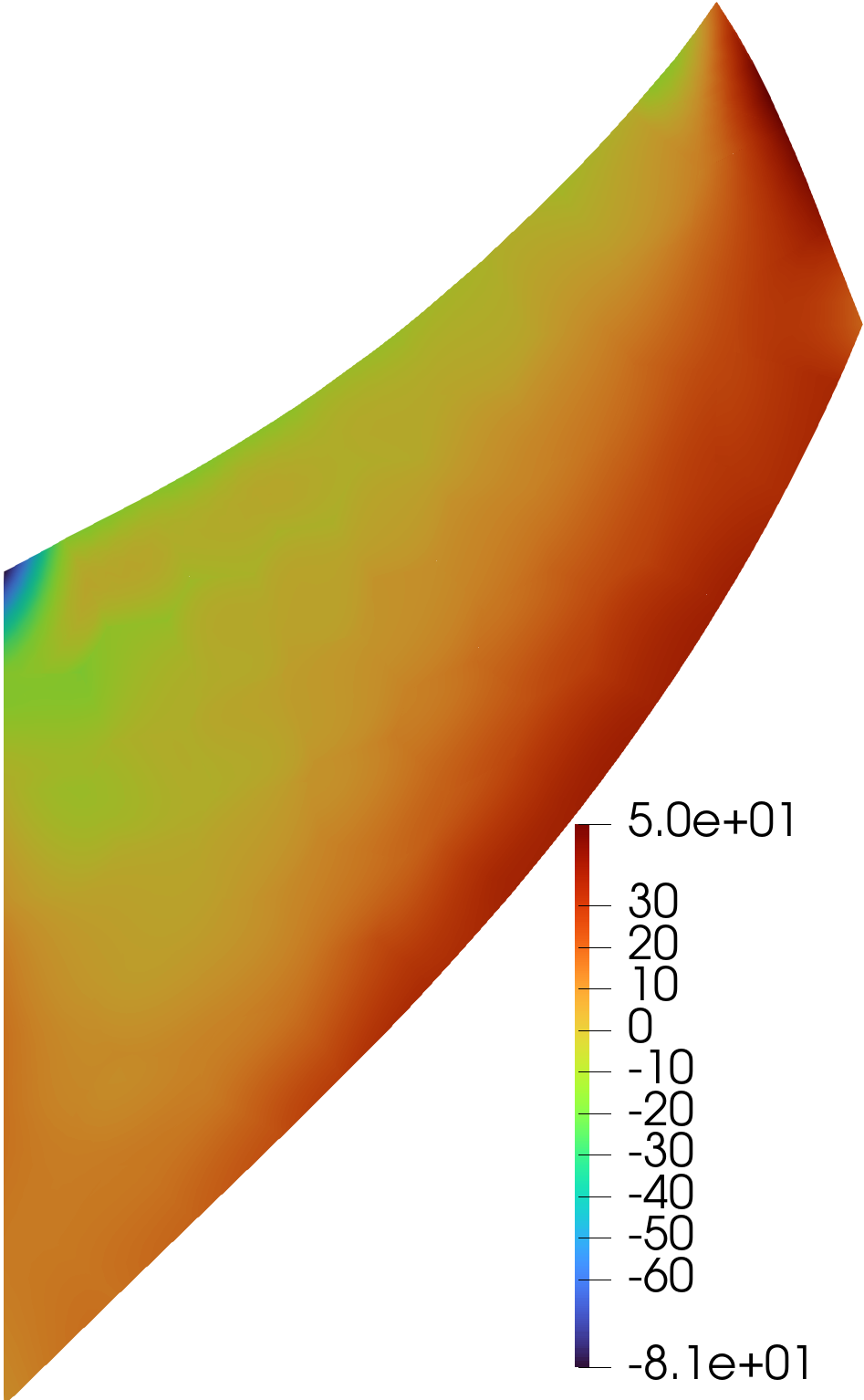}
			\subcaption{}
		\end{subfigure}
		\begin{subfigure}{0.32\textwidth}
			\includegraphics[height=0.8\textwidth]{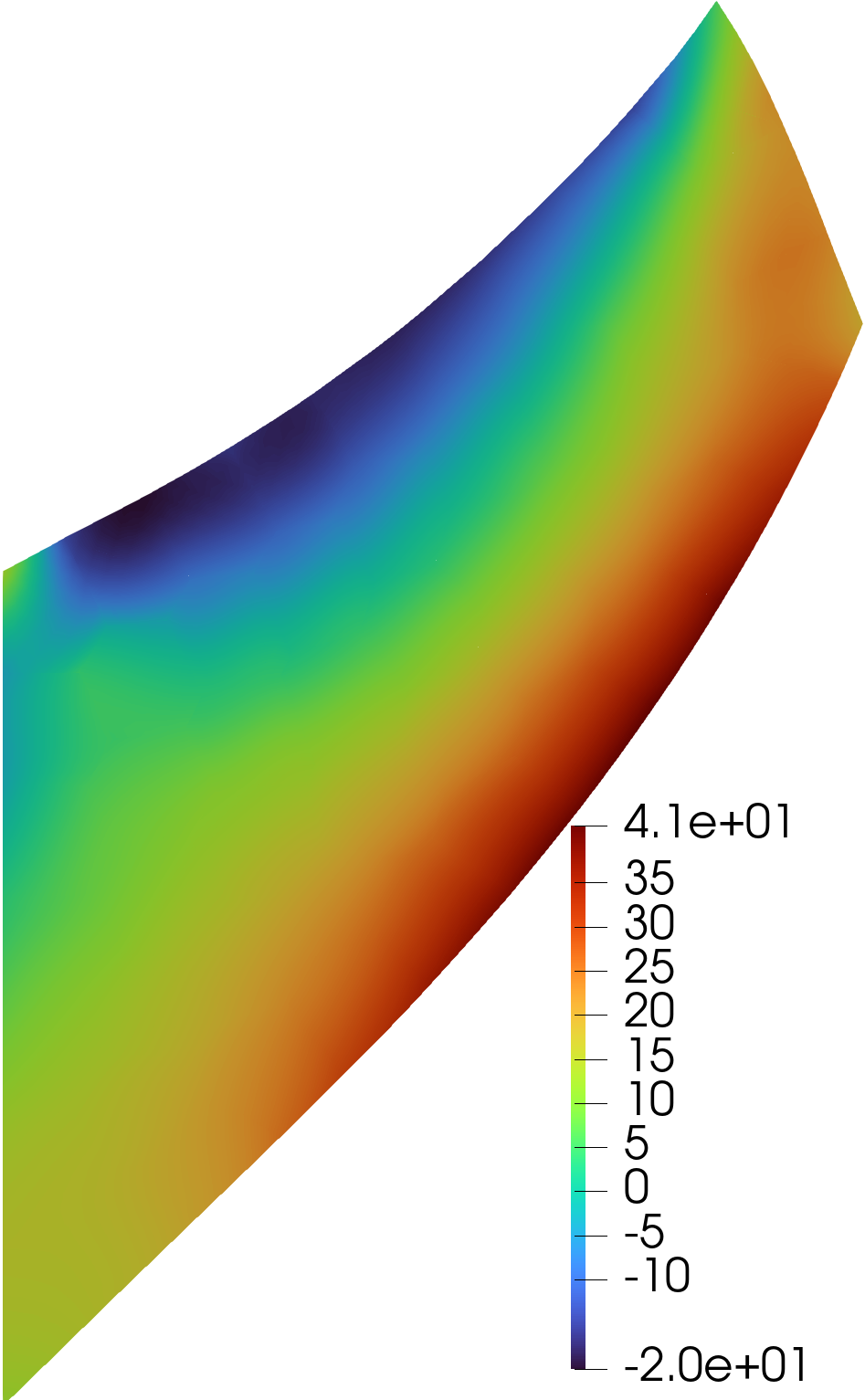}
			\subcaption{}
		\end{subfigure}
		
		\caption{Stress components for the Cook's membrane problem (a) $\sigma_{xx}$, (b) $\sig_{yy}$, (c) $\sig_{xy}$  for Example 5.}
		\label{fig:ex6_cooks_stress}
	\end{figure}

	Meshes with polygonal elements are used in the analysis (see Fig.~\ref{fig:ex6:cooksdomain-b}). The calculations are performed using approximation orders in the elements $p=$ 3, 5, 7, and 10. Numerical results are shown in Fig.~\ref{fig:ex6:convergence} in the form of convergence curves for the displacement $u_y$ at point $A$ for various $p$ and in Fig.~\eqref{fig:ex6_cooks_stress} the maps of the stress components, smoothed using the least-square technique. Convergence of the tip displacement to the value around $u_y=10.59$ is observed, which agrees with the results shown in Schr{\"o}der at al.\cite{schroder2021selection} 	In this example, the calculations are performed on the initial configuration, so the boundary and continuity condition components $\ttransfmatrix$ and $\tilde{\vm u}_b$ from~\eqref{eq:dg_discr_constr_ext_solve} are computed only once, just before the incremental-iterative procedure.
	
	\subsection{Example 6: Fourth-order biharmonic problem}\label{ssec:4order_example}
	
	In this example, we consider the fourth-order biharmonic boundary-value problem. The formulation for this problem is presented in Section~\ref{sec:pfdg_4order}. The exact solution is chosen to be trigonometric, the same as in \eqref{eq:bench_sin_exact}. The biharmonic problem presented in~\myeqref{eq:4problem} is chosen with the parameters defined as  trigonometric or polynomial functions:
	\begin{subequations}
		\begin{align}
			&
			\paramone =  \sin(1.3\pi\,x)\,\cos(1.8\pi\,y)+2 \eqcomma
			\\&
			\paramtwo =  4\,x^2\,y^2 + 2\,x y^2 + y +  \eqcomma
			\\&
			\paramthree = 3\,x y^2 - 2\,x,y+2  \eqdot
		\end{align}
	\end{subequations}

\begin{figure}\centering
		\begin{tikzpicture}
			\begin{axis}[width=14cm, height=14cm,
				xlabel=${\#dof}$,
				ylabel=$\eta$,
				scale=0.6,
				ymode=log,
				xmode=log,
				xtick = {36.331804, 52.829916, 76.830983, 111.731822, 162.480768},
				xticklabels = {1320, 2791, 5903, 12484, 26400},
				ytickten={-8,-6,...,1}
				]
				
				\addplot [solid, color=blue, thick, mark=*]  
				table [x index=3, y index=4] {figures/example4/ex3_poly_p4_vsin.txt};
				\addlegendentry{\small $p=4$}
				
				\addplot [dashed, color=red, thick, mark=square*] 
				table [x index=3, y index=4] {figures/example4/ex3_poly_p6_vsin.txt};
				\addlegendentry{\small $p=6$}
				
				\addplot [dashdotted, color=brown, thick, mark=triangle*, mark size=3pt]
				table [x index=3, y index=4] {figures/example4/ex3_poly_p8_vsin.txt};
				\addlegendentry{\small $p=8$}
				
				\addplot [dotted, color=black, thick, mark=diamond*, mark size=3pt] 
				table [x index=3, y index=4] {figures/example4/ex3_poly_p10_vsin.txt};
				\addlegendentry{\small $p=10$}
				
				\logLogSlopeTriangle{0.72}{0.15}{0.71}{5}{blue};
				\logLogSlopeTTriangle{0.37}{0.1}{0.41}{11}{black};
				
			\end{axis}
		\end{tikzpicture}
		
		\caption{Convergence study of the \PFDG{} 
			method for Example 6 on polygonal meshes. Scale is 
			$\frac{1}{2}\log$--$\log$.
		}
		\label{fig:dg_4convergence_domain_vsin}
	\end{figure}
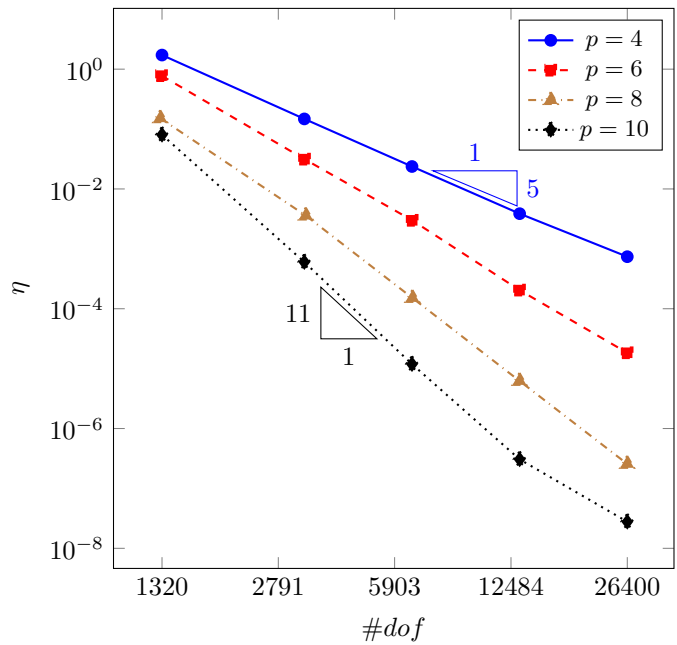

	The domain $\Domain=(-1,\,1)^2$  is covered by polygonal elements. Convergence curves for various approximation orders are presented in Fig.~\ref{fig:dg_4convergence_domain_vsin}. In all cases, optimal convergence rates are obtained, confirming the effectiveness of the proposed method.
	
	\subsection{Example 7: 3D benchmark problem}\label{ssec:3d_example}
	
	A Poisson problem in $\Domain=(-1,1)^3$ is solved with an exponential exact solution of the form:
	\begin{align}\label{eq:exp3d_example}
		u(x,y,z) = \sum_{i=1}^{50}  \exp \left(-a_i (x - x_i)^2 - b_i (y-y_i)^2 - c_i (y-z_i)^2\right) \eqcomma
	\end{align}
	where the positive values of the parameters $a$, $b$ and $c$  are chosen randomly and the points $(x_i, y_i, z_i)$ are randomly generated in the domain. This benchmark problem is the three-dimensional extension of the one presented in Section~\ref{ssec:poissons_dg_exp}.
	
\begin{figure}\centering
		\begin{tikzpicture}
			\begin{loglogaxis}[width=14cm, height=14cm,
				xlabel=${\#dof}$,
				ylabel=$\eta$,
				scale=0.6,
				legend style={at={(0.03,0.03)},anchor=south west},
				legend style={font=\Large},
				xtick = {13.177065, 16.938649, 21.774557, 27.991069, 35.981729},
				xticklabels = {2288, 4860, 10324, 21931, 46585},
				ytickten={-11,-9,...,-1}
				]

				\addplot [solid, color=blue, thick, mark=*] 
				table [x index=3, y index=4] {figures/example5_3d/ex_exp3d_p3.txt};
				\addlegendentry{\small $p=3$}
				
				\addplot [dashed, color=red, thick, mark=square*] 
				table [x index=3, y index=4] {figures/example5_3d/ex_exp3d_p4.txt};
				\addlegendentry{\small $p=4$}

				\addplot [dashdotted, color=brown, thick, mark=triangle*, mark size=3pt]
				table [x index=3, y index=4] {figures/example5_3d/ex_exp3d_p6.txt};
				\addlegendentry{\small $p=6$}

				\addplot [dotted, color=black, thick, mark=diamond*, mark size=3pt]
				table [x index=3, y index=4] {figures/example5_3d/ex_exp3d_p10.txt};
				\addlegendentry{\small $p=10$}
				
				\logLogSlopeTTriangle{0.68}{0.12}{0.10}{12}{black};
				\logLogSlopeTriangle{0.80}{0.18}{0.81}{4.5}{blue};
				
			\end{loglogaxis}
		\end{tikzpicture}
		\caption{Convergence study in the $L^2$ norm
			of the \PFDG{} 
			method for Example 7 on hexahedral meshes. Scale is 
			$\frac{1}{3}\log$--$\log$.} 
		\label{fig:example3d_convergence}
	\end{figure}
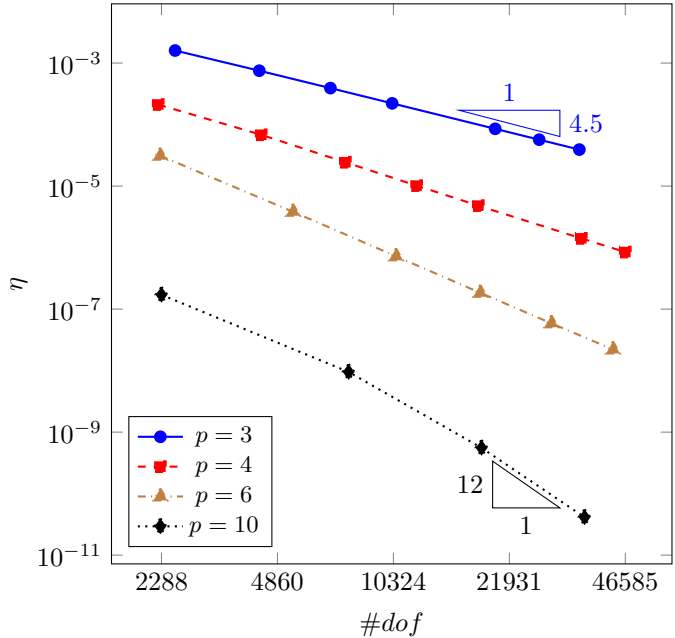

	The calculations are performed using structured hexahedral meshes for $p = 3, \ 4,\ 6,\ 10$. Convergence curves in the $L^2$ norm for different $p$ are presented in Fig.~\ref{fig:example3d_convergence}. The convergence plots once again reveal that the \PFDG{} method converges optimally for different approximation orders.
	
	\section{Conclusions}\label{sec:conclusions}
	In this paper, a new DG method is presented, in which no penalty or stability parameter is used, and hence the method is referred to as the penalty-free DG method. Generally, in DG methods, the global approximations are not continuous a priori on the mesh skeleton and continuity has to be enforced by penalty-like parameters. A similar scenario arises when imposing  Dirichlet boundary conditions. An alternative approach is proposed in this paper where the consistency of the approximation field, i.e., continuity on the mesh skeleton and satisfaction of the Dirichlet boundary conditions, is obtained by adding least-squares constraint relations. In the discrete version, the constraints require solving a singular system of equations. Such an approach allows for the construction of the discrete, continuous approximation fields over the whole domain, independent of the basis functions within the elements, and at the same time, the Dirichlet boundary conditions are satisfied. In the implementation, the consistent approximation field is augmented (especially needed on polygonal meshes), which allows for small levels of discontinuities on the mesh skeleton.
	
	The augmented consistent approximation field consists of basis functions in all finite elements, and the constraints are applied to the degrees of freedom, which guarantee the consistency of the approximation field. The constraints are expressed via a rectangular system of equations that is singular. In this paper, the method for solving such a singular system of equations is proposed, which can be solved globally, or the solution is constructed using a special assembly-free approach. The constraint relations are similar to the ones used in the DG method with Lagrange multipliers. However, in the case of the \pfdg{} method, there is no need to construct an additional space for the Lagrange multipliers.
	
	The \pfdg{} method is verified on several benchmark problems where various types of boundary-value problems are analyzed in 2D and for one problem in 3D. In each example, the convergence plots are presented for the approximation orders up to $\porder=10$, and optimal convergence rates in the $L^2$ norm and the energy seminorm are obtained. In finite-dimensional space, small jumps           (discontinuities) are observed in the approximate solution, but it has been shown via the examples that the jumps converge to zero with similar rates as the global error. The method can be especially effective in the nonlinear regime since the constraint problem is solved only once before the incremental iterative procedure. It is shown that the method works well for standard elliptic problems, such as the Poisson equation or linear and nonlinear elasticity. It was also tested for the biharmonic problem, where the approximation in a standard Galerkin method must be $C^1$ continuous. As part of future research, application of the \PFDG{} method to other classes of boundary-value problems will be pursued.
	
	\section{Data Availability Statement}
	The data that support the findings of this study are available from the corresponding author upon reasonable request.

	\appendix
	\section{Solver for the \PFDG{} method}\label{sec:dg_solver}
	
	The following  algebraic system of equations is considered:
	\begin{align}\label{eq:sq_basic_equation}
		\thematrix \rarg = \rvector \eqcomma
	\end{align}
	where $\thematrix$ is an $\matrows\times \matcols$ matrix, $\rvector$ is a $\matrows$-component column vector and $\rarg$ is an $\matcols$-component column vector. The vector $\rarg$ belongs to the set that consists of all possible solutions:
	\begin{align}
		\rarg \in \solsset = \left\{ \bm y:  \norm{\thematrix \bm y - \rvector} =  0 \right \} \eqdot
	\end{align}
	
	It is assumed that the system of equations~\eqref{eq:sq_basic_equation} is consistent, which means that $\solsset \neq \emptyset$. Since we have a linear, consistent system of equations, there exists one or infinite number of solutions. The matrix $\thematrix$ can be a square or rectangular matrix, i.e., $\matrows=\matcols$ or $\matrows<\matcols$. The matrix $\thematrix$ can be singular in the sense that the rows in the matrix are linearly dependent. The solution of the system of \eqref{eq:sq_basic_equation} can be unique, i.e., there exists only one $\rarg$, otherwise the set $\solsset$ consists of an infinite number of certain vectors.
	
	In this paper, such a general linear system of equations as defined above is solved using the Gauss-Jordan approach as presented in Valco and Vajda.\cite{valko1989advanced} Here the approach is presented in a more general way so that it can be applied to an arbitrarily large system of equations. For the sake of clarity, the nonsingular matrices are firstly considered,  which is extended to singular matrices in~\ref{ssec:singular}. 
	
	The procedure for solving the square (consistent) singular system of equations is presented in Farhat and Roux,\cite{Farhat1991} where the solver is applied to the floating subdomains as a result of domain partitioning in parallel calculations. This approach is referred to as the null-space method,~\cite{betsch2005discrete,rees2014null} which is used to solve the saddle-point problem.~\cite{scott2022null}. This method has successfully been applied to solve problems in incompressible linear and nonlinear elasticity.~\cite{munoz2008modelling} In this approach, the singular stiffness matrix is partitioned concerning the principal and redundant quantities, which is performed during the factorization process. Then the null space of the stiffness matrix can be identified, which contains the rigid-body modes of the subdomain. Then the final solution is presented as the linear combination of the null space basis vectors plus the vector constructed as the right-hand side vector multiplied by the pseudo-inverse stiffness matrix. However, the procedure presented here is more general, since it can also handle rectangular and singular matrices.
	
	The rectangular system of equations can be solved using the $LQ$ factorization, which is presented in Slobodkins and Tausch.\cite{slobodkins2023node} In this method, the matrix is decomposed as $\thematrix = \begin{bmatrix} \vm L & \vm 0\end{bmatrix} \begin{bmatrix}\vm Q^1 & \vm Q^2\end{bmatrix}\tran$, where the matrix $\vm L$ is lower triangular and matrix $\vm Q$ consists of orthogonal vectors where the vectors in $\vm Q^2$ form the basis of the null space of matrix $\thematrix$. Another approach that can be applied to solve the linear equation with a rectangular matrix is the singular value decomposition (SVD),\cite{andrilli2022elementary} and after some modifications, the SVD can also handle the solution of the singular system of equations.
	
	\subsection{Global approach}\label{sssec:glob}
	The variables in $\rarg$ are divided into basic variables vector $\rarg_b$ and free variables vector $\rarg_f$.  For the sake of clarity, it can be assumed without loss of generality, that basic and free variables are contiguous in the $\rarg$ vector. The matrix $\thematrix$ is decomposed on:
	\begin{align}\label{eq:seteq_divided}
		\rarg = \begin{bmatrix} \rarg_b \\ \rarg_f \end{bmatrix}
		\quad \Rightarrow \quad
		\begin{bmatrix}	\thematrix_b & \thematrix_f \end{bmatrix}
		\begin{bmatrix} \rarg_b \\ \rarg_f \end{bmatrix} = \rvector \eqdot
	\end{align}
	
	The division into the basic and free parts is done in such a way that the matrix $\thematrix_b$ is square and well-conditioned. On multiplication both sides of \eqref{eq:seteq_divided}  by $\thematrix_b^{-1}$, we get
	\begin{align}\label{eq:seteq_divided2}
		\begin{bmatrix}	\eye & \tilde{\thematrix}_f \end{bmatrix}
		\begin{bmatrix} \rarg_b \\ \rarg_f \end{bmatrix} = \tilde{\rvector} \eqcomma
	\end{align}
	where
	\begin{align}\label{eq:tilde_defs}
		\tilde{\thematrix}_f = \thematrix_b^{-1} \thematrix_f \quad \text{and} \quad
		\tilde{\rvector} = \thematrix_b^{-1} \rvector \eqdot
	\end{align}

	The system of equations shown in \eqref{eq:seteq_divided2}, is derived by multiplication  both sides of \eqref{eq:seteq_divided} by $\thematrix_b^{-1}$. However, the same results can be obtained using the Gauss-Jordan (GJ) elimination procedure. The GJ approach can be applied when the leading components of the matrix $\thematrix_b$ are located on the diagonal. In other cases the equations in~\eqref{eq:seteq_divided} can be rearranged to place the leading components on the diagonal or the leading component can be searched in every row of the matrix $\thematrix_b$. In this approach, the matrix $\thematrix_b$ is substituted by the matrix $\bar{\eye}$ which is the zero-one matrix, but the ones are not located on the diagonal. In this situation, it can be observed that
	\begin{align}
		\begin{bmatrix}	\bar\eye & \bar{\thematrix}_f \end{bmatrix}
		\begin{bmatrix} \rarg_b \\ \rarg_f \end{bmatrix} = \bar{\rvector}  \,\Bigr|\, 	\bar\eye\tran
		\quad \rightarrow \quad
		\begin{bmatrix}	\eye & \tilde{\thematrix}_f \end{bmatrix}
		\begin{bmatrix} \rarg_b \\ \rarg_f \end{bmatrix} = \tilde{\rvector} \eqdot
	\end{align}

	Using \eqref{eq:seteq_divided2} the basic variables can be expressed by the free variables:
	\begin{align}\label{eq:rbasic_rfree}
		\rarg_b = -\tilde{\thematrix_f}\rarg_f + \tilde{\rvector} \eqdot
	\end{align}
	Using the relation in \eqref{eq:rbasic_rfree} the vector $\rarg$ can be expressed by the free variable vector $\freevar$:
	\begin{align}\label{eq:solution}
		&\rarg = \transfmatrix\freevar + \transfvector \eqcomma
		\intertext{where}
		\label{eq:sols_defs}
		&\transfmatrix = \begin{bmatrix} -\tilde{\thematrix_f} \\ \eye \end{bmatrix}
		,\qquad
		\transfvector = \begin{bmatrix} \tilde{\rvector} \\ \vm 0 \end{bmatrix}
		,\qquad
		\freevar = \rarg_f \eqdot
	\end{align}
	
	\begin{theorem}
		If a  vector $\vm y$ is defined as $\vm y = \transfmatrix\vm z + \transfvector$, then
		\begin{align*}
			\vm y \in \solsset \ \ \forall \vm z \eqdot
		\end{align*}
	\end{theorem}
	
	\begin{proof}
		When $\vm y \in \solsset$ that means that the vector satisfies the equation 	~\eqref{eq:sq_basic_equation}. Let us suppose that there exists a vector $\vm z$ for which 	$\vm y \notin \solsset$, then there exists a nonzero vector $\residvector$  such that
		\begin{align}
			\residvector = \thematrix \vm y - \rvector = 
			\thematrix \transfmatrix \vm z + \thematrix\transfvector - \rvector \eqdot
		\end{align}
		Using the definitions of the vector $\vm y$ and in the \eqref{eq:tilde_defs} and~\eqref{eq:sols_defs} and applying the matrix division as in \eqref{eq:seteq_divided}, we have
		\begin{align}
			\residvector = \begin{bmatrix} \thematrix_b & \thematrix_f \end{bmatrix}
			\begin{bmatrix}
				-\thematrix_b^{-1} \thematrix_f \\ \eye 
			\end{bmatrix} \vm z
			+ \begin{bmatrix} \thematrix_b & \thematrix_f \end{bmatrix} 
			\begin{bmatrix} \thematrix_b^{-1} \rvector \\ \bm 0 \end{bmatrix} - \rvector \eqcomma
		\end{align}
		
		\begin{align}
			\residvector = (- \thematrix_f + \thematrix_f) \vm z + \rvector - \rvector
			\quad \Rightarrow \quad \residvector = \vm 0 \eqdot
		\end{align}
		
		It was shown that for every vector $\vm z$ the residual vector is zero.
		
	\end{proof}
	
	The solution in \eqref{eq:solution} represents the hypersurface in which $\transfvector$ represents the point on that surface and the vectors of the matrix $\transfmatrix$ are the basis vectors for that surface.  But in the case when $\matrows=\matcols$, the number of columns of the matrix $\transfmatrix$ is reduced to zero and the vector $\transfvector$ is the unique solution of the problem in \eqref{eq:sq_basic_equation}. When $\matrows<\matcols$ the vector $\rarg$ represents the points belonging to the hypersurface defined by the pair $(\transfmatrix,\transfvector)$. In this hypersurface, a vector $\rarg_m$ can be distinguished that is the shortest among all vectors in $\solsset$
	\begin{align}
		\norm{\rarg_m} = \min_{\rarg\in\solsset} \norm{\rarg} \eqdot
	\end{align}
	In order to identify the $\rarg_m$ vector the relevant $\vm z_m$ has to be found by solving the following minimization problem:
	\begin{align}
		\vm z_m \,: f(\vm z_m) = \min_{\vm z} f(\vm z) \eqcomma
		\intertext{where}
		f(\vm z) = \left(\transfmatrix\vm z + \transfvector \right)\tran \left(\transfmatrix\vm z + \transfvector \right) \eqcomma
	\end{align}
	which leads to the following algebraic system of equations:
	\begin{align}
		\transfmatrix\tran\transfmatrix \vm z_m = -\transfmatrix\tran \transfvector \eqdot
	\end{align}
	On solving for $\vm z_m$ the vector $\rarg_m$ is then given by
	\begin{align}
		\rarg_m = \transfmatrix \vm z_m + \transfvector \eqdot
	\end{align}

	\subsection{Recursive approach}\label{ssec:iterative}
	The approach presented in the previous section can be applied to obtain the recursive solution of the set of equations~\eqref{eq:sq_basic_equation}. For the recursive approach, this set of equations is divided into $\nsubdivision$ subproblems. For the sake of clarity, firstly let $\nsubdivision=2$:
	\begin{align}
		\begin{bmatrix}	\thematrix_1 \\ \thematrix_2 \end{bmatrix}	\rarg = 
		\begin{bmatrix} \rvector_1 \\ \rvector_2 \end{bmatrix} \eqcomma
	\end{align}
	where $\thematrix_i$ and $\rvector_i$ are the appropriate parts of the matrix $\thematrix$ and $\rvector$, respectively. Now, the problem defined in \eqref{eq:sq_basic_equation} is reformulated: find the vector $\rarg$ that satisfy the following equations:
	\begin{subequations}\label{eq:basic_equation_2split}
		\begin{align}
			&
			\label{eq:basic_equation_2split_a}
			\thematrix_1 \rarg = \rvector_1 \eqcomma
			\\&
			\label{eq:basic_equation_2split_b}
			\thematrix_2 \rarg = \rvector_2 \eqdot
		\end{align}
	\end{subequations}
	The solution of \eqref{eq:basic_equation_2split_a} leads to the following
	\begin{align}\label{eq:split2_sol1}
		\rarg =\transfmatrix_1 \freevar_1 + \transfvector_1 \eqdot
	\end{align}
	When applying this solution to~\eqref{eq:basic_equation_2split_b}, we have
	\begin{align}\label{eq:split2_sol2}
		\thematrix\transfmatrix_1 \freevar_1 = \rvector_2 - \thematrix_2\transfvector_2 \eqdot
	\end{align}
	The solution of \eqref{eq:split2_sol2} leads to
	\begin{align}
		\freevar_1 = \transfmatrix_2 \freevar_2 + \transfvector_2 \eqdot
	\end{align}
	When applying this equation to~\eqref{eq:split2_sol1} we have the final solution
	\begin{align}
		\rarg = \transfmatrix_1\transfmatrix_2\freevar_2 + \transfmatrix_1\transfvector_2 + \transfvector_1 \eqdot
	\end{align}

	In a case when there are multiple divisions of the set of equations~\eqref{eq:sq_basic_equation} then it is convenient to show the construction of the matrix $\transfmatrix$ and the vector $\transfvector$ using the recursive procedure:
	\begin{align}
		\begin{aligned}
			&\transfvector = \transfmatrix \transfvector_i + \transfvector \eqcomma
			\\
			&\transfmatrix = \transfmatrix\transfmatrix_i \eqcomma
		\end{aligned}
	\end{align}
	where $\transfvector_i$ and $\transfmatrix_i$ come from the $i$-th block of equations and the initial values are $\transfvector=\vm 0$ and $\transfmatrix=\eye$.
	
	The recursive procedure presented in this section has the same complexity as the global approach in \eqref{sssec:glob}, i.e., the number of multiplications in these two approaches is the same. However, the recursive approach can be applied when the system of equations is constructed by segments. Then the solution can be applied to each of the segments sequentially without the need to construct the full set of equations.

	\subsection{Singular system of equations}\label{ssec:singular}
	
	When the matrix $\thematrix$ in the problem~\eqref{eq:sq_basic_equation}, then the dependent equations should be removed from the system of equations, and then the solution shown in \eqref{eq:solution} can be achieved. The linearly dependent equations can be identified when applying the Gauss-Jordan elimination procedure to calculate the quantities in \eqref{eq:rbasic_rfree}.  The Gauss-Jordan procedure is performed using the leading position in each row to locate the pivot variable. When reaching a linear dependent equation all the elements in that record become zero, so this equation is deleted from the set of equations. There may be multiple such equations, and each time these equations are eliminated.

	\section{\PFDG{} in one dimension}\label{sec:1dexample}
	We provide the step-by-step procedure for the \PFDG{} method on two Poisson problems in 1D. In the first case, $\porder=1$ and $\bporder=1$ are chosen; in the second case, $\porder=2$ and $\bporder=1$ are selected. To clarify the exposition, in both cases, we divide the unit interval into two finite elements $[0,\,\frac{1}{2}]$ and $[\frac{1}{2},\,1]$. 
	
	Monomial basis functions are chosen to construct the approximation, which is defined in the element's local coordinates scaled by the element length. For first- and second-order finite element approximations, the $e$-th element basis are:
	\begin{align}
		\Gshapev^e = \begin{bmatrix} 1 & \frac{2}{L^e} \left(x - x_m^e \right) 
		\end{bmatrix} 
		\,,\quad
		\Gshapev^e = \begin{bmatrix} 1 & \frac{2}{L^e} \left(x - x_m^e \right) 
			& \frac{4}{{L^e}^2} \left(x - x_m^e \right)^2 \end{bmatrix} \eqcomma
	\end{align}
	where $L^e=\frac{1}{2}$ is the length of the element and $x_m^e$ is the mid-point of the element, i.e., $x_m^1=\frac{1}{4}$, $x_m^2 = \frac{3}{4}$.
	
	\subsection{Poisson problem with affine solution}
	We consider the following 1D Poisson problem: 
	\begin{subequations}\label{eq:1d_poisson}
		\begin{align}
			&
			\label{eq:1d_poisson_a}
			\temperature'' = 0 \quad \text{in } \Domain=(0,\,1),
			\\&
			\label{eq:1d_poisson_b}
			u(0)=2,\quad u(1)=5 \eqcomma
		\end{align} 
	\end{subequations}
	which has the exact solution $u(x) = 3x+2$. In this case $\porder=\bporder=1$ so the subscript in the approximation matrix $\Gshapev$ is omitted for clarity. The approximation of $\temperature$ and its derivative in  $[0,\,1]$  is: 
	\begin{align}\label{eq:1d_ex_approx}
		\temperature = \Gshapev \dftemp,\quad
		\temperature = \begin{cases}
			\begin{bmatrix} 1 & 4x-1 \end{bmatrix} \begin{bmatrix} \dftempsc_1 \\ \dftempsc_2\end{bmatrix}
			& \text{if } 0 \le x < 0.5
			\\
			\begin{bmatrix} 1 & 4x-3 \end{bmatrix} \begin{bmatrix} \dftempsc_3 \\ \dftempsc_4\end{bmatrix}
			& \text{if } 0.5 \le x \le 1 
		\end{cases} \eqcomma
	\end{align}

	\begin{align}\label{eq:ex_1d_approx}
		\temperature' = \Gshapev' \dftemp,\quad
		\temperature' = \begin{cases}
			\begin{bmatrix} 0 & 4 \end{bmatrix} \begin{bmatrix} \dftempsc_1 \\ \dftempsc_2\end{bmatrix}
			& \text{if } 0 \le x < 0.5
			\\
			\begin{bmatrix} 0 & 4 \end{bmatrix} \begin{bmatrix} \dftempsc_3 \\ \dftempsc_4\end{bmatrix}
			& \text{if } 0.5 \le x \le 1 
		\end{cases}  \eqdot
	\end{align}

	The constrained problem in \eqref{eq:dg_constr_problem} requires that \eqref{eq:discrete_conds_ext} has to be firstly solved to yield $\ttransfmatrix$ and $\ttransfvectorb$. The values of the approximation matrix $\Gshapev$ on the boundaries and its jump on the skeleton are:
	\begin{align}
		&\Gshapev(0) = \begin{bmatrix}1 & -1 & 0 & 0 \end{bmatrix},\quad
		\Gshapev'(0) = \begin{bmatrix}0& 4 & 0 & 0 \end{bmatrix} \eqcomma
		\\&
		\Gshapev(1) = \begin{bmatrix}0 & 0 & 1 & 1 \end{bmatrix},\quad
		\Gshapev'(1) = \begin{bmatrix}0& 0 & 0 & 4 \end{bmatrix} \eqcomma
		\\&
		\inbbr{\Gshapev(0.5)} = \begin{bmatrix} 0 & 0 & 1 & -1\end{bmatrix} - \begin{bmatrix} 1 & 1 & 0 & 0\end{bmatrix} = \begin{bmatrix}-1 & -1 & 1 & -1 \end{bmatrix}  \eqcomma
		\\&
		\iinbbr{\Gshapev'(0.5)} =  \frac{1}{2}\left(\begin{bmatrix} 0 & 4 & 0 & 0\end{bmatrix} + \begin{bmatrix} 0 & 0 & 0 & 4\end{bmatrix}\right) = \begin{bmatrix}0 & 2 & 0 & 2 \end{bmatrix}  \eqdot
	\end{align}
	The ${\tcmtx}$ matrix and $\tcmtvb$ vector in this case are:
	\begin{subequations}\label{eq:1d_example_constr}
		\begin{align}\label{eq:1d_example_constr_a}
			\begin{aligned}
				\tcmtx =& 	\inbbr{\Gshapev(0.5)}\tran\inbbr{\Gshapev(0.5)} + \Gshapev(0)\tran\Gshapev(0) + \Gshapev(1)\tran\Gshapev(1)
				\\=& \begin{bmatrix}
					1 &    1 &   -1  &   1\\
					1  &   1 &   -1  &   1\\
					-1  &  -1 &    1  &  -1\\
					1  &   1  &  -1   &  1
				\end{bmatrix} 
				+
				\begin{bmatrix}
					1 & -1 & 0 & 0\\
					-1 & 1 & 0 & 0\\
					0 & 0 & 0 & 0\\
					0 & 0 & 0 & 0
				\end{bmatrix}
				+
				\begin{bmatrix}
					0 & 0 & 0 & 0\\
					0 & 0 & 0 & 0\\
					0 & 0 & 1 & 1\\
					0 & 0 & 1 & 1
				\end{bmatrix}
				=
				\begin{bmatrix}
					2 & 0 & -1 & 1\\
					0 & 2 & -1 & 1\\
					-1 &-1 &  2 & 0\\
					1 & 1 &  0 & 2
				\end{bmatrix} 
			\end{aligned} \eqcomma
		\end{align}
		
		\begin{align}\label{eq:1d_example_constr_b}
			\tcmtvb =  \Gshapev(0)\tran a + \Gshapev(1)\tran b = \begin{bmatrix} 2 \\-2\\5\\5 \end{bmatrix} \eqdot
		\end{align}
	\end{subequations}
	
	The solution of the system of equations $(\tcmtx,\, \tcmtvb)$ gives the pair $(\ttransfmatrix,\, \ttransfvectorb)$, which are
	\begin{align}
		\ttransfmatrix = \begin{bmatrix} -1\\-1\\-1\\1\end{bmatrix} ,
		\quad
		\ttransfvectorb =  \begin{bmatrix} \frac{7}{2}\\\frac{3}{2}\\5\\0\end{bmatrix} \eqdot
	\end{align}
	
	The step-by-step solution procedure of the constraints problem is presented in Appendix~\ref{sec:1d_constr_example}.
	
	The matrix $\vm{K}$, vector $\vm{f}$ in \eqref{eq:dg_constr_problem}  are defined as:
	\begin{subequations}
		\begin{align}
			\begin{aligned}
				\vm{K} =& \mint_0^1 {\Gshapev'}\tran\Gshapev' \intend x  + \inbbr{\Gshapev(0.5)} \iinbbr{\Gshapev'(0.5)} 
				+ \Gshapev(0)\Gshapev'(0) - \Gshapev(1)\Gshapev'(1)
				\\=&
				\begin{bmatrix}
					0 & 0 & 0 & 0\\
					0 & 8 & 0 & 0\\
					0 & 0 & 0 & 0\\
					0 & 0 & 0 & 8\\
				\end{bmatrix}
				+
				\begin{bmatrix}
					0 & -2 & 0 & -2\\
					0 & -2 & 0 & -2\\
					0 & 2 & 0 & 2\\
					0 & -2 & 0 & -2\\
				\end{bmatrix}
				+
				\begin{bmatrix}
					0 & 4 & 0 & 0\\
					0 & -4 & 0 & 0\\
					0 & 0 & 0 & 0\\
					0 & 0 & 0 & 0\\
				\end{bmatrix}
				-
				\begin{bmatrix}
					0 & 0 & 0 & 0\\
					0 & 0 & 0 & 0\\
					0 & 0 & 0 & 4\\
					0 & 0 & 0 & 4\\
				\end{bmatrix}
				\\=&
				\begin{bmatrix}
					0 & 2 & 0 & -2\\
					0 & 2 & 0 & -2\\
					0 & 2 & 0 & -2\\
					0 & -2 & 0 & 2\\
				\end{bmatrix} \eqcomma
			\end{aligned} 
		\end{align}
		
		\begin{align}
			\vm{f} = -\mint_0^1 \Gshapev\tran 0 \intend x = \begin{bmatrix}0 \\ 0 \\ 0 \\0\end{bmatrix} \eqdot
		\end{align}
	\end{subequations}
	Now, the linear system of equations given in~\eqref{eq:dg_constr_problem} has the following form:
	\begin{align}
		\begin{bmatrix} 16 \end{bmatrix} \dftempfree = \begin{bmatrix} 12 \end{bmatrix}
		\quad\rightarrow\quad \dftempfree = \begin{bmatrix} 0.75\end{bmatrix} \eqdot
	\end{align}
	Finally, on using \eqref{eq:dg_constr} we have the vector $\dftemp$:
	\begin{align}
		\dftemp =  \begin{bmatrix} -1\\-1\\-1\\1\end{bmatrix} \begin{bmatrix} \frac{3}{4} \end{bmatrix} + 
		\begin{bmatrix}  \frac{7}{2}\\\frac{3}{2}\\5\\0\end{bmatrix} = \begin{bmatrix} \frac{11}{4}\\\frac{3}{4}\\\frac{17}{4}\\\frac{3}{4}\end{bmatrix} \eqdot
	\end{align}
	On applying the vector $\dftemp$ to the approximation in \eqref{eq:ex_1d_approx} we recover the exact solution (affine function) of the problem.

	\subsection{Poisson problem with quadratic solution}
	
	In this case, we consider the following 1D Poisson problem: 
	\begin{subequations}\label{eq:1d2p_poisson}
		\begin{align}
			&
			\label{eq:1d2p_poisson_a}
			\temperature'' = 2 \quad \text{in } \Domain=(0,\,1),
			\\&
			\label{eq:1d2p_poisson_b}
			u(0)=2,\quad u(1)=5 \eqcomma
		\end{align} 
	\end{subequations}
	which has the exact solution $u(x) = x^2 + 2x+2$. The approximation of $\temperature$ and its derivative in  $[0,\,1]$  is: 
	\begin{align}\label{eq:1d2p_ex_approx}
		\temperature = \Gshapev_\porder \dftemp,\quad
		\temperature = \begin{cases}
			\begin{bmatrix} 1 & 4x-1  & 16x^2 - 8x + 1
			\end{bmatrix} 
			\begin{bmatrix} \dftempsc_1 \\ \dftempsc_2 \\ \dftempsc_3\end{bmatrix}
			& \text{if } 0 \le x < 0.5
			\\
			\begin{bmatrix} 1 & 4x-3  & 16x^2 - 24x + 9	\end{bmatrix} 
			\begin{bmatrix} \dftempsc_4 \\ \dftempsc_5  \\ \dftempsc_6\end{bmatrix}
			& \text{if } 0.5 \le x \le 1 
		\end{cases} \eqcomma
	\end{align}

	\begin{align}\label{eq:ex_1d2p_approx}
		\temperature' = \Gshapev_\porder' \dftemp,\quad
		\temperature' = \begin{cases}
			\begin{bmatrix} 0 & 4 & 32x-8 \end{bmatrix} 
			\begin{bmatrix} \dftempsc_1 \\ \dftempsc_2 \\ \dftempsc_3\end{bmatrix}
			& \text{if } 0 \le x < 0.5
			\\
			\begin{bmatrix} 0 & 4& 32x-24  \end{bmatrix} 
			\begin{bmatrix} \dftempsc_4 \\ \dftempsc_5  \\ \dftempsc_6\end{bmatrix}
			& \text{if } 0.5 \le x \le 1 
		\end{cases}  \eqdot
	\end{align}
	
	The constrained problem in~\eqref{eq:dg_constr_problem} require that~\eqref{eq:discrete_conds_ext} has to be firstly solved to yield $\ttransfmatrix$ and $\ttransfvectorb$. The values of the approximation matrix $\Gshapev_\porder$ and $\Gshapev_\bporder$ on the boundaries and its jump on the skeleton are:
	\begin{align}
		&\Gshapev_\porder(0) = \begin{bmatrix}1 & -1 & 1 & 0 & 0 & 0 \end{bmatrix},\quad
		\Gshapev_\porder'(0) = \begin{bmatrix}0& 4 & -8 & 0 & 0 &0 \end{bmatrix} \eqcomma
		\\&
		\Gshapev_\porder(1) = \begin{bmatrix}0 & 0 & 0 & 1 & 1 & -1 \end{bmatrix},\quad
		\Gshapev_\porder'(1) = \begin{bmatrix}0& 0 & 0 & 0 & 4 & -24\end{bmatrix} \eqcomma
		\\&
		\inbbr{\Gshapev_\porder(0.5)} = \begin{bmatrix} 0 & 0 & 0 & 1 & -1 & 1\end{bmatrix} - \begin{bmatrix} 1 & 1 & 1 & 0 & 0 & 0\end{bmatrix} = \begin{bmatrix}-1 & -1 & -1 & 1 & -1 & 1 \end{bmatrix}  \eqcomma
		\\&
		\iinbbr{\Gshapev_\porder'(0.5)} =  \frac{1}{2}\left(\begin{bmatrix} 0 & 4 & 8 & 0 & 0 & 0\end{bmatrix} + \begin{bmatrix} 0 & 0 & 0 & 4 & -8\end{bmatrix}\right) = \begin{bmatrix}0 & 2 & 4 & 0 & 2 & -4 \end{bmatrix}  \eqdot
	\end{align}
	
	\begin{align}
		&\Gshapev_\bporder(0) = \begin{bmatrix}1 & -1 & 0 & 0 \end{bmatrix},\quad
		\Gshapev_\bporder'(0) = \begin{bmatrix}0& 4 & 0 & 0 \end{bmatrix} \eqcomma
		\\&
		\Gshapev_\bporder(1) = \begin{bmatrix}0 & 0 & 1 & 1 \end{bmatrix},\quad
		\Gshapev_\bporder'(1) = \begin{bmatrix}0& 0 & 0 & 4 \end{bmatrix} \eqcomma
		\\&
		\inbbr{\Gshapev_\bporder(0.5)} = \begin{bmatrix} 0 & 0 & 1 & -1\end{bmatrix} - \begin{bmatrix} 1 & 1 & 0 & 0\end{bmatrix} = \begin{bmatrix}-1 & -1 & 1 & -1 \end{bmatrix}  \eqcomma
		\\&
		\iinbbr{\Gshapev_\bporder'(0.5)} =  \frac{1}{2}\left(\begin{bmatrix} 0 & 4 & 0 & 0\end{bmatrix} + \begin{bmatrix} 0 & 0 & 0 & 4\end{bmatrix}\right) = \begin{bmatrix}0 & 2 & 0 & 2 \end{bmatrix}  \eqdot
	\end{align}
	The ${\tcmtx}$ matrix and $\tcmtvb$ vector in this case are:
	\begin{subequations}\label{eq:1d2p_example_constr}
		\begin{align}\label{eq:1d2p_example_constr_a}
			\begin{aligned}
				\tcmtx =& 	\inbbr{\Gshapev_\bporder(0.5)}\tran\inbbr{\Gshapev_\porder(0.5)} + \Gshapev_\bporder(0)\tran\Gshapev_\porder(0) + \Gshapev_\bporder(1)\tran\Gshapev_\porder(1)
				\\=& \begin{bmatrix}
					1  &    1  & 1  &   -1  &   1 & -1\\
					1   &   1  & 1  &   -1  &   1 & -1\\
					-1  &  -1  & -1 & 1 &  -1  &  1\\
					1  &   1 & 1  &  -1 &  1 & -1
				\end{bmatrix} 
				+
				\begin{bmatrix}
					1 & -1 & 1 & 0 & 0 & 0\\
					-1 & 1 & -1 & 0 & 0 & 0\\
					0 & 0 & 0 & 0 & 0 & 0\\
					0 & 0 & 0 & 0 & 0 & 0
				\end{bmatrix}
				+
				\begin{bmatrix}
					0 & 0 & 0 & 0& 0 & 0\\
					0 & 0 & 0 & 0 & 0 & 0\\
					0 & 0 & 0 & 1& 1 & 1\\
					0 & 0 &0 & 1& 1 & 1
				\end{bmatrix}
				=
				\begin{bmatrix}
					2 & 0 & 2 & -1 & 1 & -1\\
					0 & 2 &0 &  -1 & 1 & -1\\
					-1 &-1 &-1&  2 & 0& 2\\
					1 & 1 & 1&  0 & 2& 0
				\end{bmatrix} 
			\end{aligned} \eqcomma
		\end{align}
		
		\begin{align}\label{eq:1d2p_example_constr_b}
			\tcmtvb =  \Gshapev_\bporder(0)\tran a + \Gshapev_\bporder(1)\tran b = \begin{bmatrix} 2 \\-2\\5\\5 \end{bmatrix} \eqdot
		\end{align}
	\end{subequations}
	
	The solution of the system of equations $(\tcmtx,\, \tcmtvb)$ gives the pair $(\ttransfmatrix,\, \ttransfvectorb)$, which are:
	\begin{align}
		\ttransfmatrix = 
		\begin{bmatrix} 
			-1 & -1 & 0 \\
			0 & -1 & 0\\
			1 & 0 & 0\\
			0 & -1 & -1\\
			0 & 1 & 0 \\
			0 & 0 & 1
		\end{bmatrix} ,
		\quad
		\ttransfvectorb =  \begin{bmatrix} \frac{7}{2}\\\frac{3}{2}\\0\\5\\0\\0\end{bmatrix} \eqdot
	\end{align}
	
	The matrix $\vm{K}$, vector $\vm{f}$ in \eqref{eq:dg_constr_problem}  are defined as:
	\begin{subequations}
		\begin{align}
			\begin{aligned}
				\vm{K} =& \mint_0^1 {\Gshapev_\porder'}\tran\Gshapev_\porder' \intend x  + \inbbr{\Gshapev_\porder(0.5)} \iinbbr{\Gshapev_\porder'(0.5)} 
				+ \Gshapev_\porder(0)\Gshapev_\porder'(0) - \Gshapev_\porder(1)\Gshapev_\porder'(1)
				\\=&
				\begin{bmatrix}
					0 & 0 & 0 & 0 & 0 & 0\\
					0 & 8 & 0 & 0& 0 & 0\\
					0 & 0 & \frac{32}{3} & 0 & 0 & 0\\
					0 & 0 & 0 & 0& 0 & 0\\
					0 & 0 & 0& 0 & 8 & 0\\
					0 & 0 & 0 & 0& 0 & \frac{32}{3} \\
				\end{bmatrix}
				+
				\begin{bmatrix}
					0 & -2 & -4&  0 & -2& 4\\
					0 & -2 & -4&  0 & -2& 4\\
					0 & -2 & -4&  0 & -2& 4\\
					0 &  2 &  4&  0 & -2& 4\\
					0 & -2 & -4&  0 &  2& -4\\
					0 &  2 &  4&  0 &  2& -4\\
				\end{bmatrix}
				+
				\begin{bmatrix}
					0 & -4 & 8 & 0& 0 & 0\\
					0 & 4 & -8 & 0& 0 & 0\\
					0 & -4 & 8 & 0& 0 & 0\\
					0 & 0 & 0 & 0& 0 & 0\\
					0 & 0 & 0 & 0& 0 & 0\\
					0 & 0 & 0 & 0& 0 & 0\\
				\end{bmatrix}
				-
				\begin{bmatrix}
					0 & 0 & 0 & 0& 0 & 0\\
					0 & 0 & 0 & 0& 0 & 0\\
					0 & 0 & 0 & 0& 0 & 0\\
					0 & 0 & 0 & 0 & -4 & -8\\
					0 & 0 & 0 & 0 & -4 & -8\\
					0 & 0 & 0 & 0 & -4 & -8\\
				\end{bmatrix}
				\\=&
				\begin{bmatrix}
					0 & -6 & 4 & 0 & -2 & 4\\
					0 & 10 & -12 & & 0 -2 & 4\\
					0 & -6 & \frac{44}{3} && 0 -2& 4\\
					0 & 2 & 4 & 0 & 6 & 4\\
					0 & -2 & -4 & 0 & 10 & 12\\
					0 & 2 & 4 & 0 & 6 & \frac{44}{3}
				\end{bmatrix} \eqcomma
			\end{aligned} 
		\end{align}
		
		\begin{align}
			\vm{f} = -\mint_0^1 \Gshapev_\porder\tran 2 \intend x = \begin{bmatrix}-1 \\ 0 \\-\frac{1}{3}  \\-1 \\ 0 \\-\frac{1}{3}\end{bmatrix} \eqdot
		\end{align}
	\end{subequations}

	Now, the linear system of equations given in~\eqref{eq:dg_constr_problem} has the following form:
	\begin{align}
		\begin{bmatrix} 
			\frac{32}{3} & 0 & 0\\
			0 & 16 & 0 \\
			0 & 0 & \frac{32}{3}
		\end{bmatrix} \dftempfree = 
		\begin{bmatrix} \frac{2}{3}\\14\\\frac{2}{3} \end{bmatrix}
		\quad\rightarrow\quad \dftempfree = \begin{bmatrix} \frac{2}{32}\\\frac{7}{8}\\\frac{2}{32}\end{bmatrix} \eqdot
	\end{align}
	Finally, on using \eqref{eq:dg_constr} we have the vector $\dftemp$:
	\begin{align}
		\dftemp =  \begin{bmatrix} 
			-1 & -1 & 0 \\
			0 & -1 & 0\\
			1 & 0 & 0\\
			0 & -1 & -1\\
			0 & 1 & 0 \\
			0 & 0 & 1
		\end{bmatrix} \begin{bmatrix} 
			\frac{2}{32}\\\frac{7}{8}\\\frac{2}{32}
		\end{bmatrix} + 
		\begin{bmatrix} 
			\frac{7}{2}\\\frac{3}{2}\\0\\5\\0\\0\end{bmatrix} = \begin{bmatrix} \frac{41}{16}\\\frac{5}{8}\\\frac{1}{16}\\\frac{65}{16} \\ \frac{7}{8} \\ \frac{1}{16}\end{bmatrix} \eqdot
	\end{align}
	On applying the vector $\dftemp$ to the approximation in \eqref{eq:ex_1d2p_approx} we recover the exact solution (quadratic function) of the problem.

	\section{Solving the constraint equation} \label{sec:1d_constr_example}
	In this section, we present the step-by-step solution of the constraint problem in~\eqref{eq:1d_example_constr}. The problem can be solved using the procedure in Appendix~\ref{sec:dg_solver}, but here we use the procedure described in Section~\ref{sec:constriants}.
	
	The constraints for the vector $\dftemp$ are expressed as: 
	\begin{align}
		\dftemp = \transfmatrix\freevar + \transfvector \eqcomma
	\end{align}
	where at the beginning of the process $\transfmatrix=\vm I$,  $\transfvector=\vm 0$ and $\freevar=\dftemp$.
	
	According to~\eqref{eq:segm_jump_eq}, we construct the constraint equation to enforce the continuity of the approximation on the mesh skeleton, which is at $x=\frac{1}{2}$ herein:
	\begin{align}
		\segmentmatrix^{\frac{1}{2}} \dftemp = \vm 0 \eqcomma
	\end{align} 
	where
	\begin{align}
		\segmentmatrix^{\frac{1}{2}} = \inbbr{\Gshapev(1/2)}\tran\inbbr{\Gshapev(1/2)}
		= \begin{bmatrix}
			1 &    1 &   -1  &   1\\
			1  &   1 &   -1  &   1\\
			-1  &  -1 &    1  &  -1\\
			1  &   1  &  -1   &  1
		\end{bmatrix} \eqdot
	\end{align}
	It is evident that all four rows in the matrix $\segmentmatrix^{\frac{1}{2}}$ are linearly dependent. When the dependent rows are removed the matrix $\segmentmatrixred^{\frac{1}{2}}$ is obtained: 
	\begin{align}\label{eq:ex1d_sol_reduced}
		\segmentmatrixred^{\frac{1}{2}} \dftemp = \vm 0,
		\quad \text{where} \quad 
		\segmentmatrixred^{\frac{1}{2}} = \begin{bmatrix}1 &    1 &   -1  &   1\end{bmatrix} \eqdot
	\end{align}
	When solving \eqref{eq:ex1d_sol_reduced} the pivot has to be chosen. It can be any component of the vector $\dftemp$, and here is chosen as the first one, i.e.,
	\begin{align}
		\dftempsc_1 = \begin{bmatrix}-1 & 1 & -1 \end{bmatrix}
		\begin{bmatrix}\dftempsc_2 \\ \dftempsc_3 \\ \dftempsc_4 \end{bmatrix} \eqdot
	\end{align}
	The vector $\dftemp$ is now expressed as
	\begin{align} \label{eq:1d_skel_sol}
		\dftemp = \transfmatrix\freevar + \transfvector \eqcomma
	\end{align}
	where
	\begin{align}
		\transfmatrix = \begin{bmatrix}
			-1 &    1   & -1\\
			1  &   0  &   0\\
			0   &  1   &  0\\
			0   &  0   &  1
		\end{bmatrix},
		\quad
		\transfvector = \begin{bmatrix}	0\\0\\0\\0\end{bmatrix},
		\quad
		\freevar = \begin{bmatrix}\dftempsc_2 \\ \dftempsc_3 \\ \dftempsc_4 \end{bmatrix} \eqdot
	\end{align}
	
	Now we consider the boundary constraints. There are two parts of the boundary in this example, on the left and right sides of the domain. The boundary conditions on these two ends have to be met independently, \eqref{eq:1d_poisson_b}. The condition on the left side  is associated with the first finite element, so we evaluate the approximation on this element for $x=0$, and require that it has the value of $2$:
	\begin{align}
		{\Gshapev^1(0)}\tran \Gshapev^1(0) \dftemp^1 = 2 {\Gshapev^1(0)}\tran \eqcomma
	\end{align}
	\begin{align}
		\begin{bmatrix}	1 & -1 \\ -1 & 1 \end{bmatrix}
		\begin{bmatrix} \dftempsc_1 \\ \dftempsc_2 \end{bmatrix} = \begin{bmatrix}2 \\ -2 \end{bmatrix}
		\eqdot
	\end{align}
	This set of equations consists of two linearly dependent equations, so it can be reduced to a single equation:
	\begin{align}\label{eq:1d_bd1_red}
		\begin{bmatrix}	1 & -1 \end{bmatrix}
		\begin{bmatrix} \dftempsc_1 \\ \dftempsc_2 \end{bmatrix} = \begin{bmatrix}2 \end{bmatrix}
		\eqdot
	\end{align}
	The vector $\dftemp^1$ has to satisfy the constraints in \eqref{eq:1d_skel_sol}, so \eqref{eq:1d_bd1_red} is changed to
	\begin{align}\label{eq:1d_bd1_red2}
		\begin{bmatrix}	1 & -1 \end{bmatrix}
		\begin{bmatrix}
			-1 &    1   & -1\\
			1  &   0  &   0\\
		\end{bmatrix}\freevar = \begin{bmatrix} 2\end{bmatrix}
		\quad \rightarrow \quad
		\begin{bmatrix}-2 & 1 & -1 \end{bmatrix} \freevar = \begin{bmatrix}2\end{bmatrix} 
		\eqdot
	\end{align}
	The solution of \eqref{eq:1d_bd1_red2}, presented in recurrence form, is:
	\begin{align}
		\freevar = \begin{bmatrix}0.5 & -0.5 \end{bmatrix} \freevar + \begin{bmatrix}-1	\end{bmatrix}
		\eqdot
	\end{align}
	When this equation is applied to~\eqref{eq:1d_skel_sol}, we have
	\begin{align}
		\dftemp = \begin{bmatrix}
			0.5 & -0.5 \\ 0.5 & -0.5 \\ 1 & 0 \\ 0 & 1
		\end{bmatrix} \freevar
		+ \begin{bmatrix}	1 \\ -1 \\ 0 \\ 0 \end{bmatrix} \eqdot
	\end{align}
	The same procedure is performed on the second boundary:
	\begin{align}
		{\Gshapev^2(0)}\tran \Gshapev^2(0) \dftemp^2 = 5 {\Gshapev^1(0)}\tran \eqcomma
	\end{align}
	\begin{align}\label{eq:1d_bd2_red}
		&\begin{bmatrix}	1 & 1 \\ 1 & 1 \end{bmatrix}
		\begin{bmatrix} \dftempsc_3 \\ \dftempsc_4 \end{bmatrix} = \begin{bmatrix}5 \\ 5 \end{bmatrix}
		\quad \rightarrow \quad
		\begin{bmatrix}	1 & 1 \end{bmatrix}
		\begin{bmatrix} \dftempsc_3 \\ \dftempsc_4 \end{bmatrix} = \begin{bmatrix}5  \end{bmatrix}
		\quad\rightarrow\quad
		\begin{bmatrix}	1 & 1 \end{bmatrix}	\freevar = \begin{bmatrix}5  \end{bmatrix}
		\\&
		\quad\rightarrow\quad
		\freevar = \begin{bmatrix}-1\end{bmatrix}\freevar + \begin{bmatrix}5 \end{bmatrix}
		\quad\rightarrow\quad
		\dftemp = \begin{bmatrix} -1\\-1\\-1\\1\end{bmatrix}\freevar + \begin{bmatrix} 3.5\\1.5\\5\\0\end{bmatrix}
		\eqdot
	\end{align}

\bibliographystyle{plain}
\bibliography{references}

\end{document}